\documentclass[]{article}
\usepackage{lmodern}
\usepackage{amssymb,amsmath}
\usepackage{ifxetex,ifluatex}
\usepackage{fixltx2e} 
\ifnum 0\ifxetex 1\fi\ifluatex 1\fi=0 
  \usepackage[T1]{fontenc}
  \usepackage[utf8]{inputenc}
\else 
  \ifxetex
    \usepackage{mathspec}
  \else
    \usepackage{fontspec}
  \fi
  \defaultfontfeatures{Ligatures=TeX,Scale=MatchLowercase}
\fi
\IfFileExists{upquote.sty}{\usepackage{upquote}}{}
\IfFileExists{microtype.sty}{%
\usepackage{microtype}
\UseMicrotypeSet[protrusion]{basicmath} 
}{}
\usepackage[margin=1in]{geometry}
\usepackage{hyperref}
\hypersetup{unicode=true,
            pdftitle={Statistical dependence: Beyond Pearson's \textbackslash{}rho},
            pdfauthor={Dag Tjøstheim, Håkon Otneim, Bård Støve\^{}\{\textbackslash{}dagger\}},
            pdfborder={0 0 0},
            breaklinks=true}
\usepackage{longtable,booktabs}
\usepackage{graphicx,grffile}
\makeatletter
\def\maxwidth{\ifdim\Gin@nat@width>\linewidth\linewidth\else\Gin@nat@width\fi}
\def\maxheight{\ifdim\Gin@nat@height>\textheight\textheight\else\Gin@nat@height\fi}
\makeatother
\setkeys{Gin}{width=\maxwidth,height=\maxheight,keepaspectratio}
\IfFileExists{parskip.sty}{%
\usepackage{parskip}
}{
\setlength{\parindent}{0pt}
\setlength{\parskip}{6pt plus 2pt minus 1pt}
}
\providecommand{\tightlist}{%
  \setlength{\itemsep}{0pt}\setlength{\parskip}{0pt}}
\ifx\paragraph\undefined\else
\let\oldparagraph\paragraph
\renewcommand{\paragraph}[1]{\oldparagraph{#1}\mbox{}}
\fi
\ifx\subparagraph\undefined\else
\let\oldsubparagraph\subparagraph
\renewcommand{\subparagraph}[1]{\oldsubparagraph{#1}\mbox{}}
\fi

\let\rmarkdownfootnote\footnote%
\def\footnote{\protect\rmarkdownfootnote}

\usepackage{titling}

  \title{Statistical dependence: Beyond Pearson's
\(\rho\)\footnote{This work has been partly supported by the Finance Market Fund, grant number 261570}}
    \pretitle{\vspace{\droptitle}\centering\huge}
  \posttitle{\par}
    \author{Dag
Tjøstheim\footnote{Department of Mathematics, University of Bergen, P.b. 7803, 5020 Bergen, Norway},
Håkon
Otneim\footnote{Department of Business and Management Science, Norwegian School of Economics, Helleveien 30, 5045 Bergen, Norway. Corresponding author, e-mail: \texttt{hakon.otneim@nhh.no}},
Bård Støve\(^{\dagger}\)}
    \preauthor{\centering\large\emph}
  \postauthor{\par}
      \predate{\centering\large\emph}
  \postdate{\par}
    \date{\today}

\newcommand{\Cov}{\textrm{Cov}}
\newcommand{\E}{\textrm{E}}
\newcommand{\Var}{\textrm{Var}} 
\newcommand{\di}{\,\textrm{d}}
\newcommand{\Corr}{\textrm{Corr}} 

\usepackage{graphicx}
\usepackage{subfig}
\usepackage{bm}

\begin{document}
\maketitle
\begin{abstract}
Pearson's \(\rho\) is the most used measure of statistical dependence.
It gives a complete characterization of dependence in the Gaussian case,
and it also works well in some non-Gaussian situations. It is well
known, however, that it has a number of shortcomings; in particular for
heavy tailed distributions and in nonlinear situations, where it may
produce misleading, and even disastrous results. In recent years a
number of alternatives have been proposed. In this paper, we will survey
these developments, especially results obtained in the last couple of
decades. Among measures discussed are the copula, distribution-based
measures, the distance covariance, the HSIC measure popular in machine
learning, and finally the local Gaussian correlation, which is a local
version of Pearson's \(\rho\). Throughout we put the emphasis on
conceptual developments and a comparison of these. We point out relevant
references to technical details as well as comparative empirical and
simulated experiments. There is a broad selection of references under
each topic treated.
\end{abstract}

{
\setcounter{tocdepth}{2}
\tableofcontents
}
\section{Introduction}\label{intro}

Pearson's \(\rho\), the product moment correlation, was not invented by
Pearson, but rather by Francis Galton. Galton, a cousin of Charles
Darwin, needed a measure of association in his hereditary studies,
Galton (1888; 1890). This was formulated in a scatter diagram and
regression context, and he chose \(r\) (for regression) as the symbol
for his measure of association. Pearson (1896) gave a more precise
mathematical development and used \(\rho\) as a symbol for the
population value and \(r\) for its estimated value. The product moment
correlation is now universally referred to as Pearson's \(\rho\). Galton
died in 1911, and Karl Pearson became his biographer, resulting in a
massive 4-volume biography, Pearson (1922; 1930). All of this and much
more is detailed in Stigler (1989) and Stanton (2001). Some other
relevant historical references are Fisher (1915; 1921), von Neumann
(1941; 1942) and the survey paper by King (1987).

Write the covariance between two random variables \(X\) and \(Y\) having
finite second moments as
\(\Cov(X,Y) = \sigma(X,Y) = \E(X-\E(X))(Y-\E(Y))\). The Pearson's
\(\rho\), or the product moment correlation, is defined by
\[\rho = \rho(X,Y)=\frac{\sigma(X,Y)}{\sigma_X \sigma_Y}\] with
\(\sigma_{X} = \sqrt{\sigma_{X}^2} = \sqrt{\E(X-\E(X))^2}\) being the
standard deviation of \(X\) and similarly for \(\sigma_{Y}\). The
correlation takes values between and including \(-1\) and \(+1\). For a
given set of pairs of observations \((X_1,Y_1),\ldots,(X_n,Y_n)\) of
\(X\) and \(Y\), an estimate of \(\rho\) is given by

\begin{equation}
r = \widehat{\rho} = \frac{\sum_{j=1}^n (X_j-\overline{X})(Y_j-\overline{Y})}{\sqrt{\sum_{j=1}^{n}(X_j-\overline{X})^2}\sqrt{\sum_{j=1}^{n}(Y_j-\overline{Y})^2}} 
\label{eq:estimated}
\end{equation}

with \(\overline{X} = n^{-1}\sum_{j=1}^n X_j\), and similarly for
\(\overline{Y}\). Consistency and asymptotic normality can be proved
using an appropriate law of large numbers and a central limit theorem,
respectively.

The correlation coefficient \(\rho\) has been, and probably still is,
\emph{the} most used measure for statistical association. There are
several reasons for this.

\begin{enumerate}
\def\labelenumi{(\roman{enumi})}
\item
  It is easy to compute (estimate), and it is generally accepted as
  \emph{the} measure of dependence, not only in statistics, but in most
  applications of statistics to the natural and social sciences.
\item
  Linear models are much used, and in a linear regression model of \(Y\)
  on \(X\), say, \(\rho\) is proportional to the slope of the regression
  line. More precisely; if \(Y_i = \alpha+\beta X_i+\varepsilon_i,\)
  where \(\{\varepsilon_i\}\) is a sequence of zero-mean iid error terms
  whose second moment exists, then
  \[\beta = \rho(X,Y)\frac{\sigma_Y}{\sigma_X}.\] This also means that
  \(\rho\) and its estimate \(\widehat{\rho}\) appears naturally in a
  linear least squares analysis
\item
  In a bivariate Gaussian density

  \begin{align*}
  &f(x,y) = \frac{1}{2\pi \sqrt{1-\rho^2}\sigma_{X} \sigma_{Y}} \\
  & \times \exp \left\{-\frac{1}{2(1-\rho^2)}\left(\frac{(x-\mu_{X})^2}{\sigma_{X}^2}-2\rho\frac{(x-\mu_{X})(y-\mu_{Y})}{\sigma_{X} \sigma_{Y}}+\frac{(y-\mu_{Y})^2}{\sigma_{Y}^2}\right)\right\},
  \end{align*}

  the dependence between \(X\) and \(Y\) is completely characterized by
  \(\rho\). In particular, two jointly Gaussian variables \((X,Y)\) are
  independent if and only if they are uncorrelated (See e.g.~Billingsley
  (2008, 384--85) for a formal proof of this statement). For a
  considerable number of data sets, the Gaussian distribution works at
  least as a fairly good approximation. Moreover, joint asymptotic
  normality often appears as a consequence of the central limit theorem
  for many statistics, and the joint asymptotic behavior of such
  statistics are therefore generally well defined by the correlation
  coefficient.
\item
  The product moment correlation is easily generalized to the
  multivariate case. For \(p\) stochastic variables \(X_1,\ldots,X_p\),
  their joint dependencies can simply (but not always accurately) be
  characterized by their covariance matrix \(\Sigma=\{\sigma_{ij}\}\),
  with \(\sigma_{ij}\) being the covariance between \(X_i\) and \(X_j\).
  Similarly the correlation matrix is defined by
  \(\Lambda =\{\rho_{ij}\}\), with \(\rho_{ij}\) being the correlation
  between \(X_i\) and \(X_j\). Again, for a column vector
  \(x = (x_1,\ldots,x_p)^T\), the joint normality density is defined by
  \[f(x) = \frac{1}{(2\pi)^{p/2}|\Sigma|^{1/2}}\exp \left\{-\frac{1}{2}(x-\mu)^{T}\Sigma^{-1}(x-\mu)\right\}\]
  where \(|\Sigma|\) is the determinant of the covariance matrix
  \(\Sigma\) (whose inverse \(\Sigma^{-1}\) is assumed to exist), and
  \(\mu=\E(X)\). Then the complete dependence structure of the Gaussian
  vector is given by the \emph{pairwise} covariances \(\sigma_{ij}\), or
  equivalently the \emph{pairwise} correlations \(\rho_{ij}\). This is
  remarkable: the entire dependence structure is determined by pairwise
  dependencies. We will make good use of this fact later when we get to
  the local Gaussian dependence measure in Section \ref{lgc}.
\item
  It is easy to extend the correlation concept to time series. For a
  time series \(\{X_t\}\), the autocovariance and autocorrelation
  function, respectively, are defined, assuming stationarity and
  existence of second moments, by \(c(t) = \sigma(X_{t+s},X_s)\) and
  \(\rho(t) = \rho(X_{t+s},X_s)\) for arbitrary integers \(s\) and
  \(t\). For a Gaussian time series, the dependence structure is
  completely determined by \(\rho(t)\). For linear (say ARMA) type
  series the analysis as a rule is based on the autocovariance function,
  even though the entire joint probability structure cannot be captured
  by this in the non-Gaussian case. Even for nonlinear time series and
  nonlinear regression models the autocovariance function has often been
  made to play a major role. In the frequency domain all of the
  traditional spectral analysis is based again on the autocovariance
  function. Similar considerations have been made in spatial models such
  as in linear Kriging models, see Stein (1999).
\end{enumerate}

In spite of these assets, there are several serious weaknesses of
Pearson's \(\rho\). These will be briefly reviewed in Section
\ref{weaknesses}. In the remaining sections of this paper a number of
alternative dependence measures going beyond the Pearson \(\rho\) will
be described. The emphasis will be on concepts, conceptual developments
and comparisons of these. We do provide some illustrative plots of key
properties, but when it comes to technical details, empirical and
simulated experiments with numerical comparisons, we point out relevant
references instead.

\section{\texorpdfstring{Weaknesses of Pearson's
\(\rho\)}{Weaknesses of Pearson's \textbackslash{}rho}}\label{weaknesses}

We have subsumed, somewhat arbitrarily, the problems of Pearson's
\(\rho\) under three issues:

\subsection{The non-Gaussianity issue}\label{nongauss}

A natural question to ask is whether the close connection between
Gaussianity and correlation/covariance properties can be extended to
larger classes of distributions. The answer to this question is a
conditional yes. The multivariate Gaussian distribution is a member of
the vastly larger class of elliptical distributions. That class of
distributions is defined both for discrete and continuous variables, but
we limit ourselves to the continuous case. An elliptical distribution
can be defined in terms of a parametric representation of the
characteristic function or the density function. For our purposes it is
simplest to phrase this in terms of a density function.

Consider a stochastic vector \(X = (X_1,\ldots,X_p)\) and a non-negative
Lebesgue measurable function \(g\) on \([0,\infty)\) such that
\[ \int_0^{\infty}x^{\frac{p}{2}-1}g(x)\, \textrm{d}x < \infty.\]
Further, let \(\mu \in \mathbb{R}^{p}\) and let \(\Sigma\) be a positive
definite \(p \times p\) matrix, then an elliptical density function
parameterized by \(g\), \(\mu\) and \(\Sigma\) is given by

\begin{equation}
f(x; \mu, \Sigma, g) = c_p|\Sigma|^{-1/2}g\left((x-\mu)^{T}\Sigma^{-1}(x-\mu)\right),
\label{eq:elliptic1}
\end{equation}

where \(c_p\) is a normalizing factor given by
\[c_p = \frac{\Gamma(p/2)}{(2\pi)^{p/2}}\left(\int_0^{\infty} x^{\frac{p}{2}-1}g_p(x)\,\textrm{d}x\right)^{-1}.\]
The parameters \(\mu\) and \(\Sigma\) can be interpreted as location and
scale parameters, respectively, but they cannot in general be identified
with the mean \(\E(X)\) and covariance matrix \(\Cov(X)\). In fact the
parameters \(\mu\) and \(\Sigma\) in equation \eqref{eq:elliptic1} may
remain meaningful even if the mean and the covariance matrix do not
exist. If they do exist, \(\mu\) can be identified with the mean, and
\(\Sigma\) is proportional to the covariance matrix, the proportionality
factor in general depending on \(p\). A redefinition of \(c_p\) may then
make this proportionality factor equal to 1, cf. Gómez, Gómez-Villegas,
and Mari'in (2003) and Landsman and Valdez (2003).

A number of additional properties of elliptical distributions, among
other things pertaining to linear transformations, marginal
distributions and conditional distributions are surveyed in Gómez,
Gómez-Villegas, and Mari'in (2003) and Landsman and Valdez (2003). Many
of these properties are analogous to those of the multivariate normal
distribution, which is an elliptical distribution defined by
\(g(x) = \exp\{-x/2\}\).

Unfortunately, the equivalence between uncorrelatedness and independence
is generally not true for elliptical distributions. Consider for
instance the multivariate \(t\)-distribution with \(\nu\) degrees of
freedom

\begin{equation}
f(x) = 
\frac{\Gamma(\frac{p+\nu}{2})}{(\pi\nu)^{p/2}\Gamma(\nu/2)|\Sigma|^{1/2}}\left(1+\frac{(x-\mu)^{T}\Sigma^{-1}(x-\mu)}{\nu}\right)^{-\frac{p+\nu}{2}}.  
\label{eq:elliptic2}
\end{equation}

Unlike the multinormal distribution where the exponential form of the
distribution forces the distribution to factor if \(\Sigma\) is a
diagonal matrix (uncorrelatedness), this is not true for the \(t\)
distribution defined in equation \eqref{eq:elliptic2} if \(\Sigma\) is
diagonal. In other words, if two components of a bivariate \(t\)
distribution are uncorrelated, they are not necessarily independent.
This pinpoints a serious deficiency of the Pearson's \(\rho\) in
measuring dependence in \(t\) distributions, and indeed in general
elliptical (and of course non-elliptical) distributions.

\subsection{The robustness issue}\label{robustness}

As is the case for regression, it is well known that the product moment
estimator is sensitive to outliers. Even just one single outlier may be
very damaging. There are therefore several robustified versions of
\(\rho\), primarily based on ranks. The idea of rank correlation goes
back at least to Spearman (1904), and it is most easily explained
through its sample version. Given scalar observations
\(\{X_1,\ldots,X_n\}\), we denote by \(R_{i,X}^{(n)}\) the rank of
\(X_i\) among \(X_1,\ldots,X_n\). (There are various rules for treating
ties). The estimated Spearman rank correlation function given \(n\)
pairwise observations of two random variables \(X\) and \(Y\) is given
by
\[\widehat{\rho}_S = \frac{n^{-1}\sum_{i=1}^{n} R_{i,X}^{(n)}R_{i,Y}^{(n)}-(n+1)^2/4}{(n^2 - 1)/12}.\]
If \(X\) and \(Y\) have continuous cumulative distribution functions
\(F_{X}\) and \(F_{Y}\), and joint distribution function \(F_{X,Y}\),
then the population value of the Spearman's \(\rho_S\) is given by

\begin{equation}
\rho_S = 12 \int F_{X}(x)F_{Y}(y)\,\textrm{d}F_{X,Y}(x,y) - 3,
\label{eq:spearman}
\end{equation}

and hence it is a linear transformation of the correlation between the
two uniform variables \(F_{X}(X)\) and \(F_{Y}(Y)\). The rank
correlation is thought to be especially effective in picking up linear
trends in the data, but it suffers in a very similar way as the
Pearson's \(\rho\) to certain nonlinearities of the data which are
treated in the next subsection. Spearman's \(\rho\) may be modified to a
rank autocorrelation measure for time series in the obvious way, see
Knoke (1977), Bartels (1982), Hallin and Mélard (1988), and Ferguson,
Genest, and Hallin (2000).

Another way of using the ranks is the Kendall's \(\tau\) rank
correlation coefficient given by Kendall (1938). Again, consider the
situation of \(n\) pairs \((X_i,Y_i)\) of the random variables \(X\) and
\(Y\). Two pairs of observations \((X_{i},Y_{i})\) and
\((X_{j},Y_{j})\), \(i \neq j\) are said to be concordant if the ranks
for both elements agree; that is, if both \(X_{i} > X_{j}\) and
\(Y_{i} > Y_{j}\) or if both \(X_{i} < X_{j}\) and \(Y_{i} < Y_{j}\).
Similarly, they are said to be discordant if \(X_{i} > X_{j}\) and
\(Y_{i} < Y_{j}\) or if \(X_{i} < X_{j}\) and \(Y_{i} > Y_{j}\). If one
has equality, they are neither concordant nor discordant, even though
there are various rules for treating ties in this case as well. The
estimated Kendall \(\tau\) is then given by
\[ \widehat{\tau} = \frac{\textrm{(number of concordant pairs)} - (\textrm{number of discordant pairs})}{n(n-1)/2}. \]
The population value can be shown to be

\begin{equation}
\tau = 4\int F_{X,Y}(x,y)\,\textrm{d}F_{X,Y}(x,y)-1.
\label{eq:kendall}
\end{equation}

Both \(\rho_s\) and \(\tau\) are expressible in terms of the copula (see
Section \ref{copula}) associated with \(F_{X,Y}\). It is then perhaps
not surprising that both \(\rho_s\) and \(\tau\) are bivariate measures
of monotone dependence. This means that (i) they are invariant with
respect to strictly increasing (decreasing) transformations of both
variables, and (ii) they are equal to \(1\) (or \(-1\)) if one of the
variables is an increasing (or decreasing) transformation of the other
one. Property (i) does not hold for Pearson's \(\rho\), and \(\rho\) is
not directly expressible in terms of the copula of \(F_{X,Y}\) either.
The invariance property (i) is also shared by the van der Waerden (1952)
correlation based on normal scores. Some will argue that this invariance
property make them more desirable as dependence measures in case \(X\)
and \(Y\) are non-Gaussian.

The asymptotic normality of the Spearman's \(\rho_S\) and Kendall's
\(\tau\) was established early. Some of the theory is reviewed in
Kendall (1970). It can be viewed as special cases of much more general
results obtained by Hallin, Ingenbleek, and Puri (1985) and Ferguson,
Genest, and Hallin (2000). For some details in the time series case we
refer to Tjøstheim (1996). A more recent account from the copula point
of view is given by Genest and Rémillard (2004). They show in Section 3
of their paper that serialized and non-serialized versions of Spearman's
\(\rho\) and other linear rank statistics share the same limiting
distribution.

\begin{figure}[t]
\subfloat[Gaussian\label{fig:issues1}]{\includegraphics[width=0.24\linewidth]{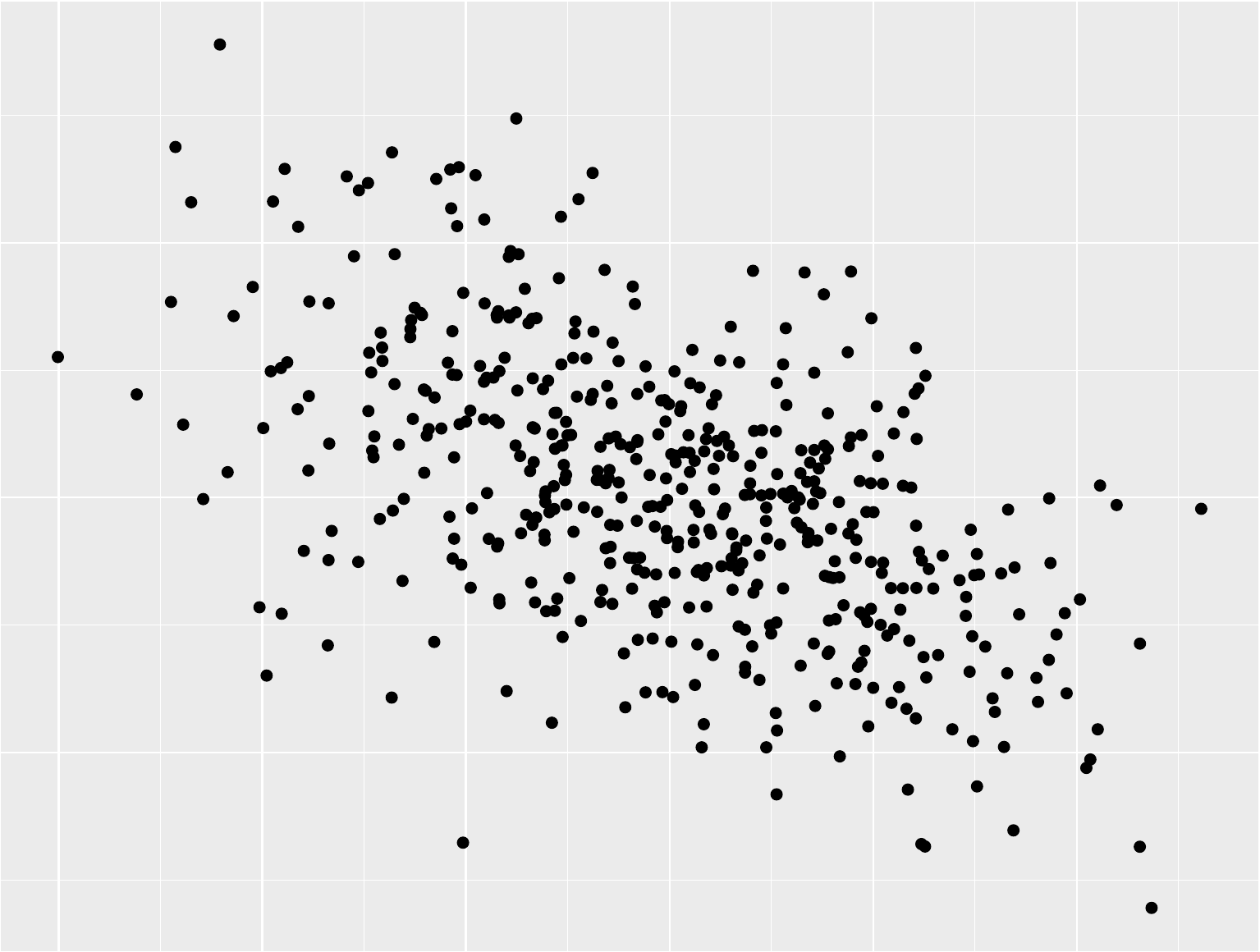} }\subfloat[Gaussian with outliers\label{fig:issues2}]{\includegraphics[width=0.24\linewidth]{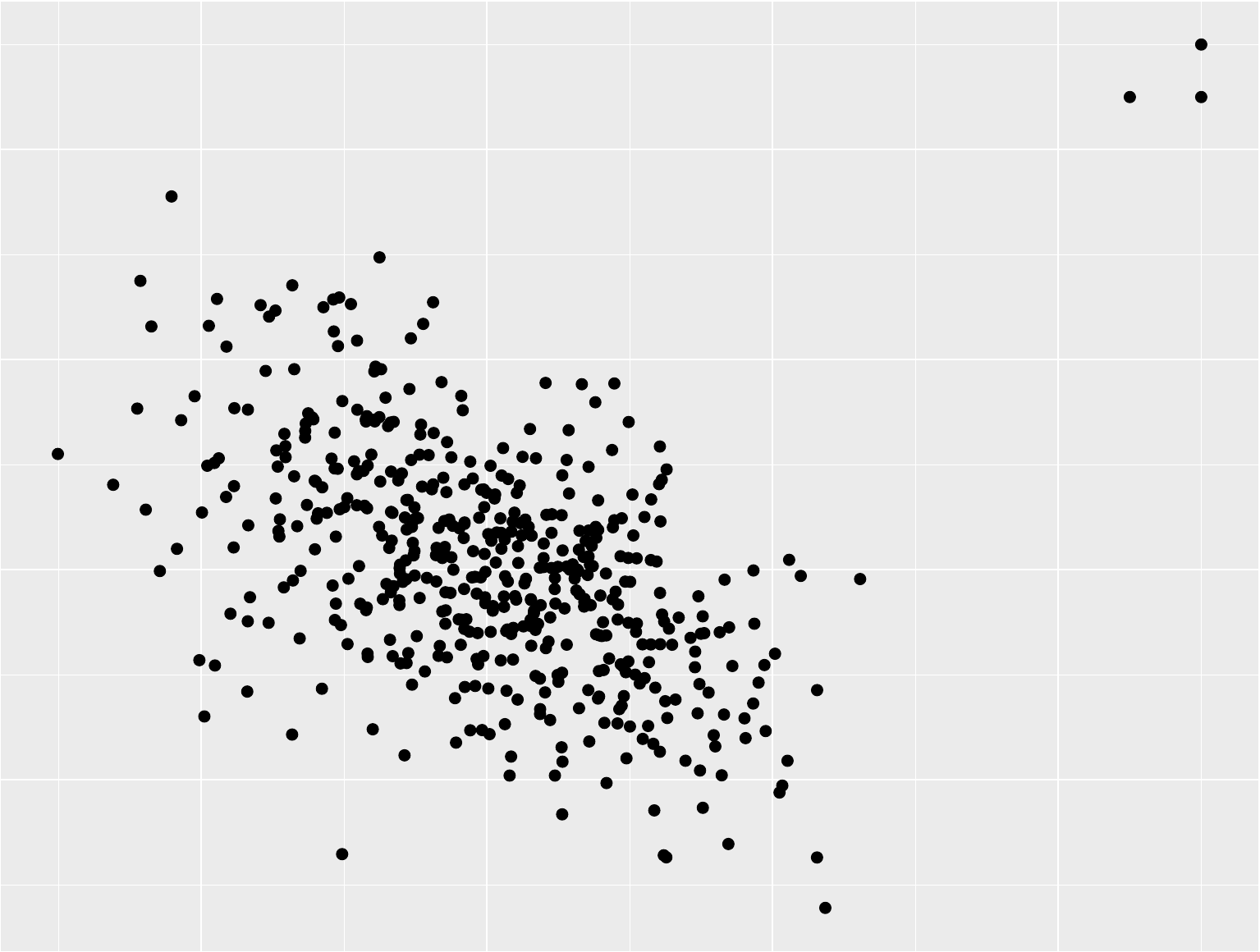} }\subfloat[Non-linear\label{fig:issues3}]{\includegraphics[width=0.24\linewidth]{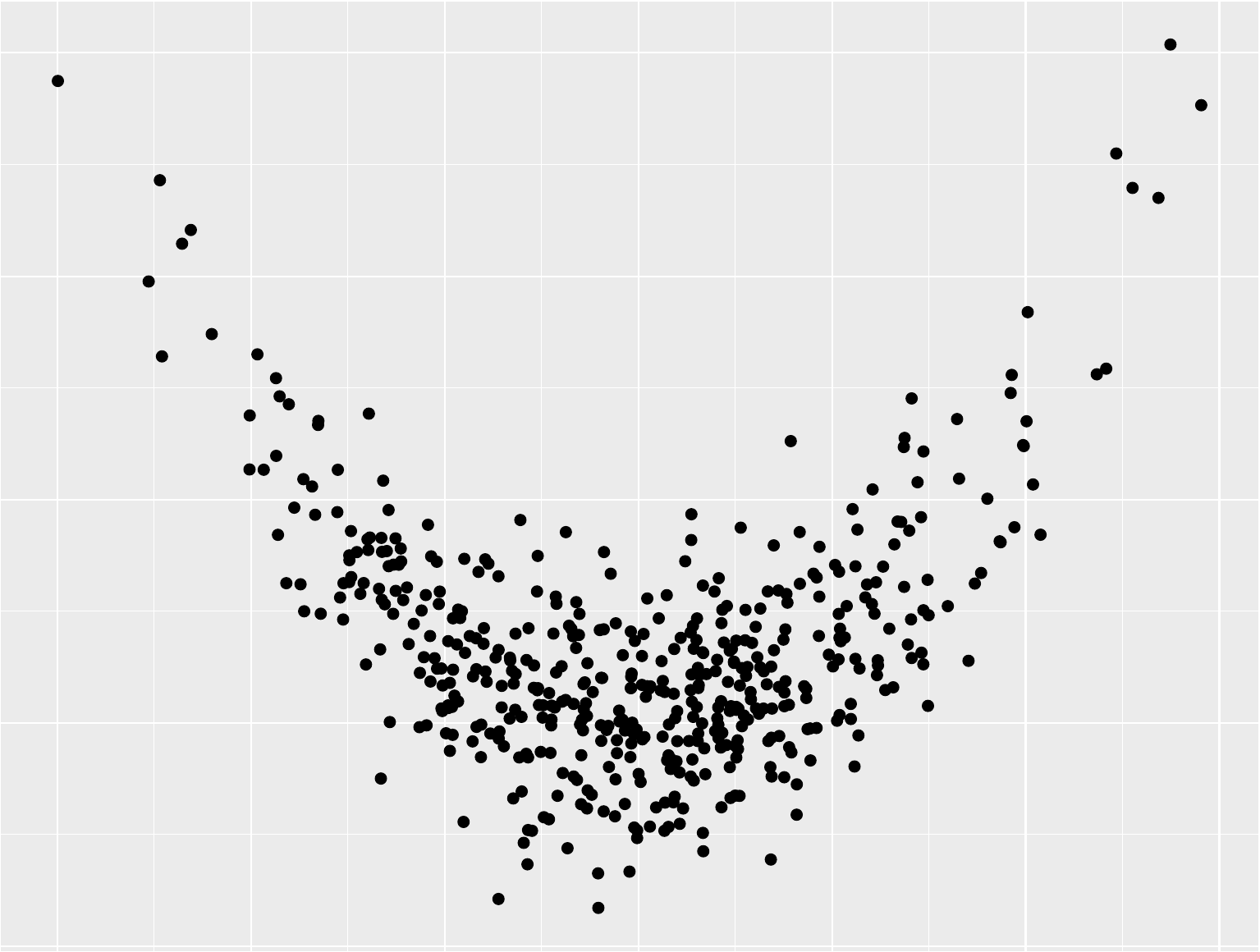} }\subfloat[Garch\label{fig:issues4}]{\includegraphics[width=0.24\linewidth]{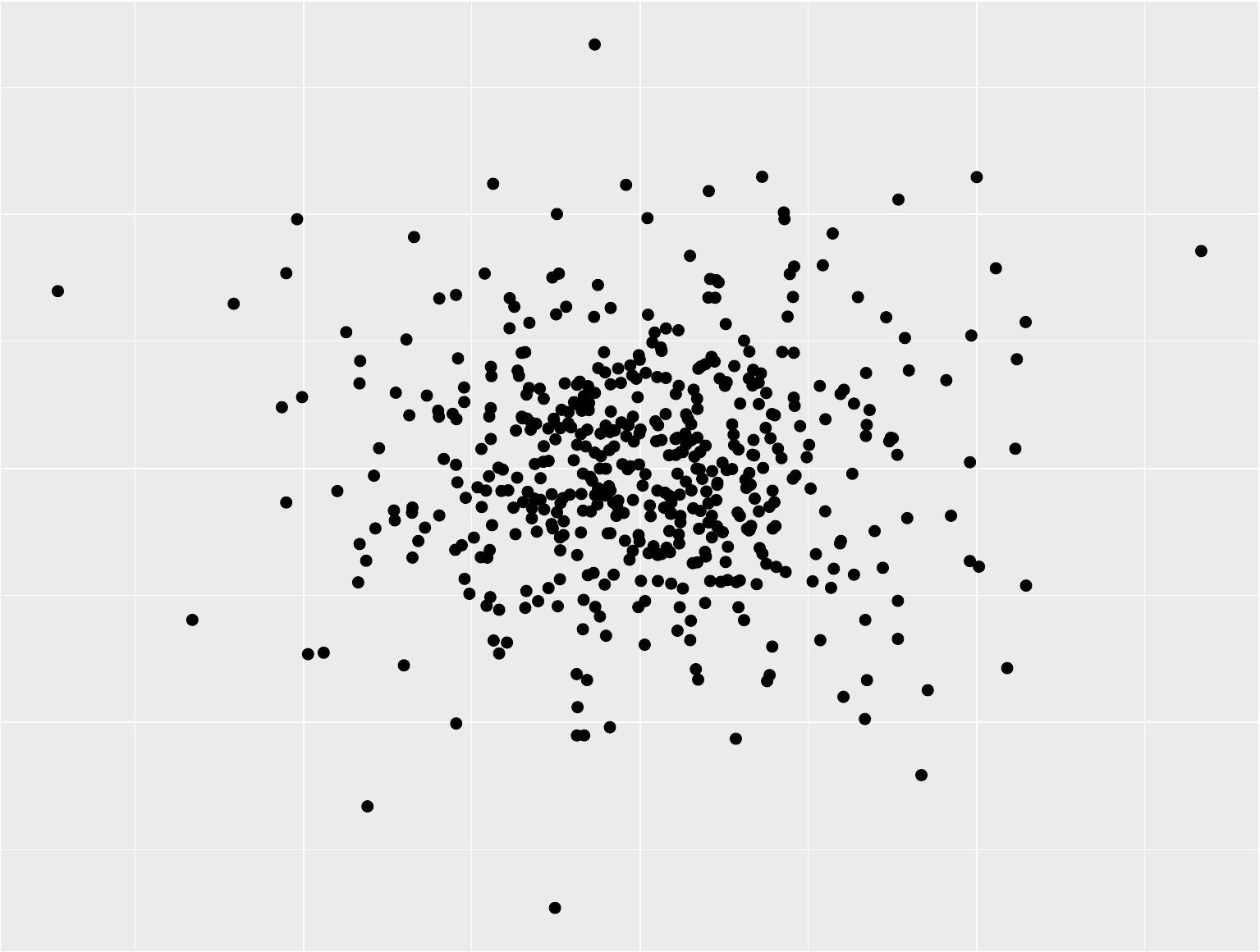} }\caption{Illustration of some problems related to the Pearson correlation coefficient}\label{fig:issues}
\end{figure}

We will illustrate the robustness issue using a simple example. In
Figure \ref{fig:issues1} we see 500 observations that have been
simulated from the bivariate Gaussian distribution having correlation
\(\rho = -0.5\). The sample value for Pearson's \(\rho\) is
\(\widehat \rho = -0.53\). If we add just three outliers to the data,
however, as shown in Figure \ref{fig:issues2}, the sample correlation
changes to \(\widehat \rho = -0.36\). The sample versions of Spearman's
\(\rho\) for the simulated data in Figures \ref{fig:issues1} and
\ref{fig:issues2} are on the other hand very similar:
\(\widehat \rho_S = -0.52\) and \(\widehat \rho_S = -0.49\), and the
corresponding values for the estimated Kendall's \(\tau\) are
\(\widehat \tau = -0.37\) and \(\widehat \tau = -0.35\).

\subsection{The nonlinearity issue}\label{nonlinearity}

This is probably the most serious issue with Pearson's \(\rho\), and it
is an issue also for the rank based correlations of Spearman and
Kendall. All of these (and similar measures), are designed to detect
rather specific types of statistical dependencies, namely those for
which large values of \(X\) tend to be associated with large values of
\(Y\), and small values of \(X\) with small values of \(Y\) (positive
dependence), or the opposite case of negative dependence in which large
values of one variable tend to be associated with small values of the
other variable. It is easy to find examples where this is not the case,
but where nevertheless there is strong dependence. A standard
introductory text book example is the case where

\begin{equation}
Y=X^2.
\label{eq:parabola}
\end{equation}

Here, \(Y\) is uniquely determined once \(X\) is given; i.e., basically
the strongest form of dependence one can have. If \(\E(X)=\E(X^3)=0\),
however, it is trivial to show that \(\rho(X,Y) = 0\), and moreover that
\(\rho_s\) and \(\tau\) will also fail spectacularly. A version of this
situation is illustrated in Figure \ref{fig:issues3}, where we have
generated 500 observations of the standard normally distributed
independent variables \(X\) and \(\epsilon\), and calculated \(Y\) as
\(Y = X^2 + \epsilon\). Still, \(\rho(X,Y) = 0\). The sample values for
Pearson's \(\rho\), Spearman's \(\rho_S\) and Kendall's \(\tau\) are
\(\widehat\rho = -0.001\), \(\widehat\rho_S = -0.03\) and
\(\widehat\tau = -0.02\) respectively, and none of them are
significantly different from zero.

Essentially the same problem will occur if \(X= UW\) and \(Y=VW\), where
\(U\) and \(V\) are independent of each other and independent of \(W\).
It is trivial to show that \(\rho(X,Y)=0\) if \(\E(U) = \E(V)=0\),
whereas \(X\) and \(Y\) are clearly dependent. This example typifies the
kind of dependence one has in ARCH/GARCH time series models: If
\(\{\varepsilon_t\}\) is a time series of zero-mean iid variables and if
the time series \(\{h_t\}\) is independent of \(\{\varepsilon_t\}\), and
\(\{X_t\}\) and \(\{h_t\}\) are given by the recursive relationship

\begin{equation}
X_t = \varepsilon_t h_t^{1/2}, \quad \quad h_t = \alpha + \beta h_{t-1} + \gamma X_{t-1}^2,
\label{eq:garch}
\end{equation}

where the stochastic process \(\{h_t\}\) is the so-called volatility
process, then the resulting model is a GARCH(1,1) model. Further,
\(\alpha>0\), and \(\beta\) and \(\gamma\) are non-negative constants
satisfying \(\beta+\gamma < 1\). This model can be extended in many ways
and the ARCH/GARCH models are extremely important in finance. A recent
book is Francq and Zakoian (2011). The work on these kind of models was
initiated by Engle (1982), and he was awarded the Nobel prize for his
work. The point as far as Pearson's \(\rho\) is concerned, is that
\(X_t\) and \(X_s\) are uncorrelated for \(t \neq s\), but they are in
fact strongly dependent through the volatility process \(\{h_t\}\),
which can be taken to measure financial risk. This is probably the best
known and most important model class where the dependence structure of
the process is not at all revealed by the autocorreletion function. The
variables are uncorrelated, but contain a dependence structure that is
very important from an economic point of view.

In Figure \ref{fig:issues4} we see some simulated data from a
GARCH(1,1)-model with \(\epsilon_t \sim\) iid \(N(0,1)\),
\(\alpha = 0.1\), \(\beta = 0.7\) and \(\gamma = 0.2\), with \(X_t\) on
the horizontal axis, and \(X_{t-1}\) on the vertical axis. In this
particular case, \(\widehat \rho(X_t, X_{t-1}) = 0.018\), despite the
strong serial dependence that is seen to exist directly from equation
\eqref{eq:garch}.

The nonlinearity issue will be analysed very extensively and quite
systematically in the following sections, but there have also been
various more ad hoc solutions to this problem. We will just mention
briefly two of them here. Slightly more details are given in the survey
paper Tjøstheim (1996), and much more details in the literature cited
there.

\begin{enumerate}
\def\labelenumi{(\alph{enumi})}
\item
  \textbf{Higher moments:} An «obvious» ad hoc solution in the nonlinear
  GARCH case is to compute the product moment correlation on squares
  \(\{X_t^2\}\) instead of \(\{X_t\}\) themselves. It is easily seen
  that the squares are autocorrelated. This is the idea behind the
  McLeod and Li (1983) test. It requires the existence of 4th moments,
  though, which will not always be fulfilled for models of financial
  time series that typically have heavy tails, see e.g.~Teräsvirta et
  al. (2010, Ch. 8).
\item
  \textbf{Frequency based tests:} These are also based on higher product
  moments, but in this instance one takes the Fourier transform of these
  to obtain the so-called bi-spectrum and tri-spectrum, on which in turn
  independence tests can be based (Subba Rao and Gabr 1980; Hinich
  1982).
\end{enumerate}

In the following sections we will look at ways of detecting nonlinear
and non-Gaussian structures by going beyond Pearson's \(\rho\).

\section{\texorpdfstring{Beyond Pearson's \(\rho\): The
copula}{Beyond Pearson's \textbackslash{}rho: The copula}}\label{copula}

For two variables one may ask, why not just take the joint density
function \(f(x,y)\) or the cumulative distribution function \(F(x,y)\)
as a descriptor of the joint dependence? The answer is quite obvious. If
a parametric density model is considered, it is usually quite difficult
to give an interpretation of the parameters in terms of the strength of
the dependence. An exception is the multivariate normal distribution of
course, but even for elliptical distributions the «correlation»
parameter \(\rho\) is not, as we have seen, necessarily a good measure
of dependence. If one looks at nonparametric estimates for multivariate
density functions, to a certain degree one may get an informal
indication of strength of dependence in certain regions from a display
of the density, but the problems increase quickly with dimension due
both to difficulties of producing a graphical display and to the lack of
precision of the estimates due to the curse of dimensionality.

Another problem in analyzing a joint density function is that it may be
difficult to disentangle effects due to the shape of marginal
distributions and effects due to dependence among the variables
involved. This last problem is resolved by the copula construction.
Sklar's (1959) theorem states that a multivariate cumulative
distribution function \(F(x) = F(x_1,\ldots,x_p)\) with marginals
\(F_i(x_i)\), \(i=1,\ldots,p\) can be decomposed as

\begin{equation}
F(x_1,\ldots,x_p) = C(F_1(x_1),\ldots, F_p(x_p))
\label{eq:copula1}
\end{equation}

where \(C(u_1,\ldots u_p)\) is a distribution function over the unit
cube \([0,1]^p\). Klaassen and Wellner (1997) point out that Hoeffding
(1940) had the basic idea of summarizing the dependence properties of a
multivariate distribution by its associated copula, but he chose to
define the corresponding function on the interval \([-1/2,1/2]\) instead
of on the interval \([0,1]\). In the continuous case, \(C\) is a
function of uniform variables \(U_1,\ldots,U_p\), using the well-known
fact that for a continuous random variable \(X_i\), \(F_i(X_i)\) is
uniform on {[}0,1{]}. Further, in the continuous case \(C\) is uniquely
determined by Sklar's (1959) theorem.

The theorem continues to hold for discrete variables under certain weak
regularity conditions securing uniqueness. We refer to Nelsen (1999) and
Joe (2014) for extensive treatments of the copula. Joe (2014), in
particular, contains a large section on copulas in the discrete case.
See also Genest and Nešlehová (2007). For simplicity and in keeping with
the assumptions in the rest of this paper we will mostly limit ourselves
to the continuous case.

The decomposition \eqref{eq:copula1} very effectively disentangle the
distributional properties of a multivariate distribution into a
dependence part measured by the copula C and a marginal part described
by the univariate marginals. Note that \(C\) is invariant with respect
to one-to-one transformations of the marginal variables \(X_i\). In this
respect it is analogous to the invariance of the Kendall and Spearman
rank based correlation coefficients.

A representation in terms of uniform variables can be said to be in
accord with a statistical principle that complicated models should
preferably be represented in terms of the most simple variables
possible, in this case uniform random variables. A possible disadvantage
of the multivariate uniform distribution is that tail behavior of
distributions may be difficult to discern on the uniform scale, as it
may result in singular type behavior in the corners of the uniform
distribution with accumulations of points there in a scatter diagram on
\([0,1]^2\) or \([0,1]^{p}\). It is therefore sometimes an advantage to
change the scale to a standard normal scale, where the uniform scores
\(U_i\) are replaced by standard normal scores \(\Phi^{-1}(U_i)\) with
\(\Phi\) being the cumulative distribution of the standard normal
distribution. This leads to a more clear representation of tail
properties. This scale is sometimes used in copula theory (see e.g. Joe
(2014)), and we have used it systematically in our work on local
Gaussian approximation described in Section \ref{lgc}.

The decomposition in \eqref{eq:copula1} is very useful in that it leads to
large classes of models that can be specified by defining the marginals
and the copula function separately. It has great flexibility in that
very different models can be chosen for the marginal distribution, and
there is a large catalog of possible parametric models available for the
copula function \(C\); it can also be estimated nonparametrically. The
simplest one is the Gaussian copula. It is constructed from a
multivariate Gaussian distribution \(\Phi_{\Sigma}\) with correlation
matrix \(\Sigma\). It is defined by\\

\begin{equation}
C_{\Sigma}(u) = \Phi_{\Sigma}\left(\Phi^{-1}(u_1),\ldots,\Phi^{-1}(u_p)\right)
\label{eq:copula2}
\end{equation}

such that \(Z_i = \Phi^{-1}(U_i)\) are standard normal variables for
\(i=1,\ldots,p\). It should be carefully noted that if one uses
\eqref{eq:copula2} in model building, one is still allowed to put in a
marginal cumulative distribution functions of one's own choice,
resulting in a joint distribution that is not Gaussian. A multivariate
Gaussian distribution with correlation matrix \(\Sigma\) is obtained if
the marginals are univariate Gaussians. If the marginals are not
Gaussians the correlation matrix in the distribution obtained by
\eqref{eq:copula1} will not in general be \(\Sigma\). Klaassen and Wellner
(1997) present an interesting optimality property of the normal scores
rank correlation coefficient, the van der Waerden correlation, as an
estimate of \(\Sigma\).

A similar construction taking as its departure the multivariate
\(t\)-distribution can be used to obtain a \(t\)-copula.

A general family of copulas is the family of Archimedean copulas. It is
useful because it can be defined in an arbitrary dimension \(p\) with
only one parameter \(\theta\) belonging to some parameter space
\(\Theta\). A copula \(C\) is called Archimedean if it has the
representation

\begin{equation}
C(u,\theta) = \psi^{[-1]}(\psi(u_1,\theta)+ \cdots + \psi(u_p,\theta);\theta)
\label{eq:copula3}
\end{equation}

where \(\psi: [0,1] \times \Theta \to [0,\infty)\) is a continuous,
strictly decreasing and convex function such that
\(\psi(1,\theta) = 0\). Moreover, \(\theta\) is a parameter within some
parameter space \(\Theta\). The function \(\psi\) is called the
generator function and \(\psi^{[-1]}\) is the pseudo inverse of
\(\psi\). We refer to Joe (2014) and Nelsen (1999) for more details and
added regularity conditions. In practice, copulas have been mostly used
in the bivariate situation, in which case there are many special cases
of the Archimedian copula \eqref{eq:copula3}, such as the Clayton, Gumbel
and Frank copulas. In particular the Clayton copula has been important
in economics and finance. It is defined by

\begin{equation}
C_C(u_1,u_2) = \max \{u_1^{-\theta}+u_2^{-\theta}-1;0]^{-1/\theta}, \quad \mbox{with} \quad \theta \in [-1,\infty)/{0}.
\label{eq:copula4}
\end{equation}

\begin{figure}[t]
\subfloat[The observed log-returns of the daily data \,\,\,\label{fig:return1}]{\includegraphics[width=0.49\linewidth]{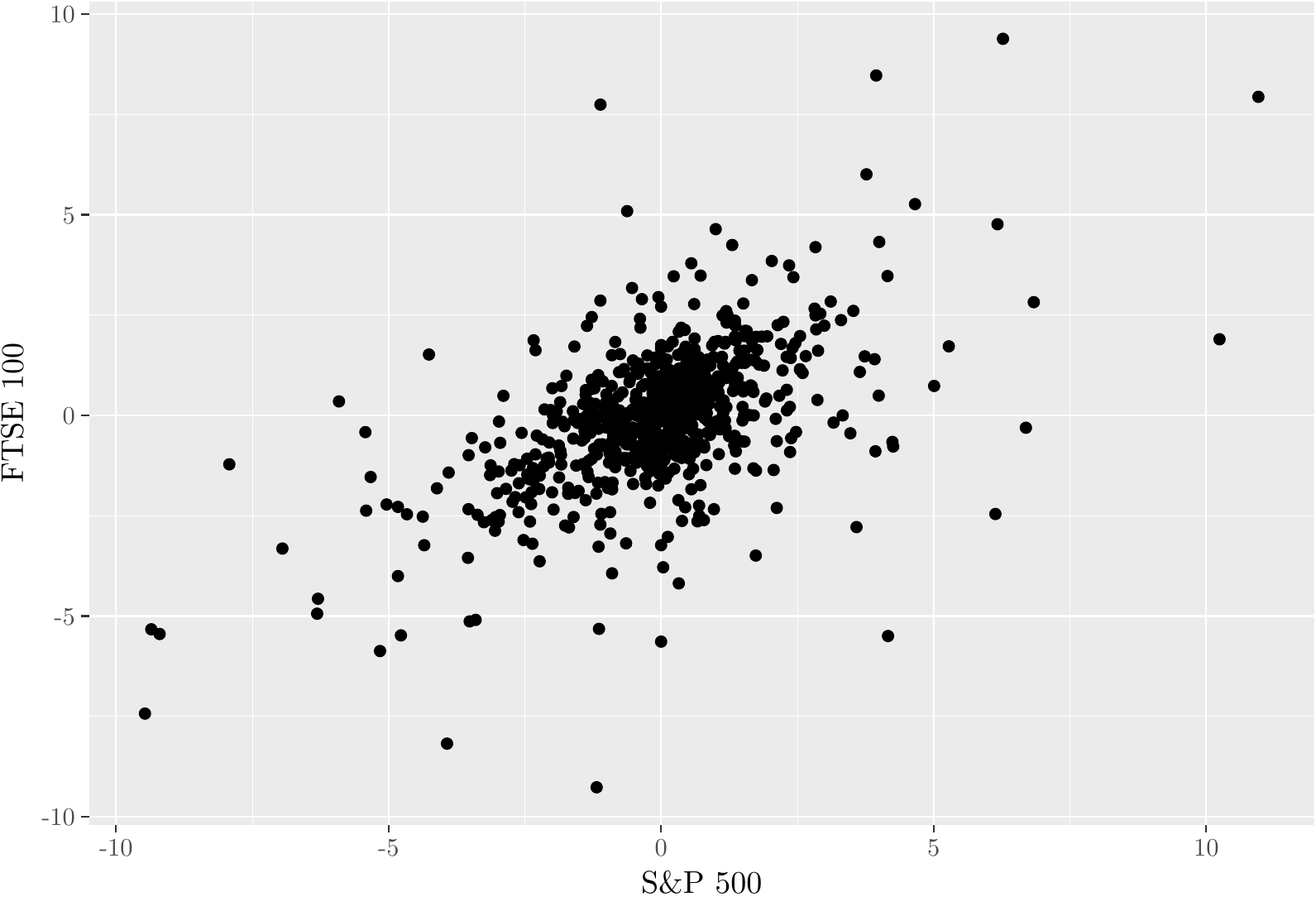} }\subfloat[Uniform scores of the financial returns data set\label{fig:return2}]{\includegraphics[width=0.49\linewidth]{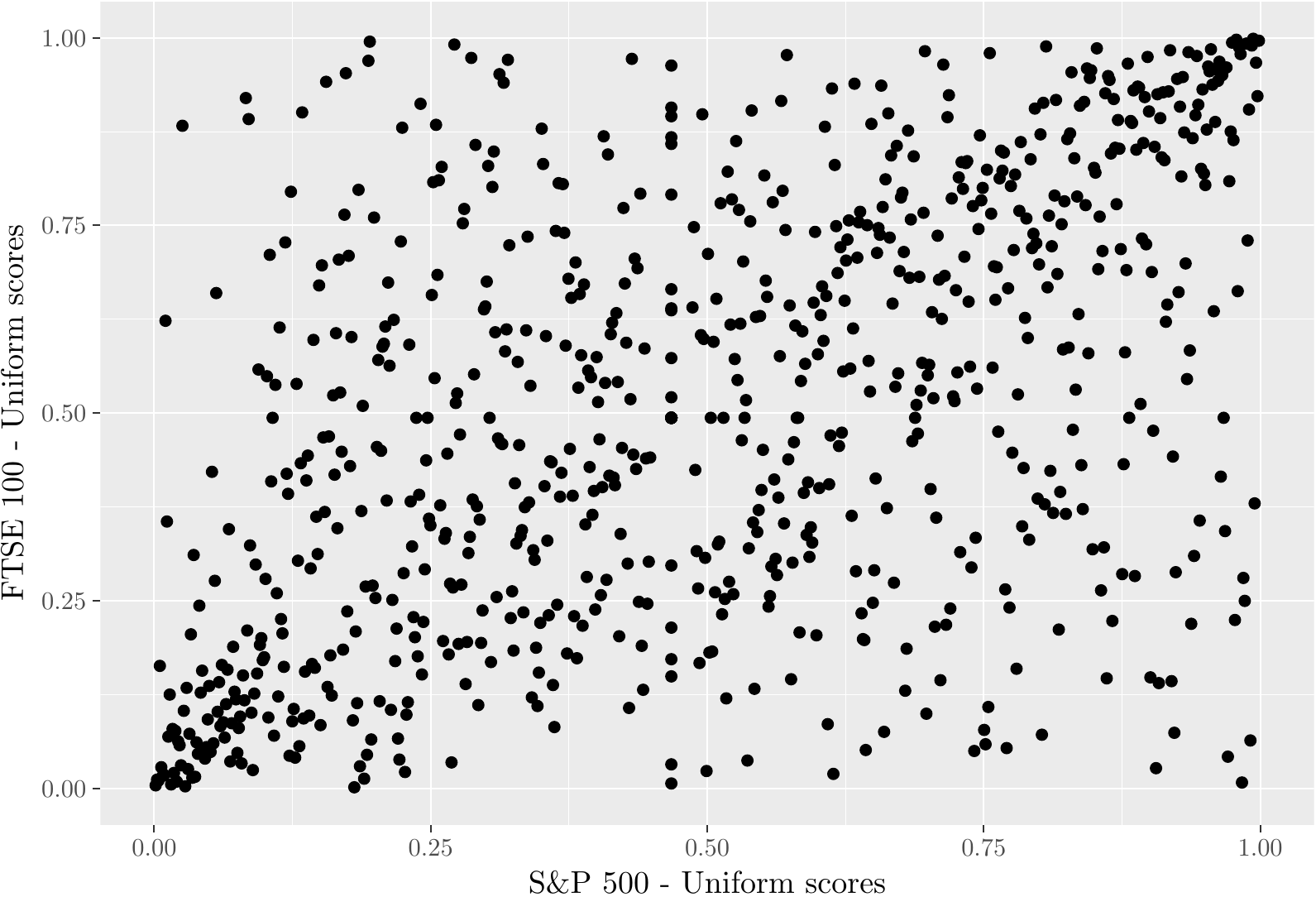} }\newline\subfloat[Normal scores of the financial returns data \,\,\, set\label{fig:return3}]{\includegraphics[width=0.49\linewidth]{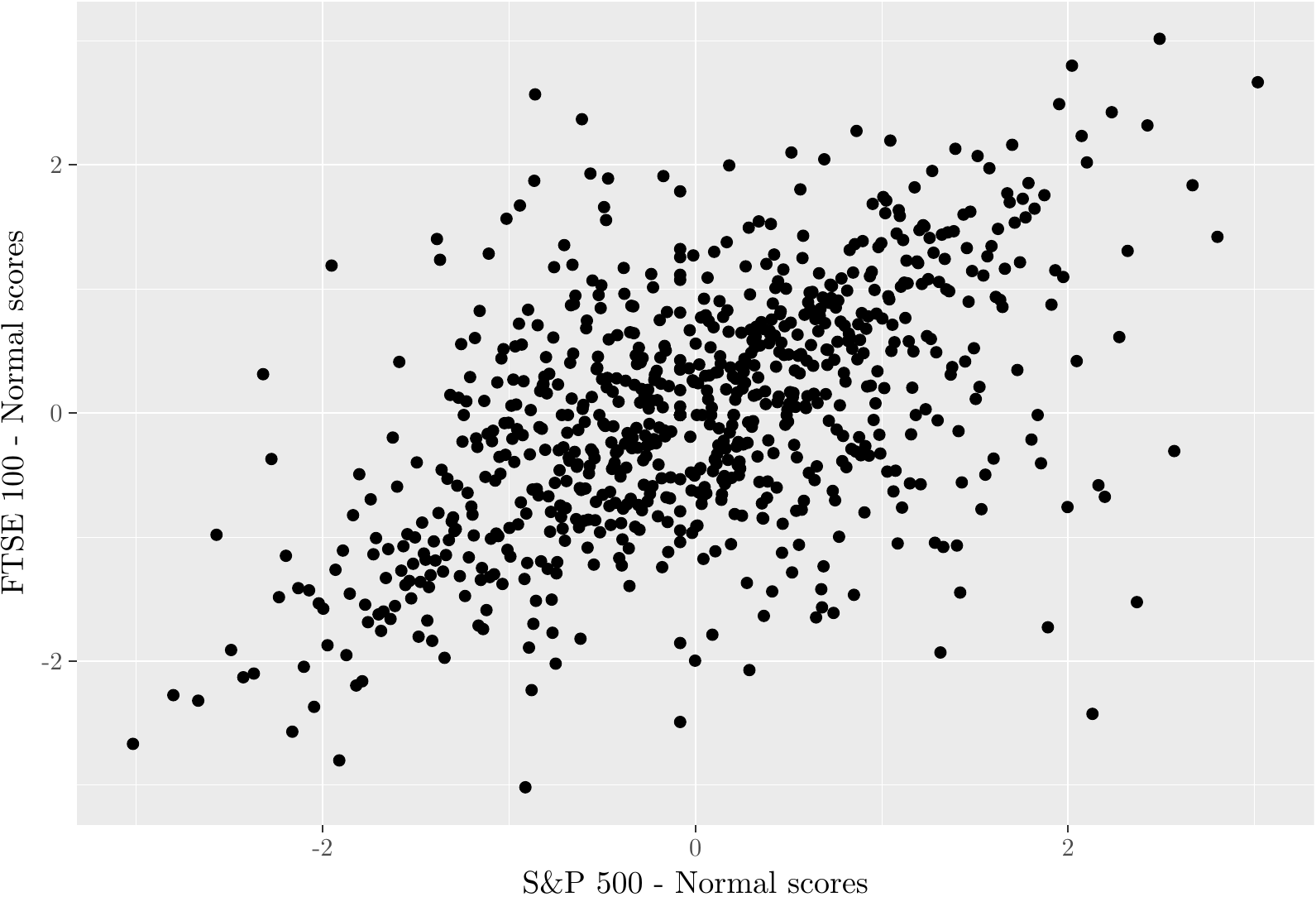} }\subfloat[Simulated data from a Clayton copula fitted to the financial returns data set\label{fig:return4}]{\includegraphics[width=0.49\linewidth]{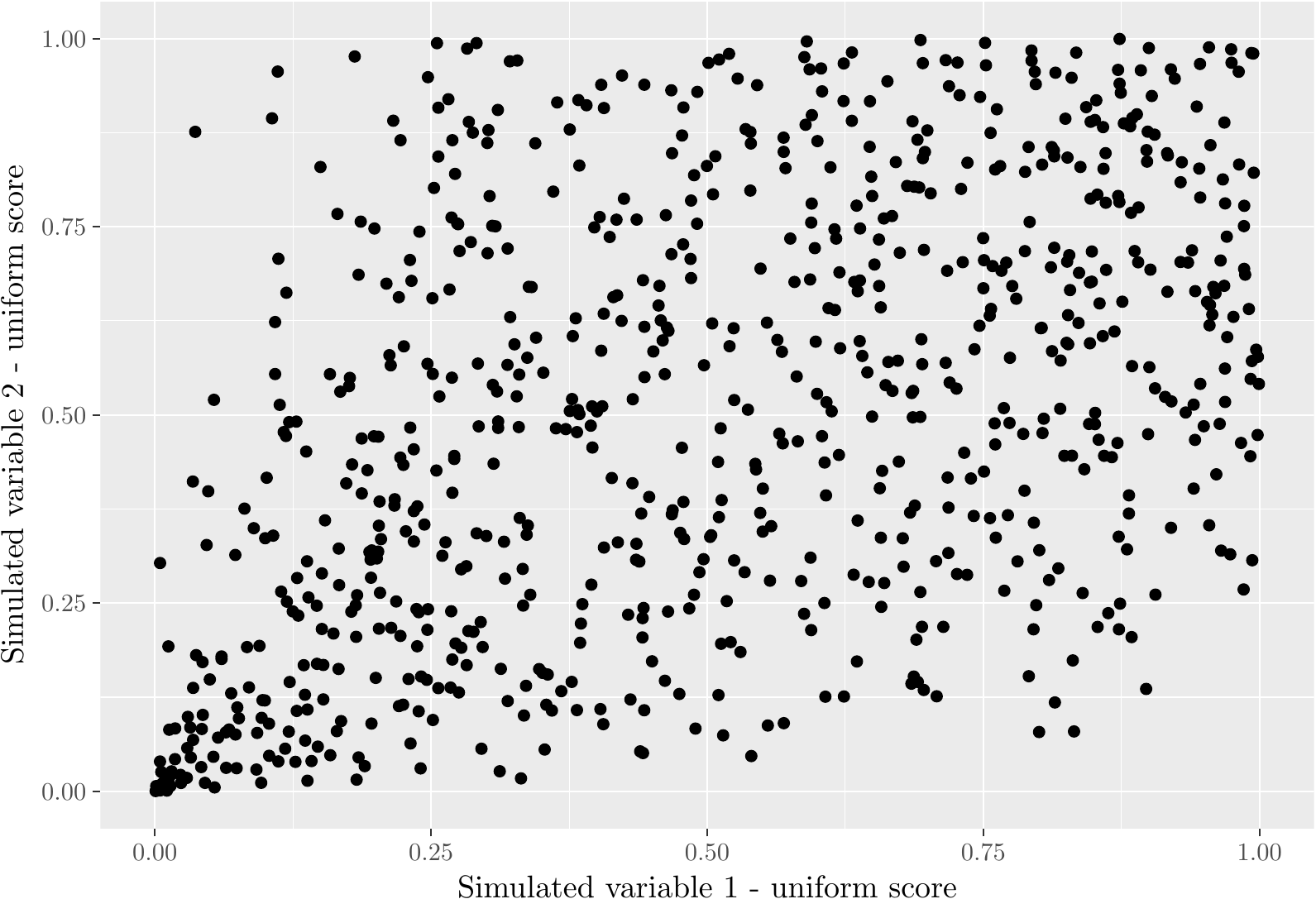} }\caption{Illustrations using the financial returns data set}\label{fig:return}
\end{figure}

We will throughout this paper illustrate several points using a
bivariate data set on some financial returns. We use daily international
equity price index data for the United States (i.e.~the S\&P 500) and
the United Kingdom (i.e.~the FTSE 100). The data are obtained from
Datastream (2018), and the returns are defined as
\[r_t = 100 \times \left(\log(p_t ) - \log(p_{t-1})\right),\] where
\(p_t\) is the price index at time \(t\). The observation span covers
the period from January 1st 2007 through December 31st 2009, in total
784 observations. In Figure \ref{fig:return}, four scatterplots are
presented.

Figure \ref{fig:return1} displays a scatterplot of the observed
log-returns, with S\&P 500 on the horizontal axis, and the FTSE 100 on
the vertical axis. Figure \ref{fig:return2} displays the uniform scores
of the same data, and we see indications of a singular behavior of the
copula density in the lower left and upper right corners of the unit
square. In Figure \ref{fig:return3} the observations have been
transformed to normal scores, which more clearly reveals the tail
properties of the underlying distribution. Finally, Figure
\ref{fig:return4} shows the scatter plot of 784 simulated pairs of
variables, on uniform scale, from a Clayton copula fitted to the return
data. This plot partially resembles Figure \ref{fig:return2}, in
particular in the lower left corner. However, there are some differences
in the upper right corner. We will look into this discrepancy in Section
\ref{lgc}.

In Figure \ref{fig:return1}, and perhaps more clearly in Figure
\ref{fig:return3}, we see that there seems to be stronger dependence
between the variables when the market is going either up or down, which
is very sensible from an economic point of view, but it is not easy to
give an interpretation of the parameter \(\theta\) of the Clayton copula
in terms of such type of dependence. In fact, in this particular case,
\(\widehat\theta = 0.96\). The difficulty of giving a clear and concrete
interpretation of copula parameters in terms of measuring strength of
dependence can be stated as a potential issue of the copula
representation. In this respect it is very different from the Pearson's
\(\rho\). We will return to this point in much detail in Section
\ref{lgc}, where we define a local correlation.

Another issue of the original copula approach has been the lack of good
practical models as the dimension increases, as it would for example in
a portfolio problem in finance. This has recently been sought solved by
the so-called pair copula construction. To simplify, in a trivariate
density \(f(x_1,x_2,x_3)\), by conditioning this can be written
\(f(x_1,x_2,x_3) = f_1(x_1)f_{23|1}(x_2,x_3|x_1)\), and a bivariate
copula construction, e.g.~a Clayton copula, can be applied to the
conditional density \(f_{23|1}(x_2,x_3|x_1)\) with \(x_1\) fixed. This
conditioning can be extended to higher dimensions under a few
simplifying assumptions, resulting in a so-called vine copula, of which
there are several types. The procedure is well described by Aas et al.
(2009), and has found a number of applications. The Clayton canonical
vine copula, for instance, allows for the occurrence of very strongly
correlated downside events and has been successfully applied in
portfolio choice and risk management operations. The model is able to
reduce the effects of extreme downside correlations and produces
improved statistical and economical performance compared to elliptical
type copulas such as the Gaussian copula \eqref{eq:copula2} and the
\(t\)-copula, see Low et al. (2013).

Other models developed for risk management applications are so-called
panic copulas to analyze the effect of panic regimes in the portfolio
profit and loss distribution, see e.g Meucci (2011). A panic reaction is
taken to mean that a number of investors react in the same way, such
that the statistical dependence becomes very strong between financial
returns from various financial objects, in this way rendering the risk
spreading of the portfolio illusory. We will return to this situation in
Section \ref{lgc} where among other things we can show that in a panic
situation the local correlation increases and approaches one as a
function of a copula parameter. The copula has also been used directly
for independence testing; see e.g. Genest and Rémillard (2004) and
Mangold (2017).

Most of the copula theory and also most of the applications are to
variables that are assumed to be iid, but there is also a growing
literature on stochastic processes such as Markov chains. The existence
of both auto dependence and cross dependence in a multivariate
stochastic process is quite challenging. Some of the mathematical
difficulties in the Markov chain case is clearly displayed in the paper
Darsow, Nguyen, and Olsen (1992). They used the ordinary copula, but it
is not obvious how the theory of Markov processes can be helped by the
concept of a copula. That work was limited to first order Markov chain.
The pair copula has also been introduced in a Markov theory framework,
and then in higher-order Markov processes, by Ibragimov (2009). Again,
so far, the impact on Markov theory has not been overwhelming. This may
partly be due to complicated technical conditions.

Two other papers using copulas (and pair copulas) in serial dependence
are Beare (2010) and Smith et al. (2010). When it comes to parametric
time series analysis, especially for multivariate time series, it has
been easier to implement the copula concept as developed for iid
variables. This is well documented in the survey paper by Patton (2012).
The reason is that the auto dependence can first be filtered out by a
marginal fit to each component series, and the copula could then be
applied to the residuals which may be assumed to be iid or at least can
be replaced by an iid vector process asymptotically. More precisely in
the framework of Patton and others,

\begin{equation}
X_{it} = \mu_i(Z_{t-1};\phi)+\sigma_i(Z_{t-1};\phi)\varepsilon_{it}, \; i=1,2,
\label{eq:copula5}
\end{equation}

where in the bivariate case primarily considered by Patton,
\[Z_{t-1} \in {\cal F}_{t-1}, \quad \varepsilon_{it} \sim F_{it}, \qquad \bm{\varepsilon}_t|{\cal F}_{t-1} \sim F_{\varepsilon t} = C_t(F_{1t},F_{2t}).\]
Here \({\cal F}_{t-1}\) can be taken as the \(\sigma\)-algebra generated
by \(X_{s}\) for \(s \leq t\). and \(Z_{t-1}\) is a stochastic vector
variable, e.g.~higher lags of \(X_t\) measurable with respect to
\({\cal F}_{t-1}\). The estimation can be done in two steps, cf.~also
Chen and Fan (2006). First the parameters \(\bm{\phi}_i\) of the
marginal processes are estimated. Then a copula modeling stage is
applied to the estimated residuals
\[ \widehat{\varepsilon}_{it} = \frac{X_{it}-\mu_i(Z_{t-1};\widehat{\phi}_i)}{\sigma_i(Z_{t-1};\widehat{\phi}_i)}.\]
In this context both parametric and nonparametric (resulting in a
semiparametric model) models have been considered for the residual
distribution \(F_{it}\). In the parametric case, time dependence can be
allowed for \(F_{it}\), whereas in the nonparametric case there is no
dependence of \(t\) permitted in \(F_{it}\). Much more details and
references are provided in Patton (2012). The modeling in
\eqref{eq:copula5} is restricted to the bivariate case. Modeling of both
cross and auto dependence, including use of vine copulas, in a
multivariate time series or Markov process is given in Smith (2015).
Time dependent risk is treated using a dynamic copula model by Oh and
Patton (2018).

\section{\texorpdfstring{Beyond Pearson's \(\rho\): Global dependence
functionals and tests of
independence}{Beyond Pearson's \textbackslash{}rho: Global dependence functionals and tests of independence}}\label{global}

Studies of statistical dependence may be said to center mainly around
two problems: (i) definition and estimation of measures of dependence
and (ii) tests of independence. Of course these two themes are closely
related. Measures of association such as the Pearson's \(\rho\) can also
be used in tests of independence, or more precisely: tests of
uncorrelatedness. On the other hand, test functionals for tests of
independence can also in many, but not all, cases be used as a measures
of dependence. A disadvantage with measures derived from tests is that
they are virtually always based on a distance function and therefore
non-negative. This means that they cannot distinguish between negative
and positive dependence, whatever this may mean in the general nonlinear
case. We will return to this later in the paper.

Most of the test functionals are based on the definition of independence
in terms of cumulative distribution functions or in terms of density
functions. Consider \(p\) stochastic variables \(X_1,\ldots,X_p\). These
variables are independent if and only if their joint cumulative
distribution function is the product of the marginal distribution
functions:
\(F_{X_1,\ldots,X_p}(x_1,\ldots,x_p) = F_1(x_1) \cdots F_p(x_p)\), and
the same is true for all subsets of variables of \((X_1,\ldots,X_p)\).
If the variables are continuous, this identity can be phrased in terms
of the corresponding density functions instead. A typical test
functional is then designed to measure the distance between the
estimated joint distributions/densities and the product of the estimated
marginals. This is not so easily done for parametric densities, since
the dependence on parameters in the test functional may be very complex,
and tests of independence may be more sensibly stated in terms of the
parameters themselves, as is certainly the case for the Gaussian
distribution. Therefore one would usually estimate the involved
distributions nonparametrically, which, for joint distributions, may be
problematic for moderate and large \(p\)'s due to the curse of
dimensionality. We will treat these problems in some detail in Sections
\ref{cumulative}-\ref{density}.

Before starting on the description of the various dependence measures,
let us remark that Rényi (1959) proposed that a measure of dependence
between two stochastic variables \(X\) and \(Y\), \(\delta(X,Y)\),
should ideally have the following 7 properties:

\begin{enumerate}
\def\labelenumi{(\roman{enumi})}
\tightlist
\item
  \(\delta(X,Y)\) is defined for any \(X,Y\) neither of which is
  constant with probability 1.
\item
  \(\delta(X,Y) = \delta(Y,X)\).
\item
  \(0 \leq \delta(X,Y) \leq 1\).
\item
  \(\delta(X,Y) = 0\) if and only if \(X\) and \(Y\) are independent.
\item
  \(\delta(X,Y) = 1\) if either \(X=g(Y)\) or \(Y=f(X)\), where \(f\)
  and \(g\) are measurable functions.
\item
  If the Borel-measurable functions \(f\) and \(g\) map the real axis in
  a one-to-one way to itself, then \(\delta(f(X),g(Y)) = \delta(X,Y)\).
\item
  If the joint distribution of \(X\) and \(Y\) is normal, then
  \(\delta(X,Y) = |\rho(X,Y)|\), where \(\rho(X,Y)\) is Pearson's
  \(\rho\).
\end{enumerate}

The product moment correlation \(\rho\) satisfies only (i), (ii) and
(vii).

One can argue that the rules i) - vii) do not take into account the
difference between positive and negative dependence; it only looks at
the strength of the measured dependence. If this wider point of view
were to be taken into account, (iii) could be changed into (iii') :
\(-1 \leq \delta(X,Y) \leq 1\), (v) into (v'): \(\delta(X,Y) = 1\) or
\(\delta(X,Y) = -1\) if there is a deterministic relationship between
\(X\) and \(Y\). Finally, (vii) should be changed into (vii') requiring
\(\delta(X,Y) = \rho(X,Y)\). Moreover, some will argue that property
(vi) may be too strong to require. It means that the strength of
dependence is essentially independent of the marginals as for the copula
case.

We will discuss these properties as we proceed in the paper. Before we
begin surveying the test functionals as announced above, we start with
the maximal correlation which, it will be seen, is intertwined with at
least one of the test functionals to be presented in the sequel.

\subsection{Maximal correlation}\label{maximal}

The maximal correlation is based on the Pearson \(\rho\). It is
constructed to avoid the problem demonstrated in Section
\ref{nonlinearity} that Pearson's \(\rho\) can easily be zero even if
there is strong dependence.

It seems that the maximal correlation was first introduced by Gebelein
(1941). He introduced it as \[S(X,Y) = \sup_{f,g}\rho(f(X),g(Y)),\]
where \(\rho\) is the Pearson's \(\rho\). Here the supremum is taken
over all Borel-measurable functions \(f,g\) with finite and positive
variance for \(f(X)\) and \(g(Y)\). The measure \(S\) gets rid of the
nonlinearity issue of \(\rho\). It is not difficult to check that
\(S=0\) if and only if \(X\) and \(Y\) are independent. On the other
hand \(S\) cannot distinguish between negative and positive dependence,
and it is in general difficult to compute.

The maximal correlation \(S(X,Y)\) cannot be evaluated explicitly except
in special cases, not the least because there does not always exist
functions \(f_0\) and \(g_0\) such that
\(S(X,Y) = \rho(f_0(X),g_0(Y))\). If this equality holds for some
\(f_0\) and \(g_0\), it is said that the maximal correlation between
\(X\) and \(Y\) is attained. Rényi (1959) gave a characterization of
attainability.

Czáki and Fischer (1963) studied mathematical properties of the maximal
correlation and computed it for a number of examples. Abrahams and
Thomas (1980) considered maximal correlation in the context of
stochastic processes. A multivariate version of maximal correlation was
proposed in Koyak (1987). In a rather influential paper, at least at the
time, Breiman and Friedman (1985) presented the ACE (alternating
conditional expectation) algorithm for estimating the optimal functions
\(f\) and \(g\) in the definition of the maximal correlation. They
applied it both to correlation and regression. Some curious aspects of
the ACE algorithm is highlighted in Hastie and Tibshirani (1990,
84--86).

Two more recent publications are Huang (2010), where the maximal
correlation is used to test for conditional independence, and Yenigün,
Székely, and Rizzo (2011), where it is used to test for independence in
contingency tables. The latter paper introduces a new example where
\(S(X,Y)\) can be explicitly computed. See also Yenigün and Rizzo
(2014).

\subsection{Measures and tests based on the distribution
function}\label{cumulative}

We start with, and in fact put the main emphasis on, the bivariate case.
Let \(X\) and \(Y\) be stochastic variables with cumulative distribution
functions \(F_X\) and \(F_Y\). The problem of measuring the dependence
between \(X\) and \(Y\) can then be formulated as a problem of measuring
the distance between the joint cumulative distribution function
\(F_{X,Y}\) of \((X,Y)\) and the distribution function \(F_XF_Y\) formed
by taking the product of the marginals. Let \(\Delta(\cdot,\cdot)\) be a
candidate for such a distance functional. It will be assumed that
\(\Delta\) is a metric, and it is natural to require, Skaug and
Tjøstheim (1996), that

\begin{align}
& \Delta(F_{X,Y},F_XF_Y) \geq 0 \quad \textrm{ and } \quad \Delta(F_{X,Y},F_XF_Y) = 0 \nonumber\\ 
& \textrm{if and only if} \quad F_{X,Y} = F_XF_Y.
\label{eq:zero}
\end{align}

Clearly, such a measure is capable only of measuring the strength of
dependence, not its direction. Corresponding to requirement (vi) in
Rényi's scheme listed in the beginning of this section, one may want to
require invariance under transformations, or more precisely

\begin{equation}
\Delta(F_{X,Y}^{*},F_X^{*}F_Y^{*}) = \Delta(F_{X,Y},F_XF_Y) 
\label{eq:invariance}
\end{equation}

where \(F_X^{*}(x) = F_X(g(x))\), \(F_Y^{*}(y) = F_Y(h(y))\) and
\(F_{X,Y}^{*}(x,y) = F_{X,Y}\{g(x),h(y)\}\). Here \(g\) and \(h\) are
increasing functions, and \(F_X^{*}, F_Y^{*}, F_{X,Y}^{*}\) are the
marginals and bivariate distribution functions of the random variables
\(g^{-1}(X), h^{-1}(Y)\) and \(\{g^{-1}(X),h^{-1}(Y)\}\), respectively.

For distance functionals not satisfying equation \eqref{eq:invariance},
one can at least obtain scale and location invariance by standardizing
such that \(E(X) = E(Y) = 0\) and \(\Var(X) = \Var(Y)=1\), assuming that
the second moment exists. In practice, empirical averages and variances
must be employed, but asymptotically the difference between using
empirical and theoretical quantities is a second order effect. In Skaug
and Tjøstheim (1996) and Tjøstheim (1996) such a standardization has
been employed for all functionals.

Pearson's \(\rho\) can be expressed as a functional on \(F_X\), \(F_Y\)
and \(F_{X,Y}\) although not generally as a distance functional
depending on \(F_{X,Y}\) and \(F_XF_Y\). For instance, with \(X\) and
\(Y\) standardized, the Pearson correlation squared can be expressed,
\[\rho^2 = \left\{\int xy \di F_{X,Y}(x,y) - \int x\di F_X(x)\int y \di F_Y(y)\right\}^2.\]
Clearly \(\rho^2\) does not satisfy either of the conditions
\eqref{eq:zero} or \eqref{eq:invariance}. Similarly, the square of the
population values of the Spearman rank correlation and Kendall's
\(\tau\) are obtained by squaring in the formulas \eqref{eq:spearman} and
\eqref{eq:kendall}. For these measures the requirement \eqref{eq:zero} is
not fulfilled, whereas the invariance property does hold.

A natural estimate \(\widehat{\Delta}\) of a distance functional
\(\Delta\) is obtained by setting
\[\widehat{\Delta}(F_{X,Y},F_XF_Y) = \Delta(\widehat{F}_{X,Y},\widehat{F}_X\widehat{F}_Y),\]
where \(\widehat{F}\) may be taken to be the empirical distribution
functions given by
\[\widehat{F}_X(x) = \frac{1}{n}\sum_{j=1}^{n}1(X_{j} \leq x) \quad \widehat{F}_Y(y) = \frac{1}{n} \sum_{j=1}^{n}1(Y_{j} \leq y)\]
and
\[\quad F_{X,Y}(x,y) = \frac{1}{n}\sum_{j=1}^{n}1(X_{j} \leq x)1(Y_{j} \leq y),\]
or a normalized version with \(n^{-1}\) replaced by \((n+1)^{-1}\) for
given observations \(\{(X_{1},Y_{1},\ldots,(X_{n},Y_{n})\}\). Similarly,
for a stationary time series \(\{X_t\}\) at lag \(k\),
\[\widehat{F}_k(x_1,x_2) \doteq \widehat{F}_{X_t,X_{t-k}}(x_1,x_2) = \frac{1}{n-k} \sum_{t=k+1}^{n}1(X_t \leq x_1)1(X_{t-k} \leq x_2).\]

Conventional distance measures between two distribution functions \(F\)
and \(G\) are the Kolmogorov-Smirnov distance
\[\Delta_1(F,G) = \sup_{(x,y)}|F(x,y)-G(x,y)|\] and the Cramér-von Mises
type distance of a distribution \(G\) from a distribution \(F\)
\[\Delta_2(F,G) = \int \left\{F(x,y)-G(x,y)\right\}^2 \di F(x,y)\] Here
\(\Delta_1\) satisfies \eqref{eq:zero} and \eqref{eq:invariance}, whereas
\(\Delta_2\) satisfies \eqref{eq:zero} but not \eqref{eq:invariance}.

Most of the work pertaining to measuring dependence and testing of
independence has been done in terms of the Cramér-von Mises distance.
This work started already by Hoeffding (1948) who looked at iid pairs
\((X_{i},Y_{i})\), and studied finite sample distributions in some
special cases. With considerable justification it has been named the
Hoeffding-functional by some. This work was continued by Blum, Kiefer,
and Rosenblatt (1961) who provided an asymptotic theory, still for the
iid case. It was extended to the time series case with a resulting test
of serial independence in Skaug and Tjøstheim (1993a). A recent paper
using a copula framework is Kojadinovic and Holmes (2009). We will
briefly review the time series case because it illustrates some of the
problems, and because some of the same essential ideas as for the
Hoeffding-functional have been used in more recent work on the distance
covariance, which we treat in Section \ref{distance}. Since
\(F_{X_{t}}(x)\) does not depend on \(t\) for a stationary time series,
we simply write \(F(x)\) in the following.

In the time series case the Cramér-von Mises distance at lag \(k\) is
given by
\[D_k = \int \Big(F_k(x_1,x_2)-F(x_1)F(x_2)\Big)^2 \,\textrm{d}F_k(x_1,x_2),\]
where \(F_k\) and \(F\) are the joint and marginal distributions of
\((X_t,X_{t-k})\) and \(X_t\), respectively. Replacing theoretical
distributions by empirical ones leads to the estimate

\begin{equation}
\widehat{D}_k = \frac{1}{n-k} \sum_{t=k+1}^{n} \Big(\widehat{F}(X_t,X_{t-k})-\widehat{F}(X_t)\widehat{F}(X_{t-k})\Big)^2.
\label{eq:Hoeffding1}
\end{equation}

Assuming \(\{X_t\}\) to be ergodic, we have that
\(\widehat{D}_k \to D_k\) almost surely as \(n \to \infty\).

To construct a test of serial independence we need the distribution of
\(\widehat{D}_k\) under the assumption of the null hypothesis of
\(\{X_t\}\) being iid. Let
\(Z_t = (Z_t^{(1)},Z_t^{(2)}) \doteq (X_t,X_{t-k})\). Then it is
possible to represent \(\widehat{D}_k\) as \[
\widehat{D}_k = \frac{1}{n^2}\sum_{s,t=1}^{n}h(Z_s,Z_t)+O_p(n^{-3/2}),
\] where \(V_n \doteq n^{-2}\sum_{s,t=1}^{n}h(Z_s,Z_t)\) is a von Mises
U-statistic in the technical sense of Denker and Keller (1983) with a
degenerate symmetric kernel function. Using asymptotic theory, Carlstein
(1988), Denker and Keller (1983) and Skaug (1993), of this statistic or
the related U-statistic one has (Skaug and Tjøstheim 1993a, Theorem 2)
the convergence in distribution

\begin{equation}
n\widehat{D}_k \buildrel {{\cal L}} \over {\to} \sum_{i,j=1}^{\infty} \eta_i\eta_j W_{ij}^2 \quad \textrm{as} \quad n \to \infty
\label{eq:Hoeffding2}
\end{equation}

where \(\{W_{ij}\}\) is an independent identically distributed sequence
of \({\cal N}(0,1)\) variables, and where the \(\{\eta_m\}\) are the
eigenvalues of the eigenvalue problem
\[\int g(x,y)h(y)\,\textrm{d}F(y) = \eta h(x)\] with
\[g(x,y) = \int \{1(x \leq w)-F(w)\}\{1(y \leq w)-F(w)\}\,\textrm{d}F(w).\]
If the distribution of each \(X_t\) is continuous, then
\(\widehat{D}_k\) is distribution free, i.e., its distribution does not
depend on \(F\). Then all calculations can be carried out with \(F\)
being the uniform \([0,1]\) distribution, in which case
\(g(x,y) = -\max(x,y)+(x^2+y^2)/2 + 1/3\) and \(\eta_m = (m\pi)^{-2}\),
and the distribution of the right hand side of \eqref{eq:Hoeffding2} can
be tabulated by truncating it for a large value of the summation index.
Similar distribution results will be seen to hold for test functionals
in Sections \ref{distance} and \ref{hsic}.

A test of the null hypothesis of independence, or rather pairwise
independence at lag \(k\), can now be constructed. It is then reasonable
to reject \(H_0\) if large values of \(\widehat{D}_k\) is observed. Thus
a test of level \(\varepsilon\) is:
\[\textrm{reject} \,\, H_0 \,\, \textrm{if} \,\, n\widehat{D}_k > u_{n,\varepsilon},\]
where \(u_{n,\varepsilon}\) is the upper \(\varepsilon\)-point in the
null distribution of \(n\widehat{D}_k\). Since the exact distribution of
\(\widehat{D}_k\) is unknown, we can use the asymptotic approximation
furnished by retaining a finite number of terms in \eqref{eq:Hoeffding2}.
For \(n=50,100\) and \(k\) small this works well. However, as \(k\)
increases, in general (Skaug and Tjøstheim 1993a) the level is severely
overestimated. The results of Skaug and Tjøstheim (1993a) have since
been very considerably extended and improved by Hong (1998).

Under the hypothesis of \(\{X_t\}\) being iid the bootstrap is a natural
tool to use for constructing the null distribution and critical values.
For moderate and large \(k\)'s the bootstrapping yields a substantially
better approximation to the level.

Under the alternative hypothesis that \(X_t\) and \(X_{t-k}\) are
dependent, the test statistic \(\widehat{D}_k\) will in general be
asymptotically normal with a different rate from that of
\eqref{eq:Hoeffding2}, but the power function will be complicated; see
e.g. Hong (2000).

To extend the scope to testing of serial dependence among
\((X_t,\ldots,X_{t-k})\), or alternatively between a set of several
random variables for which there are iid observations of the set, one
might use a functional

\begin{align}
& \widehat{D}_{1,\ldots,k} \nonumber \\ 
& = \frac{1}{n}\sum_{t=k+1}^{n} \{\widehat{F}_{1,\dots,k}(X_t,X_{t-1},\ldots,X_{t-k})-\widehat{F}(X_t)\widehat{F}(X_{t-1})\cdots\widehat{F}(X_{t-k})\}^2.
\label{eq:Hoeffding3}
\end{align}

The asymptotic theory under the null hypothesis of independence for such
a test has been examined by Delgado (1996) using empirical process
theory, but due to the curse of dimensionality, problems in practice
might be expected for moderately large \(k\)'s. As an alternative Skaug
and Tjøstheim (1993a) used a ``Box-Pierce-Ljung'' analogy, testing for
pairwise independence in all of the pairs
\((X_t,X_{t-1}),(X_t,X_{t-2}),\ldots,(X_t,X_{t-k})\) using the statistic

\begin{equation}
\widehat{D}^{(k)} = \sum_{i=1}^{k} \widehat{D}_i.
\label{eq:Box-test0}
\end{equation}

The asymptotic theory of such a test is given in Skaug (1993). Hong
(1998) noted that

\begin{equation}
\widehat{D}^{(k) *} = \sum_{i=1}^{k}(n-i)\widehat{D}_i
\label{eq:Box-test1}
\end{equation}

has better size properties for large \(k\).

There have been several contributions to the limit theory of statistics
such as \eqref{eq:Hoeffding1} and \eqref{eq:Hoeffding3} using empirical
process theory. Delgado did that based on developments in Blum, Kiefer,
and Rosenblatt (1961), but the limiting Gaussian process is complicated
with a complex covariance matrix which makes it difficult to tabulate
critical values. Ghoudi, Kulperger, and Rémillard (2001) based their
work on the so-called Möbius transformation promoted by Deheuvels
(1981a; 1981b) in his papers on independence testing. This
transformation takes explicitly into account the joint distribution
function of all subsets of \((X_1, \ldots, X_p)\) mentioned in the
second paragraph of Section \ref{copula}. To explain the Möbius
transformation, let \(X_1, \ldots, X_k\), \(k \geq 2\) be random
variables. For \(1 \leq j \leq k\) let \(F_j\) denote the marginal
cumulative distribution function of \(X_j\) and let \(F_{1,\ldots,k}\)
be the corresponding joint cumulative distribution function. Consider a
subset \(A \subset I_k = \{1,\ldots,k\}\), and for any
\(x \in \mathbb{R}^{k}\), define
\[\mu_A(x) = \sum_{B \subset A} (-1)^{|A-B|}F_{1,\ldots,k}(x^{B})\prod_{j \in A-B}F_j(x_j),\]
where \(|C|\) is the number of elements in a set \(C\),
\(A-B = A \cap B^{c}\), where by convention \(\prod_{\emptyset}=1\), and
where \(x^{B}\) is the vector whose \(i\)th component is defined by
\(x_i\) if \(i \in B\), and \(x_i = \infty\) otherwise. Then one can
state the following criterion of independence: \(X_1,\ldots,X_k\) are
independent if and only if \(\mu_A = 0\) for any \(A \subset I_k\). This
is shown in e.g. Ghoudi, Kulperger, and Rémillard (2001). In that paper
it is also shown how this transformation leads to a Gaussian empirical
process limit with a relatively simple covariance function, making it
easier to tabulate critical values. The authors manage to do this both
for the Cramér-von Mises statistic and the Kolmogorov-Smirnov statistic,
and they consider three cases: iid vector samples, time series samples
and residuals in time series models. See also Beran, Bilodeau, and
Lafaye de Micheaux (2007).

Genest and Rémillard (2004) use the Möbius transformation in testing of
independence in a copula framework, and Ghoudi and Rémillard (2018) use
it to obtain tests for independence of residuals in a parametric model
\[X_i = \mu_i(\theta)+\sigma_i(\theta)\varepsilon_i.\]\\
where the iid innovations \(\{\varepsilon_i\}\) have mean 0 and variance
1.

As an alternative to testing uncorrelatedness in time series using the
Pearson \(\rho\) at accumulated lags as in the Box-Ljung-Pierce
statistic, one could test for a constant spectral density. Quite general
specification tests in terms of the spectral density \(f(\omega)\) has
been considered by Anderson (1993), who looked at tests of
\(H_0: f(\omega) =\) constant by using both the Cramér-von Mises and the
Kolmogorov-Smirnov criteria. Hong (2000) introduced a spectral
counterpart for the independence tests based on the marginal
distribution function \(F\) for a time series and the lag \(k\)
distribution function \(F_k\) introduced earlier in this section. This
was achieved by replacing the ordinary autocorrelation function
\(\rho_k\) by the dependence measure
\[\rho_k^{*}(u,v) \doteq F_k(u,v)-F(u)F(v) = \Cov\{1(X_t \leq u),1(X_{t-k} \leq v)\}\]
and then taking Fourier transforms, which leads to

\begin{equation}
h(\omega,u,v) = (2\pi)^{-1}\sum_{k=-\infty}^{\infty}\rho_{k}(u,v)e^{-ik\omega}
\label{eq:Hong1}
\end{equation}

where \(\rho_k(u,v) = \rho_{|k|}^{*}(u,v)\). It is shown in Hong (2000)
how this can be used to construct a test of independence that has power
in cases where the tests based on the ordinary spectrum has not.
However, as explained in a note in that paper, weak power can be
expected for this test against ARCH/GARCH type dependence.

Instead of stating independence in terms of cumulative distribution
functions this can alternatively be expressed in terms of the
characteristic function. Székely, Rizzo, and Bakirov (2007) and Székely
and Rizzo (2009), as will be seen in Section \ref{distance}, make
systematic use of this in their introduction of the distance covariance
test. Two random variables \(X\) and \(Y\) are independent if and only
if the characteristic functions satisfy
\[\phi_{X,Y}(u,v) = \phi_{X}(u)\phi_{Y}(v), \quad \forall(u,v)\] where
\[ \phi_{X,Y}(u,v) = \E\left(e^{iuX+ivY}\right), \quad \phi_{X}(u) = \E\left(e^{iuX}\right), \quad \phi_{Y}(v) = \E\left(e^{ivY}\right).\]
This was exploited by Csörgö (1985) and Pinkse (1998) to construct tests
for independence based on the characteristic function in the iid and
time series case, respectively. Further work on testing of conditional
independence was done by Su and White (2007). See also Fan et al.
(2017). Hong (1999) put this into a much more general context by
replacing \(\rho_k(u,v)\) in \eqref{eq:Hong1} by
\[\sigma_k(u,v) = \phi_{X_{t},X_{t-|k|}}(u,v)-\phi_{X_{t}}(u)\phi_{X_{t-|k|}}(v).\]
By taking Fourier transform of this quantity one obtains

\begin{equation}
f(\omega,u,v) = \frac{1}{2\pi} \sum_{k= -\infty}^{\infty} \sigma_k(u,v)e^{-ik\omega}.
\label{eq:Hong2}
\end{equation}

Hong (1999) called \eqref{eq:Hong2} the generalized spectral density
function. Here \(f(\omega, u,v)\) can be estimated by
\[\widehat{f}_n(\omega,u,v) = \frac{1}{2\pi} \sum_{k=-n+1}^{n-1} (1-|k|/n)^{1/2}w(k/b)\widehat{\sigma}_k(u,v)e^{-ik\omega},\]
where \(w\) is a kernel weight function, \(b\) is a bandwidth or lag
order, and
\[\widehat{\sigma}_k(u,v) = \widehat{\phi}_k(u,v)-\widehat{\phi}_k(u,0)\widehat{\phi}_k(0,v)\]
with
\[\widehat{\phi}_k(u,v) = (n-|k|)^{-1}\sum_{t=|k|+1}^{n} e^{i(uX_t+vX_{t-|k|})}.\]
Under the null hypothesis of serial independence \(f(\omega,u,v)\)
becomes a constant function of frequency \(\omega\):
\[f_0(\omega,u,v) = \frac{1}{2\pi}\sigma_0(u,v)\] with
\(\sigma_0(u,v) = \phi(u+v)-\phi(u)\phi(v)\), where
\(\phi(\cdot) = \phi_{X_{t}}(\cdot)\). In order to test for independence
one can compare \(\widehat{f}_n(\omega,u,v)\) and
\(\widehat{f}_0(\omega,u,v)\) using e.g.~an \(L^2\)-functional. More
work related to this has been done by Hong and Lee (2003) and Escanciano
and Velasco (2006).

\subsection{Distance covariance}\label{distance}

We have seen that there are at least two ways of constructing
functionals that are consistent against all forms of dependence, namely
those based on the empirical distribution function initiated by
Hoeffding (1948) and briefly reviewed above, and those based on the
characteristic function represented by Csörgö (1985) in the iid case and
Pinkse (1998) in the serial dependence case, and continued in Hong
(1999; 2000) in a time series spectrum approach. Se also Fan et al.
(2017). Both Pinkse and Hong use a kernel type weight function in their
functionals. Thus, Pinkse uses a weight function \(g\) in the functional
\[ \int g(u)g(v) |\phi(u,v)|^2\di u \di v \] where for a pair of two
random variables \((X,Y)\)
\[\phi(u,v) = \E\left(e^{i(uX+vY)}\right)-\E\left(e^{iuX}\right)\E\left(e^{ivY}\right).\]
Let \(h(x) = \int e^{iux}g(u) \di u\). Pinkse in his simulation
experiments chose the weight functions \(h(x) = \exp(-\frac{1}{2}x^2)\)
and \(h(x) = 1/(1+x^2)\).

The authors of two remarkable papers, Székely, Rizzo, and Bakirov (2007)
and Székely and Rizzo (2009), take up the characteristic function test
statistic again in the non-time series case. But what distinguishes
these from earlier papers is an especially judicious choice of weight
function reducing the empirical characteristic function functional to
empirical moments of differences between the variables, or distances in
the vector case, this leading to covariance of distances. Some of these
ideas go back to what the authors term an ``energy statistic''; see
Székely (2002), Székely and Rizzo (2013) and also Székely and Rizzo
(2012). It has been extended to time series and multiple dependencies by
Davis et al. (2018), Fokianos and Pitsillou (2017), Zhou (2012), and
Dueck et al. (2014), Dueck, Edelman, and Richards (2015) and Yao, Zhang,
and Shao (2018). In the locally stationary time series case there is
even a theory, see Jentsch et al. (2018). The distance covariance, dcov,
seems to work well in a number of situations, and it has been used as a
yardstick by several authors writing on dependence and tests of
independence. In particular it has been used as a measure of comparison
in the work on local Gaussian correlation to be detailed in Section
\ref{lgc}. There are also points of contacts, as will be seen in Section
\ref{hsic}, with the HSIC measure of dependence popular in the machine
learning community.

The central ideas and derivations are more or less all present in
Székely, Rizzo, and Bakirov (2007). The framework is that of pairs of
iid vector variables \((X,Y)\) in \(\mathbb{R}^p\) and \(\mathbb{R}^q\),
respectively, and the task is to construct a test functional for
independence between \(X\) and \(Y\). Let
\(\phi_{X,Y}(u,v) = \E\left(e^{i(\langle {X,u \rangle + \langle Y, v} \rangle)}\right)\),
\(\phi_X(u) = \E\left(e^{i \langle X, u \rangle }\right)\) and
\(\phi_{Y}(v) = \E\left(e^{i \langle Y,v \rangle}\right)\) be the
characteristic functions involved, where
\(\langle \cdot , \cdot \rangle\) is the inner product in
\(\mathbb{R}^p\) and \(\mathbb{R}^q\), respectively. The starting point
is again the weighted characteristic functional

\begin{equation}
{\cal V}^{2}(X,Y;w) = \int_{\mathbb{R}^{p+q}} |\phi_{X,Y}(u,v)-\phi_{X}(u)\phi_{Y}(v)|^2w(u, v)\di u \di v,
\label{eq:Szekely1}
\end{equation}

where \(w\) is a weight function to be chosen, Note that it is easy to
choose \(w\) so that \({\cal V}^2(X,Y) = 0\) if and only if \(X\) and
\(Y\) are independent. Similarly, one defines

\begin{equation}
{\cal V}^2(X;w) = \int_{\mathbb{R}^{2p}} |\phi_{X,X}(u,v)-\phi_{X}(u)\phi_{X}(v)|^2w(u,v)\di u \di v
\label{eq:Szekely1a}
\end{equation}

and \({\cal V}^2(Y;w)\). The distance correlation, dcor, is next defined
by, assuming \({\cal V}^{2}(X){\cal V}^{2}(Y) > 0\),
\[{\cal R}^2(X,Y) = \frac{{\cal V}^{2}(X,Y)}{\sqrt{{\cal V}^{2}(X){\cal V}^{2}(Y)}}.\]
These quantities can be estimated by the empirical counterparts given
\(n\) observations of the vector pair \((X,Y)\) with

\begin{equation}
{\cal V}_n^{2}(X,Y;w) =  \int_{\mathbb{R}^{p+q}} |\phi_{X,Y}^{n}(u,v)-\phi_{X}^{n}(u)\phi_{Y}^{n}(v)|^2w(u,v)\di u \di v,
\label{eq:Szekely1b}
\end{equation}

where, for a set of observations \(\{(X_1,Y_1),\ldots,(X_n,Y_n)\}\) the
empirical characteristic functions are given by
\[\phi_{X,Y}^{n}(u,v) = \frac{1}{n} \sum_{k=1}^{n}\exp\{i(\langle X_k,u \rangle + \langle Y_k,v \rangle)\}\]
and \[
\phi_{X}^{n}(u) = \frac{1}{n} \sum_{k=1}^{n} \exp \{i \langle X_k,u \rangle\}, \quad  \phi_{Y}^{n}(v) = \frac{1}{n} \sum_{k=1}^{n} \exp \{i \langle Y_k, v \rangle \}.\]
It turns out that it is easier to handle the weight function in the
framework of the empirical characteristic functions. It will be seen
below that

\begin{equation}
w(u,v) = (c_pc_q |u|_p^{1+p}|v|_q^{1+q})^{-1}
\label{eq:Szekely2}
\end{equation}

is a good choice. Here \(|\cdot|_p\) is the Euclidean norm in
\(\mathbb{R}^{p}\) and similarly for \(|\cdot|_q\). Moreover, the
normalizing constants are given by
\(c_j = \pi^{(1+j)/2}/\Gamma((1+j)/2)\), \(j=p,q\). For it to make sense
to introduce the weight function on the empirical characteristic
function one must show that the empirical functionals \({\cal V}_n\)
converges to the theoretical functionals \({\cal V}\) for this weight
function. This is not trivial because of the singularity at \(0\) for
\(w\) given by \eqref{eq:Szekely2}. A detailed argument is given in the
proof of Theorem 2 in Székely, Rizzo, and Bakirov (2007).

The advantage of introducing the weight function for the empirical
characteristic functions is that one can compute the squares in
\eqref{eq:Szekely1b} and then interchange summation and integration. The
resulting integrals can be computed using trigonometric identities, in
particular the odd symmetry of products of cosines and sines which makes
corresponding integrals disappear. The details are given in the proof of
Theorem 1 in Székely, Rizzo, and Bakirov (2007) and in Lemma 1 of the
Appendix of Szekely and Rizzo (2005) who in turn refer to Prudnikov,
Brychkov, and Marichev (1986) for the fundamental lemma \[
\int_{\mathbb{R}^{d}} \frac{1-\cos \langle x,u \rangle}{|u|_d^ {u+\alpha}} \di u = C(d,\alpha)|x|_d^{\alpha}
\] for \(0 < \alpha < 2\) with

\begin{equation}
C(d, \alpha) = \frac{2\pi^{d/2} \Gamma(1-\alpha/2)}{\alpha 2^{\alpha}\Gamma((d+\alpha)/2)},
\label{eq:Szekely}
\end{equation}

and where the weight function considered above corresponds to
\(\alpha = 1\) and \(d=p\) or \(d=q\) in \eqref{eq:Szekely2}. The general
\(\alpha\)-case corresponds to a weight function
\[w(u,v;\alpha) = (C(p,\alpha)C(q, \alpha) |u|_p^{p+\alpha}|v|_q^{q+\alpha})^{-1}.\]
With the simplification \(\alpha = 1\) all of this implies that
\({\cal V}_n^{2}\) as defined in \eqref{eq:Szekely1b}, can be computed as
\[{\cal V}_n^2(u,v) = S_1+S_2-2S_3\] where
\[S_1 = \frac{1}{n^2}\sum_{k,l=1}^{n}|X_k-X_l|_p|Y_k-Y_l|_q,\] \[
S_2 = \frac{1}{n^2}\sum_{k,l = 1}^{n}|X_k-X_l|_p \frac{1}{n^2}\sum_{k,l=1}^{n}|Y_k- Y_l|_q,
\]

\begin{equation}
S_3 = \frac{1}{n^3}\sum_{k=1}^{n}\sum_{l,m=1}^{n}|{\bf X}_k-{\bf X}_l|_p|{\bf Y}_k-{\bf Y}_m|_q
\label{eq:Szekely4}
\end{equation}

which explains the appellation distance covariance. In fact, it is
possible to further simplify this by introducing
\[a_{kl} = |X_k-X_l|_p, \quad \overline{a}_{k.} = \frac{1}{n}\sum_{l=1}^{n}a_{kl}, \quad \overline{a}_{.l} = \frac{1}{n}\sum_{k=1}^{n}a_{kl},\]
\[\overline{a}_{..} = \frac{1}{n^2}\sum_{k,l = 1}^{n}a_{kl}, \quad A_{kl} = a_{l}-\overline{a}_{k.}-\overline{a}_{.l}+a_{..},\]
for \(k,l = 1,\ldots,n\). Similarly, one can define
\(b_{kl}=|Y_k- Y_l|_q\) and
\(B_{kl}=b_{kl}-\overline{b}_{k.}-\overline{b}_{.l}+\overline{b}_{..}\)
and \[{\cal V}_n^{2}(X,Y) = \frac{1}{n^2} \sum_{k,l=1}^{n}A_{kl}B_{kl}\]
and
\[{\cal V}_n^2(X) = {\cal V}_n^{2}(X,X) = \frac{1}{n^2} \sum_{k,l = 1}^{n} A_{kl}^{2}\]
and similarly for \({\cal V}_n^{2}(Y)\). From this one can easily
compute \({\cal R}_n^{2}(X,Y)\). The computations are available in an
R-package: Rizzo and Szekely (2018).

As is the case of the empirical joint distribution functional it can be
expected that the curse of dimensionality will influence the result for
large and moderate values of \(p\) and \(q\). Obviously, in the time
series case, it is possible to base oneself on pairwise distances as in
\eqref{eq:Box-test0} or \eqref{eq:Box-test1}, which has been done in Yao,
Zhang, and Shao (2018).

Letting \(n \to \infty\), it is not difficult to prove that an
alternative expression for \({\cal V}({\bf X},{\bf Y})\) is given by
(assuming \(\E|X|_p < \infty\) and \(\E|Y|_q < \infty\))

\begin{align}
{\cal V}^{2}(X,Y) &= \E_{X,X',Y,Y'}\{|X-X'|_p|Y-Y'|_q\}\nonumber \\ &+\E_{X,X'}\{|X-X'|_p\}\E_{Y,Y'}\{|Y-Y'|_q\} \nonumber \\ &-2\E_{X,Y}\{\E_{X'}|X-X'|_p\E_{Y'}|Y-Y'|_q\}
\label{eq:alternative}
\end{align}

where \((X, Y)\), \((X', Y')\) are iid. This expression will be useful
later in Section \ref{hsic} in a comparison with the HSIC statistic.
Properly scaled \({\cal V}_n^2\) has a limiting behavior under
independence somewhat similar to that of \eqref{eq:Hoeffding2} in Section
\ref{cumulative}. Namely, under the condition of existence of first
moment, \(n{\cal V}_n^2/S_2\), converges in distribution to a quadratic
form
\[Q \buildrel {{\cal L}} \over {\to} \sum_{j=1}^{\infty} \lambda_j Z_j^2\]
where \(\{Z_j\}\) are independent standard normal variables,
\(\{\lambda_j\}\) are non-negative constants that depend on the
distribution of \((X,Y)\). One can also obtain an empirical process
limit theorem, Theorem 5 of Székely, Rizzo, and Bakirov (2007). In the
R-package, as for the case of the empirical distribution function, it
has been found advantageous to rely on re-sampling via permutations.
This is quite fast since the algebraic formulas \eqref{eq:Szekely4} are
especially amenable to permutations. Both Székely, Rizzo, and Bakirov
(2007) and Székely and Rizzo (2009) in their experiments only treat the
case of \(\alpha = 1\) in \eqref{eq:Szekely}.

Turning to the properties (i) - (vii) of Rényi (1959) listed in the
beginning of this section, it is clear that (i) - (iv) are satisfied by
\({\cal R}\). Moreover, according to Székely, Rizzo, and Bakirov (2007),
if \({\cal R}_n(x;y) = 1\), then there exists a vector \(\alpha\), a
non-zero real number \(\beta\) and an orthogonal matrix \(C\) such that
\(Y = \alpha + \beta X C\), which is not quite the same as Rényi's
requirement (v). Also, the general invariance in his property (vi) does
not seem to hold. The final criterion (vii) of Rényi is that the
dependent measure should reduce to the absolute value of Pearson's
\(\rho\) in the bivariate normal case. This is not quite the case for
the dcov, but it comes close, as is seen from Theorem 6 of Székely and
Rizzo (2009). In fact, if \((X,Y)\) is bivariate normal with
\(\E(X)=\E(Y)=0\) and \(\Var(X) = \Var(Y)=1\) and with correlation
\(\rho\), then \({\cal R}(X,Y) \leq |\rho|\) and
\[\inf_{\rho \neq 0} \frac{{\cal R}(X,Y)}{|\rho|} = \lim_{\rho \to 0} \frac{
{\cal R}(X,Y)}{|\rho|} = \frac{1}{2(1+\pi/3-\sqrt{3})^{1/2}} \approx 0.891.\]

\subsection{The HSIC measure of dependence}\label{hsic}

Recall the definition and formula for the maximal correlation. This, as
stated in Section \ref{maximal} gives rise to a statistic \(S(X,Y)\),
where \(S(X,Y) = 0\) if and only if \(X\) and \(Y\) are independent. But
it is difficult to compute since it requires the supremum of the
correlation \(\rho(f(X),g(Y))\) taken over Borel-measurable \(f\) and
\(g\). In the framework of reproducing kernel Hilbert spaces (RKHS) it
is possible to pose this problem, or an analogous one, much more
generally, and one can compute an analogue of \(S\) quite easily. This
is the so-called HSIC (Hilbert-Schmidt Independence Criterion).

Reproducing kernel Hilbert spaces are very important tools in
mathematics as well as in statistics. A general reference to
applications in statistics is Berlinet and Thomas-Agnan (2004). In the
last decade or so there has also been a number of uses of RKHS in
dependence modeling. These have often been published in the machine
learning literature. See e.g Gretton, Herbrich, et al. (2005), Gretton
and Györfi (2010), Gretton and Györfi (2012) and Sejdinovic et al.
(2013).

We have found the quite early paper by Gretton, Bousquet, et al. (2005)
useful both for a glimpse of the general theory and for the HSIC
criterion in particular.

A reproducing kernel Hilbert space is a separable Hilbert space
\({\cal F}\) of functions \(f\) on a set \({\cal X}\), such that the
evaluation functional \(f \to f(x)\) is a continuous linear functional
on \({\cal F}\) for every \(x \in {\cal X}\). Then, from the Riesz
representation theorem, Muscat (2014), chapter 10, there exists an
element \(k_x \in {\cal F}\) such that \(\langle f,k_x \rangle = f(x)\),
where \(\langle \cdot, \cdot \rangle\) is the inner product in
\({\cal F}\). Applying this to \(f = k_x\) and another point
\(y \in {\cal X}\) we have \(\langle k_x, k_y \rangle = k_x(y)\). The
function \((x,y) \to k_x(y)\) from \(\cal{X} \times \cal{X}\) to
\(\mathbb{R}\) is the kernel of the RKHS \({\cal F}\). It is symmetric
and positive definite because of the symmetry and positive definiteness
of the inner product in \({\cal F}\). We use the notation \(k(x,y)\) for
the kernel.

The next step is to introduce another set \({\cal Y}\) with a
corresponding RKHS \({\cal G}\) and to introduce a probability structure
and probability measures \(p_X\), \(p_Y\) and \(p_{X,Y}\) on
\({\cal X}\), \({\cal Y}\) and \({\cal X} \times {\cal Y}\),
respectively. With these probability measures and function spaces
\({\cal F}\) and \({\cal G}\) one can introduce correlation of functions
of stochastic variables on \({\cal X}\), \({\cal Y}\) and
\({\cal X} \times {\cal Y}\). This is an analogy of the functions used
in the definition of the maximal correlation. In RKHS setting the
covariance (or cross covariance) is an \emph{operator} on the function
space \({\cal F}\). Note also that this has a clear analogy in
functional statistics, see e.g. Ferraty and Vieu (2006).

It is time to introduce the Hilbert-Schmidt operator: A linear operator
\(C: {\cal G} \to {\cal F}\) is called a Hilbert-Schmidt operator if its
Hilbert-Schmidt (HS) norm \(||C||_{HS}\)
\[||C||_{HS}^2 \doteq \sum_{i,j} \langle Cv_j,u_i \rangle_{{\cal F}} < \infty\]
where \(u_i\) and \(v_j\) are orthonormal bases of \({\cal F}\) and
\({\cal G}\), respectively. The HS-norm generalizes the Froebenius norm
\(||A||_{F} = (\sum_i \sum_j a_{ij}^2)^{1/2}\) for a matrix
\(A = (a_{ij})\). Finally, we need to define the tensor product in this
context: If \(f \in {\cal F}\) and \(g \in {\cal G}\), then the tensor
product operator \(f \otimes g: {\cal G} \to {\cal F}\) is defined by
\[(f \otimes g)h \doteq f\langle g , h \rangle_{{\cal G}}, \quad h \in {\cal G}.\]
Moreover, by using the definition of the HS norm it is not difficult to
show that
\[||f \otimes g||_{HS}^2 = ||f||_{{\cal F}}^2||g||_{{\cal G}}^2.\] We
can now introduce an expectation and a covariance on these function
spaces. Again, the analogy with corresponding quantities in functional
statistics will be clear. We assume that \(({\cal X},\Gamma)\) and
\(({\cal Y},\Lambda)\) are furnished with probability measures
\(p_X, p_Y\), and with \(\Gamma\) and \(\Lambda\) being
\(\sigma\)-algebras of sets on \({\cal X}\) and \({\cal Y}\). The
expectations \(\mu_X \in {\cal F}\) and \(\mu_Y \in {\cal G}\) are
defined by, \(X\) and \(Y\) are stochastic variables in
\(({\cal X}, \Gamma)\) and \(({\cal Y}, \Lambda)\), respectively,
\[\langle \mu_X,f \rangle_{{\cal F}} = \E_X[f(X)]\] and
\[\langle \mu_Y, g \rangle_{{\cal G}} = \E_Y[g(Y)]\] where \(\mu_X\) and
\(\mu_Y\) are well-defined as elements in \({\cal F}\) and \({\cal G}\)
because of the Riesz representation theorem. The norm is obtained by
\[||\mu_X||_{{\cal F}}^2 = \E_{X,X'}[k(X,X')],\] where as before \(X\)
and \(X'\) are independent but have the same distribution \(p_X\), and
where \(||\mu_Y||\) is defined in the same way. With given
\(\phi \in {\cal F}\), \(\psi \in {\cal G}\) we can now define the cross
covariance operator as
\[C_{X,Y} \doteq \E_{X,Y}[(\phi(X)-\mu_X)\otimes (\psi(Y)-\mu_Y)] = \E_{X,Y}[\phi(X) \otimes \phi(Y)]-\mu_X \otimes \mu_Y.\]
Now, take \(\phi(X)\) to be identified with \(k_X \in{\cal F}\) defined
above as a result of the Riesz representation theorem, and
\(\psi(Y) \in{\cal G}\) defined in exactly the same way. The
Hilbert-Schmidt Information Criterion (HSIC) is then defined as the
squared HS norm of the associated cross-covariance operator \[
\textrm{HSIC}(p_{XY},{\cal F},{\cal G}) \doteq ||C_{XY}||_{HS}^2.
\] Let \(k(x,x')\) and \(l(y,y')\) be kernel functions on \({\cal F}\)
and \({\cal G}\). Then (Gretton, Bousquet, et al. 2005, Lemma 1), the
HSIC criterion can be written in terms of these kernels as

\begin{align}
& \textrm{HSIC}(p_{XY},{\cal F},{\cal G}) \nonumber \\
&= \E_{X,X',Y,Y'}[k(X,X')l(Y,Y')]+\E_{X,X'}[k(X,X')]\E_{Y,Y'}[l(Y,Y')] \nonumber \\
&\qquad\qquad\qquad\qquad\qquad\qquad\qquad-2\E_{X,Y}\{\E_{X'}[k(X,X')]\E_{Y'}[l(Y,Y')]\}
\label{eq:HSIC0}
\end{align}

Existence is guaranteed if the kernels are bounded. The similarity in
structure to \eqref{eq:alternative} for the distance covariance should be
noted (partly due to the identity \((a-b)^2 = a^2+b^2-2ab\) but going
deeper as will be seen below when HSIC is compared to dcov). Note that
the kernel functions depend on the way the spaces \({\cal F}\) and
\({\cal G}\) and their inner products are defined. In fact it follows
from a famous result by Moore-Aronszajn, see Aronszajn (1950), that if
\(k\) is a symmetric, positive definite kernel on a set \({\cal X}\),
then there is a unique Hilbert space of functions on \({\cal X}\) for
which \(k\) is a reproducing kernel. Hence as will be seen next, in
practice when applying the HSIC criterion, the user has to choose a
kernel.

With some restrictions the HSIC measure is a proper measure of
dependence in the sense of the Rényi (1959) criterion (iv): From Theorem
4 of Gretton, Bousquet, et al. (2005) one has that if the kernels \(k\)
and \(l\) are universal (universal kernel is a mild continuity
requirement on the kernel) on compact domains \({\cal X}\) and
\({\cal Y}\), then \(||C_{XY}||_{HS} = 0\) if and only if \(X\) and
\(Y\) are independent. The compactness assumption results from the
application of an equality for bounded random variables taken from
Hoeffding (1963), that is being used actively in the proof.

A big asset of the HSIC measure it that its empirical version is easily
computable. In fact if we have independent observations
\(X_1,\ldots,X_n\) and independent observations \(Y_1,\ldots,Y_n\), then

\begin{equation}
\textrm{HSIC}_n(X,Y,{\cal F},{\cal G}) = (n-1)^{-2}\textrm{tr}\{KHLH\}
\label{eq:HSIC1}
\end{equation}

where tr is the trace operator and the \(n \times n\) matrices
\(H, K, L\) are defined by \[
K = \{K_{ij}\} = \{k(X_i,X_j)\}, \;\;L = \{L_{ij}\} = \{l(Y_i,Y_j)\},\;\;
H = \{H_{ij}\} = \{\delta_{ij}-n^{-1}\},
\] where \(\delta_{ij}\) is the Kronecker delta. It is shown in Gretton,
Bousquet, et al. (2005) that this estimator converges towards
\(||C_{XY}||_{HS}^2\). The convergence rate is \(n^{-1/2}\). There is
also a limit theorem for the asymptotic distribution, which under the
null hypothesis of independence and scaled with \(n\) converges in
distribution to the random variable
\(Q = \sum_{i,j = 1}^{\infty} \lambda_i\eta_j N_{ij}^2\), where the
\(N_{ij}\) are independent standard normal variables, and \(\lambda _i\)
and \(\eta_j\) are eigenvalues of integral operators associated with
centralized kernels derived from \(k\) and \(l\) and integrating using
the probability measures \(p_X\) and \(p_Y\), respectively. Again, this
should be compared to the limiting variable for the statistic in the
Cramér- von Mises functional \eqref{eq:Hoeffding2}. Critical values can be
obtained for \(Q\), but as a rule one seems to rely more on resampling
as is the case for most independence test functionals.

It is seen from \eqref{eq:HSIC1} that computation of the empirical HSIC
criterion requires the evaluation of \(k(X_i,X_j)\) and \(l(Y_i,Y_j)\).
Then appropriate kernels have to be chosen. Two commonly used kernels
are the Gaussian kernel given by
\[k(x,y) = e^{\frac{|x-y|^{2}}{2\sigma^2}}, \quad \sigma > 0\] and the
Laplace kernel \[k(x,y) = e^{\frac{|x-y|}{\sigma}}, \quad \sigma > 0.\]
Pfister and Peters (2017) describe a recent R-package involving HSIC.
Gretton, Bousquet, et al. (2005) use these kernels in comparing the HSIC
test with several other tests, including the dcov test in, among other
cases, an independent component setting. Both of these tests do well,
and one of these tests does not decisively out-compete the other one.
This is perhaps not so unexpected because there is a strong relationship
between these two tests. This is demonstrated by Sejdinovic et al.
(2013). They look at both the dcov test and the HSIC test in a
generalized setting of semimetric spaces, i.e with kernels and distances
defined on such spaces \({\cal X}\) and \({\cal Y}\). For a given
distance function they introduce a distance-induced kernel, and under
certain regularity conditions they establish a relationship between
these two quantities. There is a related paper by Lyons (2013) which
obtains similar results but not in an explicit RKHS context, in fact in
a general dcov context.

Let \(\rho_{\cal X}\) and \(\rho_{\cal Y}\) be distance measures on the
semi-metric spaces \({\cal X}\) and \({\cal Y}\), respectively. Then a
generalized dcov distance functional can be defined as, compare again to
\eqref{eq:alternative},

\begin{align*}
& \cal{V}^{\textrm{2}}_{\rho_{\cal{X}},\rho_{\cal{Y}}}(X,Y) = \\ & \qquad  \E_{XY}\E_{X'Y'}\rho_{\cal{X}}(X,X')\rho_{\cal{Y}}(Y,Y')+\E_X\E_{X'}\rho_{\cal{X}}(X,X')\E_Y\E_{Y'}\rho_{\cal{Y}}(Y,Y') \\
&\qquad\qquad\qquad\qquad\qquad\qquad\qquad\qquad -2\E_{XY}\{\E_{X'}\rho_{\cal{X}}(X,X')\E_{Y'}\rho_{\cal{Y}}(Y,Y')\}.
\end{align*}

This distance covariance in metric spaces characterizes independence,
that is, \({\cal V}_{\rho_{\cal X},\rho_{\cal Y}}(X,Y) = 0\) if and only
if \(X\) and \(Y\) are independent, and if the metrics \(\rho_{\cal X}\)
and \(\rho_{\cal Y}\) satisfy an additional property, termed strong
negative type. See Sejdinovic et al. (2013) for more details. An asset
of the RKHS formulation is that it is very general. As was seen from the
introduction of HSIC above, the sets \({\cal X}\) and \({\cal Y}\) can
have a metric space structure, and probability measures \(p_X\), \(p_Y\)
and \(p_{X,Y}\) can still be introduced, and the definition of HSIC
given in the beginning of this section and the accompanying
decomposition \eqref{eq:HSIC0} still make sense in this generalized
framework. It can then be shown that, Theorem 24 in Sejdinovic et al.
(2013), one has the following equivalence: Let \(k_{\cal X}\) and
\(k_{\cal Y}\) be any two kernels on \({\cal X}\) and \({\cal Y}\) that
generate \(\rho_{\cal X}\) and \(\rho_{\cal Y}\), respectively, and let
\(k((x,y),(x',y')) = k_{\cal X}(x,x')k_{\cal Y}(y,y')\), then
\({\cal V}_{\rho_{\cal X},\rho_{\cal Y}}^2 = 4\textrm{HSIC}^{2}(p_{XY}, {\cal F},{\cal G})\).
Among the regularity conditions required for this result is the
assumption of ``negative type'', which is satisfied in standard
Euclidean spaces.

However, it is not possible to find a direct RKHS representation of the
characteristic function representation \eqref{eq:Szekely1a} of
\({\cal V}^{2}\).

Lately there have been other extensions of both the dcov and HSIC to
conditional dependence, partial distance and to time series. A few
references are Szekely and Rizzo (2014), Chwialkowski and Gretton
(2014), Zhang et al. (2012) and Pfister et al. (2018). A recent tutorial
on RKHS is Gretton (2017).

\subsection{Density based tests of independence}\label{density}

Intuitively, one might think that knowing that the density exists should
lead to increased power of the independence tests due to more
information. This is true, at least for some examples (see e.g.
Teräsvirta, Tjøstheim, and Granger (2010), Chapter 7.7). As in the
preceding sections one can construct distance functionals between the
joint density under dependence and the product density under
independence. A number of authors have considered such an approach; both
in the iid and time series case, see e.g. Rosenblatt (1975), Robinson
(1991), Skaug and Tjøstheim (1993b; 1996), Granger, Maasoumi, and Racine
(2004), Hong and White (2005), Su and White (2007) and Berrett and
Samworth (2017). For two random variables \(X\) and \(Y\) having joint
density \(f_{X,Y}\) and marginals \(f_X\) and \(f_Y\) the degree of
dependence can be measured by \(\Delta(f_{X,Y},f_Xf_Y)\), where
\(\Delta\) is now the distance measure between two bivariate density
functions. The variables are normalized with \(\E(X)=\E(Y)=0\) and
\(\Var(X)=\Var(Y)=1\). It is natural to consider the Rényi (1959)
requirements again, in particular the requirements (iv) and (vi).

All of the distance functionals considered will be of type

\begin{equation}
\Delta = \int B\{f_{X,Y}(x,y),f_X(x),f_Y(y)\}f_{X,Y}(x,y)\di x \di y
\label{eq:density0}
\end{equation}

where \(B\) is a real-valued function such that the integral exists. If
\(B\) is of the form \(B(z_1,z_2,z_3) = D(z_1/z_2z_3)\), we have

\begin{equation}
\Delta = \int D\left\{\frac{f_{X,Y}(x,y)}{f_X(x)f_Y(y)}\right\}f_{X,Y}(x,y) \di x \di y
\label{eq:density1}
\end{equation}

which by the change of variable formula for integrals is seen to have
the Rényi property (vi). Moreover, if \(D(w) \geq 0\) and \(D(w) = 0\)
if and only if \(w=1\), then Rényi property (iv) is fulfilled. Several
well-known distance measures for density functions are of this type. For
instance, letting \(D(w) = 2(1-w^{-1/2})\) we obtain the Hellinger
distance

\begin{align*}
H & = \int \left\{\sqrt{f_{X,Y}(x,y)}-\sqrt{f_X(x)f_Y(y)}\right\}^2 \di x \di y \\
  & = 2 \int \left\{1-\sqrt{\frac{f_X(x)f_Y(y)}{f_{X,Y}(x,y)}}\right\}f_{X,Y}(x,y)\di x \di y
\end{align*}

between \(f_{X,Y}\) and \(f_Xf_Y\). The Hellinger distance is a metric
and hence satisfies the Rényi property (iv).

Chung et al. (1989) defined the so-called directed divergence of degree
\(\gamma\) \((0 < \gamma < 1)\), which is also related to the Rényi
divergence, Rényi (1961),

\begin{equation}
\Delta_{\gamma}(f_{X,Y},f_Xf_Y) = \frac{1}{1-\gamma} \int \left[1-\left\{\frac{f_X(x)f_Y(y)}{f_{X,Y}(x,y)}\right\}^{\gamma}\right]f_{X,Y}(x,y)\di x \di y.
\label{eq:density2}
\end{equation}

It is seen that the Hellinger distance is a special case
\((\gamma = 1/2)\). Clearly, the measure \(\Delta_{\gamma}\) satisfies
(vi), and for \(0 < \gamma < 1\), using Hölder's inequality,
\[(\gamma - 1)\Delta_{\gamma}(f_{X,Y},f_Xf_Y) = \int \{f_X(x)f_Y(y)\}^{\gamma}\{f_{X,Y}(x,y)\}^{1-\gamma}\di x \di y-1 \leq 0\]
with equality if and only if \(f_{X,Y} = f_Xf_Y\). Hence,
\(\Delta_{\gamma}\) satisfies (iv) for \(0 < \gamma < 1\).

The familiar Kullback-Leibler information (entropy) distance is obtained
as a limiting case as \(\gamma \to 1\),

\begin{equation}
I= \int \ln \left\{\frac{f_{X,Y}(x,y)}{f_X(x)f_Y(y)}\right\}f_{X,Y}(x,y)\di x \di y.
\label{eq:Kullback}
\end{equation}

Since this distance is of type \eqref{eq:density1}, it satisfies (vi), and
it can also be shown to satisfy (iv). Going outside the range
\(0 < \gamma < 1\), for \(\gamma = 2\) in \eqref{eq:density2}, the
test-of-fit distance in Bickel and Rosenblatt (1973) emerges. See also
Granger and Lin (1994). A very recent paper linking \(I\) with other
recent approaches to independence testing is Berrett and Samworth
(2017).

All of the above measures are trivially extended to two arbitrary
multivariate densities. However, estimating such densities in high or
moderate dimensions may be difficult due to the curse of dimensionality.
A functional built up from pairwise dependencies can be considered
instead such as in \eqref{eq:Box-test0} and \eqref{eq:Box-test1}.

For a given functional \(\Delta = \Delta(f,g)\) depending on two
densities \(f\) and \(g\), \(\Delta\) may be estimated by
\(\widehat{\Delta} = \Delta(\widehat{f},\widehat{g})\). There are
several ways of estimating the densities, e.g.~the kernel density
estimator, \[\widehat{f}_X(x) = \frac{1}{n}\sum_{i=1}^{n}K_b(x-X_i)\]
for given observations \(\{X_1,\ldots,X_n\}\). Here
\(K_b(x-X_i) = b^{-p}K\{b^{-1}(x-X_i)\}\), where \(b\) is the bandwidth
(generally a matrix), \(K\) is the kernel function, and \(p\) is the
dimension of \(X_i\). The kernel function is usually taken to be a
product of one-dimensional kernels; i.e., \(K(x) = \prod K_i(x_i)\),
where each \(K_i\) generally is non-negative and satisfies
\[\int K_i(u) \di u = 1,\quad \int u^2K_i(u) \di u < \infty.\]

Once estimators for \(f_{X,Y}\), \(f_X\) and \(f_Y\) in the integral
expression \eqref{eq:density0} for \(\Delta\) have been obtained, the
integral can be computed by numerical integration or by empirical
averages using the ergodic theorem (or law of large numbers in the iid
case). Consequently for a given lag \(k\) in the time series case,
\[\widehat{\Delta}_k = \frac{1}{n-k}\sum_{t=k+1}^{n} B\{\widehat{f}_k(X_t,X_{t-k}),\widehat{f}(X_t),\widehat{f}(X_{t-k})\}w(X_t,X_{t-k}).\]
Here \(f_k\) is the joint density of \((X_t,X_{t-k})\), and \(w\) is a
weight function, e.g.,
\(w(u,v) = 1\{|u| \leq c \sigma_X\}1\{|v| \leq c \sigma_X\}\) for some
chosen constant \(c\).

Under regularity conditions (see e.g. Skaug and Tjøstheim (1996)\},
consistency and asymptotic normality can be obtained for the estimated
test functionals. It should be noted that the leading term in an
asymptotic expansion of the standard deviation of \(\widehat{\Delta}\)
for the estimated Kullback-Leibler functional \(\widehat{I}_{k,w}\) and
the estimated Hellinger functional \(\widehat{H}_{k,w}\) is of order
\(O(n^{-1/2})\). This is of course the same as for the standard
deviation of a parametric estimate in a parametric estimation problem.
In that situation the next term of the Edgeworth expansion is of order
\(O(n^{-1})\), and for moderately large values of \(n\) the first order
term \(n^{-1/2}\) will dominate. However, for the functionals considered
above, due to the presence of an \(n\)-dependent bandwidth, the next
terms in the Edgeworth expansion are much closer, being of order
\(O(n^{-1/2}b)\) and \(O(\{nb\}^{-1})\), and since typically
\(b=O(n^{-1/6})\) or \(O(n^{-1/5})\), \(n\) must be very large indeed to
have the first term dominate in the asymptotic expansion. As a
consequence, first order asymptotics in terms of the normal
approximation cannot be expected to work well unless \(n\) is
exceedingly large. Hence, basing a test of independence directly on the
asymptotic theory may be hazardous as the real test size will typically
deviate substantially from the nominal size. In this sense the situation
is quite different from the empirical functionals treated in the
previous sections, where there is no bandwidth parameter involved.

All of this suggests the use of the bootstrap or permutations as an
alternative for constructing the null distribution. One may anticipate
that it picks up higher-order terms of the Edgeworth expansion (Hall
1992, Chapters 3 and 4), although no rigorous analysis to confirm this
has been carried out for the functionals discussed here.

The fact that the permutation test yields an exact size, and that
resampling tests generally perform much better, underscores a major
point. In tests involving a bandwidth it is absolutely essential to use
resampling in practice. The asymptotic theory is far too inaccurate
except possibly in cases where the sample size is extremely large. In
the empirical functional case treated in Sections \ref{cumulative} and
\ref{distance}, the asymptotic theory is incomparably more accurate, but
even in this case the experience so far seems to be that resampling
generally does slightly better.

It is quite difficult to undertake local asymptotic power analyses for
the functionals based on estimated density functions. For reasons
mentioned above, asymptotic studies could be expected to be unreliable
unless \(n\) is very large. It has therefore been found more useful to
carry out comparative simulation studies against a wide choice of
alternatives using a modest sample size, see e.g.~Skaug and Tjøstheim
(1993; 1996) and Hong and White (2005). These references also contain
applications to real data.

\subsection{Test functionals generated by local dependence
relationships}\label{global-local}

If one has bivariate normal data with standard normal marginals and
\(\rho = 0\) one gets observations scattered in a close to circular
region around zero, and most test functionals will easily recognize this
as a situation of independence. However, as pointed out by Heller,
Heller, and Gorfine (2013), if data are generated along a circle, e.g.
\(X^2+Y^2 = N\) for some stochastic noise variable \(N\), then \(X\) and
\(Y\) are dependent, but the dcov, and undoubtedly other test
functionals, among which \(\rho\)-based tests, fail. They give some
other examples as well of similar failures for test functionals to
detect symmetric geometric patterns. Heller, Heller, and Gorfine (2013)
point a way out of this difficulty, namely by looking at dependence
locally (along the circle) and then aggregate the dependence by
integrating or by other means over the local regions. There are of
course several ways of measuring local dependence and we will approach
this problem more fundamentally in Section \ref{local}.

Before coming to the paper by Heller, Heller, and Gorfine (2013) and
papers using similar approaches, partly for historic reasons, we look
briefly at the correlation integral and the so-called BDS test named
after its originators Brock, Dechert and Scheinkman. This test has a
local flavor at its basis, but the philosophy is a bit different from
the other tests presented in the current subsection. The BDS test
attracted much attention among econometricians in the 1990s, and it has
since been improved.

The starting point is the correlation integral introduced in Grassberger
and Procaccia (1983) as a means of measuring fractal dimension of
deterministic data. It measures serial dependence patterns in the sense
that it keeps track of the frequency with which temporal patterns are
repeated in a data sequence. Let \(\{x_1,\ldots ,x_n\}\) be a sequence
of numbers and let
\[x_t^{(k)} = [x_t,x_{t-1},\ldots,x_{t-k+1}], \quad k \leq t \leq n.\]
Then the correlation integral for embedding dimension \(k\) is given by
\[C_{k,n}(\varepsilon) = \frac{2}{n(n-1)} \sum_{1 \leq s \leq t \leq n}1\left(||x_t^{(k)}-x_{s}^{(k)}|| < \varepsilon\right).\]
Here, for \(x=[x_1,\ldots,x_k]\),
\(||x|| = \max_{1 \leq i \leq k}|x_i|\), where \(1(\cdot)\) is the
indicator function and \(\varepsilon > 0\) is a cut off threshold which
could be a multiple of the standard deviation in the case of a
stationary process. The parameter \(\varepsilon\) may also be considered
to be a tuning parameter. Let
\[C_k(\varepsilon) = \lim_{n \to \infty}C_{k,n}(\varepsilon).\] If
\(\{X_t\}\) is an absolutely regular (Bradley, 1986, p.169) stationary
process the above limit exists and is given by
\[C_k(\varepsilon) = \int1(||x_1^{(k)}-x_2^{(k)}|| < \varepsilon)\di F_k(x_1^{(k)})\di F_k(x_2^{(k)}),\]
where \(F_k\) is the joint cumulative distribution function of
\(X_t^{(k)}\). Since
\[1(||x_1^{(k)}-x_2^{(k)}|| < \varepsilon) = \prod_{i=1}^{k} 1(|x_{1i}-x_{2i}| < \varepsilon),\]
it is easily seen that if \(\{X_t\}\) consists of iid random variables,
then \[C_k(\varepsilon) = \{C_1(\varepsilon)\}^{k}.\]

This theory is the starting point for the BDS test. Under the hypothesis
of independence, and excluding the case of uniformly distributed random
variables, Broock et al. (1996) have established asymptotic normality
\({\cal N}(0,1)\) of
\(\sqrt{n}[C_{k,n}(\varepsilon)-\{C_{1,n}(\varepsilon)\}^k]/V_{k,n}\)
under an appropriate scaling factor \(V_{k,n}\). As mentioned, the test
has found considerable use among econometricians, but it suffers from
some limitations, the arbitrariness of the choice of \(\varepsilon\),
the probability of rejecting independence does not always approach 1 as
\(n \to \infty\), and finally, and probably most critically, the
convergence to the asymptotic normal distribution may be very slow,
making bootstrapping a possible alternative. These problems were pointed
out by Genest, Ghoudi, and Rémillard (2007) who proposed a rank based
extension, where these difficulties are to a large degree eliminated.
Because the limiting distribution of the rank-based test statistics is
margin-free, their finite-sample p-values can be easily calculated by
simulations.

Next, returning to the test of Heller, Heller, and Gorfine (2013), note
that if \(X\) and \(Y\) (in \(\mathbb{R}^p\) and \(\mathbb{R}^q\), say)
are continuous and are dependent, then there exists a point
\((x_0,y_0)\) in the sample space of \((X,Y)\), and radii \(R_x\) and
\(R_y\) around \(x_0\) and \(y_0\), respectively, such that the joint
distribution of \((X,Y)\) is different from the product of the marginals
in the neighborhood defined by \(R_x\) and \(R_y\). The next step is the
introduction of a distance functions \(d\) in the sample spaces of \(X\)
and \(Y\), and following Heller, Heller, and Gorfine (2013) we do not
distinguish between these distance functions in our notation. Consider
the indicator functions \(1\{d(x_0,X) \leq R_x\}\) and
\(1\{d(y_0,Y) \leq R_y\}\). For a sample
\(\{(X_1,Y_1),\ldots,(X_n,Y_n)\}\) one gets \(n\) pairs of values of
zeros and ones from the indicator functions that can be set up in a
contingency table structure. Evidence against independence may then be
quantified by Pearson's chi-square test statistic or the likelihood
ratio test for \(2 \times 2\) contingency tables.

The data is used to guide in the choice of \((x_0,y_0)\), \(R_x\) and
\(R_y\). For every sample point \((x_i,y_i)\), that point is in its turn
chosen to be \((x_0,y_0)\) and for every sample point \((x_j,y_j)\),
\(j \neq i\), it is chosen in turn to define \(R_x = d(x_i,x_j)\) and
\(R_y = d(y_i,y_j)\) (thus defining the locality of the test). The
remaining \(n-2\) observations are then inserted in the indicator
functions. For every pair \((i,j)\) based on this one can construct a
classic test statistic \(S(i,j)\) for a Pearson chi-square test for a
\(2 \times 2\) contingency table. To test for independence these
quantities are aggregated in a test statistic
\(T=\sum_{i,j; j\neq i}S(i,j)\). Critical values are obtained by
resampling. The test, taking local properties into account, works very
well for the circle example and several other examples, both similar to
the circle example and not. See also Heller et al. (2016).

The next paper in this category, Reshef et al. (2011), is published in
Science. The idea behind their MIC (Maximal Information Coefficient)
statistic consists in computing the mutual information \(I\) as defined
in \eqref{eq:Kullback} \emph{locally} over a grid in the data set and then
take as statistic the maximum value of these local information measures
as obtained by maximizing over a suitable choice of grid. The authors
compare with several other classifiers on simulated and real data with
apparently good results. Some limitations of the method are identified
in a later article by Reshef et al. (2013). Another follow-up article is
Chen et al. (2016). See also Kinney and Atwal (2014), Reshef et al.
(2014) and Murrell, Murrell, and Murrell (2014). There are, however,
publications where the results are more mixed. See Simon and Tibshirani
(2014) and Gorfine, Heller, and Heller (2012). In particular these two
papers give several examples where the MIC is clearly inferior to dcov.

The two final papers in this category are Wang et al. (2015) and Wang et
al. (2017). In both papers the authors defines the locality by means of
a neighborhood of \(X\) and then consider suitable \(Y\)-values.
Consequently, as remarked by Wang et al. (2015) their test may be best
suited to nonlinear regression alternatives of the form
\(Y_i = f(X_i)+\varepsilon_i\). Wang et al. (2015) denote by
\((X_i,Y_i)\) \(i=1,\ldots,n\), \(n\) observations of the stochastic
variables \(X\) and \(Y\) in the construction of their CANOVA test
statistic. They define the \emph{within neighborhood sum of squares}
statistic as
\[W = \sum_{i,j} (Y_i-Y_j)^2, \;\; j < i,: \quad |\textrm{rank}(X_i)-\textrm{rank}(X_j)| < K\]
where \(K\) is an integer constant which is supposed to be chosen by the
user. Then \(|\textrm{rank}(X_i)-\textrm{rank}(X_j)|\) defines the
\(X\)-neighborhood structure. The assumption of CANOVA is that
dependence should imply that ``similar/neighbor \(X\)-values lead to
similar \(Y\)-values''. Thus when \(X\) and \(Y\) are dependent, small
values of \(W\) are expected. Critical values of \(W\) is determined by
permutations. The test for 4 values of \(K\) (2, 4, 8, 12) is compared
to a number of other tests among them Pearson's \(\rho\), the Kendall
and Spearman correlation coefficients, dcov, MIC, and the Hoeffding-test
based on the empirical distribution function. The CANOVA tests does not
do particularly well for linear models and it fails for the circle data
with weak noise, but its performance on the tested nonlinear examples is
very good.

The paper by Wang et al. (2017) follows much of the same pattern and
ideas. The \(X\)-values are first used to construct bagging
neighborhoods. And then they get an out-of-bag estimator of \(Y\), based
on the bagging neighborhood structure. The square error is calculated to
measure how well \(Y\) is predicted by \(X\). Critical values are again
obtained by permutations in the resulting statistic. In a comparison
with other methods five out of eight examples consist of various
sinusoidal function with added noise, where the new test does very well.

\section{\texorpdfstring{Beyond Pearson's \(\rho\): Local
dependence}{Beyond Pearson's \textbackslash{}rho: Local dependence}}\label{local}

The test functionals treated in Section \ref{global} deal with the the
second aspect of modelling dependence stated in the beginning of that
section, namely that of \emph{testing} of independence. These
functionals all do so by the computation of one non-negative number,
which is derived from local properties in Section \ref{global-local}.
This number properly scaled may possibly be said to deal with the the
first aspect stated, namely that of \emph{measuring} the strength of
dependence. But, as such, it can be faulted in several ways. Unlike the
Pearson \(\rho\), these functionals do not distinguish between positive
and negative dependence, and it is not local, thus not allowing for
stronger dependence in multivariate tails as is felt intuitively is the
case for data in finance for example.

In Section \ref{lgc} the main story will be the treatment of local
Gaussian correlation which in a sense returns to the Pearson \(\rho\)
but a local version of \(\rho\) which satisfies many of Rényi (1959)'s
requirements. But first, in the present section, we go back to some
other attempts to define local dependence, starting with a remarkable
paper by Lehmann (1966), who manages to define positive and negative
dependence in a quite general nonlinear situation.

\subsection{Quadrant dependence}\label{quadrant}

Lehmann's theory is based on the concept of quadrant dependence.
Consider two random variables \(X\) and \(Y\) with cumulative
distribution \(F_{X,Y}\). Then the pair \((X,Y)\) or its distribution
function \(F_{X,Y}\) is said to be positively quadrant dependent if

\begin{equation}
P(X \leq x, Y \leq y) \geq P(X \leq x)P(Y \leq y) \quad \textrm{for all} \;\; (x,y).
\label{eq:pos-quad1}
\end{equation}

Similarly, \((X,Y)\) or \(F_{X,Y}\) is said to be negatively quadrant
dependent if \eqref{eq:pos-quad1} holds with the central inequality sign
reversed.

The connection between quadrant dependence and Pearson's \(\rho\) is
secured through a lemma of Hoeffding (1940). The lemma is a general
result and resembles the result by Székely (2002) in his treatment of
the so-called Cramér functional, a forerunner of the Cramér - von Mises
functional. If \(F_{X,Y}\) denotes the joint and \(F_X\) and \(F_Y\) the
marginals, then assuming that the necessary moments exist,
\[E(XY)-E(X)E(Y) = \int_{-\infty}^{\infty}\int_{-\infty}^{\infty}\left(F_{XY}(x,y)-F_X(x)F_Y(y)\right)\di x \di y.\]
It follows immediately from definitions that if \((X,Y)\) is positively
quadrant dependent (negatively quadrant dependent), then for the
Pearson's \(\rho\), \(\rho \geq 0\) (\(\rho \leq 0\)). Similarly, it is
shown by Lehmann that if \(F_{X,Y}\) is positively quadrant dependent,
then Kendall's \(\tau\), Spearman's \(\rho_S\), and the quadrant measure
\(q\) defined by Blomqvist (1950) are all non-negative. The paper by
Blomqvist is an even earlier paper where positive and negative
dependence were considered in a nonlinear case, and using quadrants
centered at the median. An analogous result holds in the negatively
quadrant dependent case.

Lehmann (1966) introduces two additional and stronger concepts of
dependence. The first is regression dependence. Definition
\eqref{eq:pos-quad1} can be written
\[P(Y \leq y|X \leq x) \geq P(Y \leq y)\] but following Lehman, it may
be felt that the intuitive concept of positive dependence is better
represented by the stronger condition

\begin{equation}
P(Y \leq y|X \leq x) \geq P(Y \leq y|X \leq x') \quad \textrm{for all} \quad x < x' \;\; \textrm{and all} \;\; y.
\label{eq:pos-quad2}
\end{equation}

Rather than \eqref{eq:pos-quad2}, Lehmann considers the stronger condition
\[P(Y \leq y|X=x) \quad \textrm{is non-increasing in} \;\;x\] which was
discussed earlier by Tukey (1958). The concept of negative regression
dependence is defined by an obvious analog.

Finally, Lehmann introduces a stronger type of dependence still, by
requiring the conditional density of \(Y\) given \(x\) to have a
monotone likelihood ratio. Assuming the existence of a density
\(f = f_{X,Y}\), the condition may be written formally as
\[f(x,y')f(x',y) \leq f(x,y)f(x',y') \quad \textrm{for all} \quad x < x', \;\; y < y'.\]
If the inequality is reversed \((X,Y)\) is said to be negatively
likelihood ratio dependent. The bivariate normal distribution is
positively or negatively likelihood ratio dependent according to
\(\rho \geq 0\) or \(\rho \leq 0\). We will briefly return to these
dependence concepts when we get to the local Gaussian correlation in
Section \ref{lgc}.

\subsection{Local measures of dependence}\label{local-measures}

As mentioned already, econometricians have long looked for a formal
statistical way of describing the shifting region-like dependence
structure of financial markets. It is obvious that when the market is
going down there is a stronger dependence between financial objects, and
very strong in case of a panic. Similar effects, but perhaps not quite
so strong, appear when the market is going up. But how should it be
quantified and measured? This is important in finance, not the least in
portfolio theory, where it is well-known, see e.g. Taleb (2007), that
ordinary Gaussian description does not work, and if used, may lead to
catastrophic results. Mainly two approaches have been used among
econometricians. The first is non-local and consists simply in using
copula theory, but it may not always be so easy to implement in a time
series and portfolio context. The other approach is local and is to use
``conditional correlation'' as in Silvapulle and Granger (2001) and
Forbes and Rigobon (2002). One then computes an estimate as in
\eqref{eq:estimated} of Pearson's \(\rho\) but in various regions of the
sample space, e.g.~in the tail of two distributions. Let \(R\) be such a
region. A conditional correlation estimate is then given by (we let
\(n_R\) be the number of observed pairs \((X_i,Y_i) \in R\) and
\((\overline{X}_R = n_R^{-1}\sum_{(X_i,Y_i) \in R}X_i\), and similarly
for \(\overline{Y}_R\)),
\[\widehat{\rho}_R = \frac{\sum_{(X_i,Y_i)\in R}(X_i-\overline{X}_R)(Y_i-\overline{Y}_R)}{\sqrt{\sum_{(X_i,Y_i) \in R}(X_i-\overline{X}_r)^2}\sqrt{\sum_{(X_i,Y_i) \in R}(Y_i-\overline{Y}_R)^2}}.\]
However, this estimate suffers from a serious bias, which is obvious by
using the ergodic theorem or the law of large numbers, in the sense that
for a Gaussian distribution it does not converge to \(\rho\). This is
unfortunate because if the data happen to be Gaussian, one would like
estimated correlations to be close or identical to \(\rho\) in order to
approximate the classic Gaussian portfolio theory of Markowitz (1952).
This requirement is consistent with Rényi's property (vii).

The bias is examined in Boyer, Gibson, and Loretan (1999). Consider, as
an example, a bivariate Gaussian distribution with correlation
\(\rho = 0.5\). The conditional correlation when one of the variables is
large, for example larger than its 75\% quantile, \(X > q_{75}\), is
reduced to 0.27, and as the quantile increases \(\rho_R\) converges to
zero. In the finance literature one has tried to correct for this for
instance in contagion studies, Forbes and Rigobon (2002). What is wrong
here, one may think, is that one tries to use a product moment
estimator, which is a linear concept, on a quantity \(\rho\) that in the
nonlinear case is better thought of as a distributional parameter. We
will return to this in Section \ref{lgc}, where a distributional
approach yields an estimate without bias.

Statisticians have also tried various other ways of describing local
dependence. We will report on two such attempts: Bjerve and Doksum
(1993) had the idea of trying to extend the relationship between
correlation and regression coefficients in a linear regression model to
a nonlinear situation. Recall that in a linear model
\(Y= \alpha + \beta X +\varepsilon\),
\[\rho = \beta \frac{\sigma_X}{\sigma_Y} = \frac{\beta \sigma_X}{\sqrt{(\beta \sigma_X)^2+\sigma_{\varepsilon}^2}}\]
where \(\sigma_X^2\), \(\sigma_Y^2\) and \(\sigma_{\varepsilon}^2\) are
the variances of \(X\), \(Y\) and \(\varepsilon\). Based on this formula
Bjerve and Doksum suggested a local measure of dependence, the
correlation curve, based on localizing \(\rho\) by conditioning on \(X\)
(note that methods of Section \ref{global-local} for aggregating local
dependence also conditioned on \(X\)). Consider a generalization of the
linear model to \(Y = f(X)+ g(X)\varepsilon\) where \(f\) and \(g\) are
continuous functions, \(f\) is in addition continuously differentiable,
and \(\varepsilon\) has zero mean and is independent of \(X\). The
correlation curve is defined by
\[\rho(x) = \rho_{X,Y}(x) = \frac{\beta(x)\sigma_X}{\sqrt{(\beta(x)\sigma_X)^2+ \sigma_{\varepsilon}^2(x)}}\]
where \[\beta(x) = \mu'(x) \quad \mbox{with} \quad \mu(x) = \E(Y|X=x),\]
and \(\sigma_{\varepsilon}^2(x) = \Var(Y|X=x)\). It is trivial to check
that \(\rho(x)\) reduces to \(\rho\) in the linear case. The quantities
\(\beta(x)\) and \(\sigma_{\varepsilon}^2(x)\) can be estimated by
standard nonparametric methods. The correlation curve inherits many of
the properties of \(\rho\), but it succeeds in several of the cases
where \(\rho\) fails to detect dependence, such as the parabola
\eqref{eq:parabola} in Section \ref{nonlinearity}. However, unlike
\(\rho\), it is not symmetric in \((X,Y)\). In fact, it depends only on
\(x\), whereas one would want it to depend on \((x,y)\) in such a way
that \(\rho_{X,Y}(x,y) = \rho_{Y,X}(y,x)\). This is of course due to the
conditioning on and regression on \(X\). Conditioning and regression on
\(Y\) would in general produce a different result. This brings out the
difference between regression analysis and multivariate analysis, where
\(\rho\) is a concept of the latter, which accidentally enters into the
first. Bjerve and Doksum do propose a solution to this dilemma, but it
is an ad hoc one. Moreover, it is not so difficult to find examples
where the correlation curve is zero even though there is dependence.
Some further references are Wilcox (2005; 2007).

In Section \ref{global-local} we saw that Heller, Heller, and Gorfine
(2013) used local contingency type arguments to construct a global test
functional. Such reasoning goes further back in time. Holland and Wang
(1987) consider continuous stochastic variables \((X,Y)\) defined on
\(\mathbb{R}^2\). Let \(R_{x,y}\) denote the rectangle containing the
point \((x,y)\) having sides of length \(\Delta x\) and \(\Delta y\).
Then, approximately, as the sides become small,
\[P_{x,y} \doteq P(X,Y) \in R_{x,y}) \approx f(x,y)\Delta x \Delta y\]
where \(f\) is the joint density function. We now imagine that the
sample space of \((X,Y)\) are covered by such non-overlapping rectangles
(cells). For each cell we pick one point \((x,y)\) contained in that
cell. Based on all these pairs \((x,y)\), construct a contingency table
with indices \((i,j)\) with the elements \(P_{x,y} = P_{ij}\). Now
consider four neighboring cells \((i,k)\), \((i,l)\), \((j,k)\) and
\((j,l)\) with \(i < j\) and \(k < l\) and with a point \((x,y)\) in the
cell defined by \((i,k)\) and using the simplified notation \(\Delta i\)
for \(\Delta x_i\). The cross-product ratio is

\begin{align*}
\alpha((i,k),(j,l)) &= \frac{P_{ik}P_{jl}}{P_{il}p_{jk}} \\
& \approx \frac{f(i,k)\Delta i\Delta k \cdot f(j,l)\Delta j\Delta l}{f(i,l)\Delta
i\Delta l \cdot f(j,k)\Delta j \Delta k} = \frac{f(i,k)f(j,l)}{f(i,l)f(j,k)}.
\end{align*}

Let \(\theta((i,k),(j,l)) = \ln \alpha ((i,k),(j,l))\). Letting the
sides of all four cells tend to zero and then taking limits, one obtains
\[\gamma(x,y) = \lim_{\Delta x \to 0,\Delta y \to 0} \frac{\theta[(x,y),(x + \Delta x, y + \Delta y)]}{\Delta x \Delta y} = \frac{\partial^2}{\partial x \partial y} \ln f(x,y),\]
which is the local dependence function. Implicitly it is assumed here
that both the mixed second order partial derivatives exist and are
continuous. For an alternative derivation based on limiting arguments of
local covariance functions and for properties and extensions we refer to
Jones (1996), Jones (1998), Jones and Koch (2003), Sankaran and Gupta
(2004) and Inci, Li, and McCarthy (2011).

The local dependence function does not take values between -1 and 1, and
it does not reduce to \(\rho\) in the Gaussian bivariate case. Actually,
in that case
\[\gamma(x,y) = \frac{\rho}{1-\rho^2}\frac{1}{\sigma_X \sigma_Y}.\] Both
the correlation curve and the local dependence function works well for
the example \(Y=X^2 + \varepsilon\), where \(\rho = 0\), see
\eqref{eq:parabola}, producing dependence proportional to \(x\) and
producing a sign of the dependence in accordance with intuition.

\section{\texorpdfstring{Beyond Pearson's \(\rho\): Local Gaussian
correlation}{Beyond Pearson's \textbackslash{}rho: Local Gaussian correlation}}\label{lgc}

The Pearson \(\rho\) gives a complete characterization of dependence in
a bivariate Gaussian distribution but, as has been seen, not for a
general density \(f(x,y)\) for two random variables \(X\) and \(Y\). The
idea of the Local Gaussian Correlation (LGC), introduced in Tjøstheim
and Hufthammer (2013) is to approximate \(f\) locally in a neighborhood
of a point \((x,y)\) by a bivariate Gaussian distribution
\(\psi_{x,y}(u,v)\), where \((u,v)\) are running variables. In this
neighborhood we get close to a complete local characterization of
dependence, its precision depending on the size of the neighborhood and
of course on the properties of the density at the point \((x,y)\). In
practice it has to be reasonably smooth. This section of the paper gives
a suvey of some of the results obtained so far.

\subsection{Definition and examples}\label{definition-and-examples}

For notational convenience in this section we write \((x_1,x_2)\)
instead of \((x,y)\), and, by a slight inconsistency of notation,
\(x=(x_1,x_2)\). Similarly, \((u,v)\) is replaced by \(v=(v_1,v_2)\).
Then, in this notation, letting \(\mu(x) = (\mu_1(x),\mu_2(x))\) be the
mean vector of \(\psi\), \(\sigma(x) = (\sigma_1(x),\sigma_2(x))\) the
vector of standard deviations and \(\rho(x)\) the correlation of
\(\psi\), the approximating density \(\psi\) is then given by

\begin{align*}
&\psi(v,\mu_1(x),\mu_2(x),\sigma_1^2(x),\sigma_2^2(x),\rho(x)) = \frac{1}{2\pi\sigma_1(x)\sigma_2(x)\sqrt{1-\rho^2(x)}} \\
& \qquad \times \exp\bigg[-\frac{1}{2} \frac{1}{1-\rho^2(x)}\bigg(\frac{(v_1-\mu_1(x))^2}{\sigma_1^2(x)}-2\rho(x)\frac{(v_1-\mu_1(x))(v_2-\mu_2(x))}{\sigma_1(x)\sigma_2(x)} \\ & \qquad\qquad\qquad\qquad\qquad\qquad\qquad\qquad\qquad\qquad\qquad\qquad\qquad +\frac{(v_2-\mu_2(x))^2}{\sigma_2^2(x)}\bigg) \bigg].
\end{align*}

Moving to another point \(y=(y_1,y_2)\) in general gives another
approximating normal distribution \(\psi_y\) depending on a new set of
parameters\(\{\mu_1(y), \allowbreak \mu_2(y), \allowbreak\sigma_1(y), \allowbreak\sigma_2(y), \allowbreak\rho(y)\}\).
An exception is the case where \(f\) itself is Gaussian with parameters
\(\{\mu_1,\allowbreak\mu_2, \allowbreak\sigma_1, \allowbreak\sigma_2 ,\allowbreak\rho\}\),
in which case
\(\{\mu_1(x),\allowbreak\mu_2(x), \allowbreak\sigma_1(x), \allowbreak\sigma_2(x), \allowbreak\rho(x))\} \allowbreak\equiv \{\allowbreak\mu_1, \allowbreak\mu_2, \allowbreak\sigma_1, \allowbreak\sigma_2, \allowbreak\rho\}\).
This means that the bias of the conditional correlation described in
Section \ref{local} is avoided and it means that the property (vii) in
Rényi (1959)'s scheme is satisfied.

To make this into a construction that can be used in practice one must
define the vector population parameter
\(\theta(x) \doteq \{\mu_1(x),\mu_2(x),\sigma_1(x), \sigma_2(x), \rho(x)\}\)
and estimate it. Fortunately, this is a problem that has been treated in
larger generality by Hjort and Jones (1996) and Loader (1996). They
looked at the problem of approximating \(f(x)\) with a general
parametric family of densities, the Gaussian being one such family. Here
\(x\) in principle can have a dimension ranging from 1 to \(p\), but
with \(p=1\) mostly covered in those publications. They were concerned
with estimating \(f\) rather than the local parameters, one of which is
the LGC \(\rho(x)\). But their estimation method using local likelihood
is applicable also for estimating local parameters, and will be followed
here.

But first we need a more precise definition of \(\theta(x)\). This can
be done in two stages using a neighborhood defined by bandwidths
\(b=(b_1,b_2)\) in the \((x_1,x_2)\) direction, and then letting
\(b \to 0\) component-wise. (Alternatively one could use a set of
smoothness conditions requiring not only \(f\) and \(\psi\) to coincide
at \(x\), but also first and second order derivatives as indicated in
Berentsen et al. (2017), but the resulting equations are in general
difficult to solve).

A suitable function measyring the difference between \(f\) and \(\psi\)
is defined by

\begin{equation}
q = \int K_b(v-x)[\psi(v,\theta(x))-\ln \psi\{v,\theta(x)\}f(v)]\di v
\label{eq:Kullback-distance}
\end{equation}

where
\(K_b(v-x) = (b_1b_2)^{-1}K_1(b_1^{-1}(v_1-x_1))K_2(b_2^{-1}(v_2-x_2))\)
is a product kernel. As is seen in Hjort and Jones (1996), the
expression in \eqref{eq:Kullback-distance} can be interpreted as a locally
weighted Kullback-Leibler distance from \(f\) to
\(\psi(\cdot,\theta(x))\). We then obtain that the minimizer
\(\theta_b(x)\) (also depending on \(K\)) should satisfy

\begin{equation}
\int K_b(v-x)\frac{\partial}{\partial \theta_j} [\ln\{\psi(v,\theta(x))\}f(v)-\psi(v,\theta(x)] \di v = 0, \;\;j=1,\ldots,5.
\label{eq:score-equation}
\end{equation}

In the first stage we define the population value \(\theta_b(x)\) as the
minimizer of \eqref{eq:Kullback-distance}, assuming that there is a unique
solution to \eqref{eq:score-equation}. The definition of \(\theta_b(x)\)
and the assumption of uniqueness are essentially identical to those used
in Hjort and Jones (1996) for more general parametric families of
densities. A trivial example where \eqref{eq:score-equation} is satisfied
with a unique \(\theta_b(x)\) is when \(X \sim {\cal N}(\mu,\Sigma)\)
where \(\Sigma\) is the covariance matrix of \(X\). Another is the step
function of a Gaussian variable as given in equation (5) of Tjøstheim
and Hufthammer (2013).

In the next stage we let \(b \to 0\) and consider the the limiting value
\(\theta(x) = \lim_{b\rightarrow0}\theta_b(x)\). This is in fact
considered indirectly by Hjort and Jones (1996) on pp.~1627-1630 and
more directly in Tjøstheim and Hufthammer (2013), both using Taylor
expansion arguments. In the following we will assume that a limiting
value \(\theta(x)\) independent of \(b\) and \(K\) exists. (It is
possible to avoid the problem of a population value altogether if one
takes the view of some of the publications cited in Section
\ref{global-local} by just estimating a suitable dependence function).
Excepting the Gaussian or the Gaussian step model, it seems difficult to
find an explicit expression for \(\theta(x).\) We will return, however,
to a useful partial result for copulas later in this section.

In estimating \(\theta(x)\) and \(\theta_b(x)\) a neighborhood with a
finite bandwidth has to be used in analogy with nonparametric density
estimation. The estimate \(\widehat{\theta}(x) = \widehat{\theta}_b(x)\)
is obtained from maximizing a local likelihood. Given observations
\(X_1,\ldots,X_n\) the local log likelihood is determined by
\[L(X_1,\ldots,X_n,\theta(x)) = n^{-1}\sum_i K_b(X_i-x)\ln \psi(X_i,\theta(x))-\int K_b(v-x)\psi(v,\theta(x)) \di v.\]
The last (and perhaps somewhat unexpected) term is essential, as it
implies that \(\psi(x,\theta_b(x))\) is not allowed to stray far away
from \(f(x)\) as \(b \to 0\). It is also discussed at length in Hjort
and Jones (1996). (When \(b \to \infty\), the last term has 1 as its
limiting value and the likelihood reduces to the ordinary global
likelihood). Using the notation
\[u_j(\cdot,\theta) \doteq \frac{\partial}{\partial \theta_j} \ln \psi(\cdot,\theta),\]
by the law of large numbers, or by the ergodic theorem in the time
series case, assuming
\(\E\{K_b(X_i-x)\ln \psi(X_i,\theta_b(x))\} < \infty\), we have almost
surely

\begin{align}
\frac{\partial L}{\partial \theta_j} &= n^{-1}\sum_iK_b(X_i-x)u_j(X_i,\theta_b(x))-\int K_b(v-x)u_j(v,\theta_b(x))\psi(v,\theta_b(x)) \di v \nonumber \\
&\to \int K_b(v-x)u_j(v,\theta_b(x))[f(v)-\psi(v,\theta_b(x))] \di v.
\label{eq:LGC-emp}
\end{align}

Putting the expression in the first line of \eqref{eq:LGC-emp} equal to
zero yields the local maximum likelihood estimate
\(\widehat{\theta}_b(x)=\widehat{\theta}(x)\) of the population value
\(\theta_b(x)\) which satisfies \eqref{eq:score-equation}.

We see the importance of the additional last term in the local
likelihood by letting \(b \to 0\), Taylor expanding and requiring
\(\partial L/\partial \theta_j = 0\), which leads to
\[u_j(x,\theta_b(x))[f(x)-\psi(x,\theta_b(x))] + O(b^{T}b) = 0\] where
\(b^{T}\) is the transposed of \(b\). It is seen that ignoring solutions
that yield \(u_j(x,\theta_b(x)) = 0\) requires \(\psi(x,\theta_b(x))\)
to be close to \(f(x)\).

An asymptotic theory has been developed in Tjøstheim and Hufthammer
(2013) for \(\widehat{\theta}_b(x)\) for the case that \(b\) is fixed
and for \(\widehat{\theta}(x)\) in the case that \(b \to 0\). The first
case is much easier to treat than the second one. In fact for the first
case the theory of Hjort and Jones (1996) can be taken over almost
directly, although it is extended to the ergodic time series case in
Tjøstheim and Hufthammer (2013). In the case that \(b\rightarrow 0\),
this leads to a slow convergence rate of \((n(b_1b_2)^{3})^{-1/2}\),
which is the same convergence rate as for the the estimated dependence
function treated in Jones (1996).

The local correlation is clearly dependent on the marginal distributions
of \(X_1\) and \(X_2\) as is Pearson's \(\rho\). This marginal
dependence can be removed by scaling the observations to a standard
normal scale. As mentioned in Section \ref{copula} about the copula, the
dependence structure is disentangled from the marginals by Sklar's
theorem. For the purpose of measuring local dependence in terms of the
LGC, at least for a number of purposes it is advantageous to replace a
scaling with uniform variables \(U_i = F_i(X_i)\) by standard normal
variables

\begin{equation}
Z =(Z_1,Z_2) = (\Phi^{-1}(F_1(X_1)),\Phi^{-1}(F_2(X_2))),
\label{eq:LGC-normal}
\end{equation}

where \(\Phi\) is the cumulative distribution of the standard normal
distribution. The local Gaussian correlation on the \(Z\)-scale will be
denoted by \(\rho_Z(z)\). Of course the variable \(Z\) cannot be
computed via the transformation \eqref{eq:LGC-normal} without knowledge of
the margins \(F_1\) and \(F_2\), but these can be estimated by the
empirical distribution function. Extensive use has been made of
\(\rho_Z(z_1,z_2)\), or rather \(\rho_{\widehat{Z}}(z_1,z_2)\) with
\(\widehat{Z}_i = \Phi^{-1}(\widehat{F}_i)\). Under certain regularity
conditions, as in the copula case, the difference between \(Z\) and
\(\widehat{Z}\) can be ignored in limit theorems. Using the sample of
pairs of Gaussian pseudo observations
\(\{\Phi^{-1}(\widehat{F}_1(X_{1i}),\Phi^{-1}(\widehat{F}_2(X_{2i})\}\),
\(i=1,\ldots,n\) one can estimate \(\rho_Z(z_1,z_2)\) by local log
likelihood as described above. Under regularity conditions the
asymptotic theory will be the same as in Tjøstheim and Hufthammer
(2013). In Otneim and Tjøstheim (2017; 2018) a further simplification is
made by taking \(\mu_{Z_i}(z) = 0\) and \(\sigma_{Z_i}(z) = 1\), in
which case the asymptotic theory simplifies and one obtains the familiar
nonparametric rate of \(O((nb_1b_2)^{-1/2})\) for
\(\widehat{\rho}_{\widehat{Z}}(z)\).

The choice of Gaussian margins in the transformation \eqref{eq:LGC-normal}
is not made without a purpose. Since the copula of \((X_1,X_2)\) is
defined as the distribution function of
\((U_1,U_2) = (F_1(X_1),F_2(X_2))\) one could in principle consider the
local Gaussian correlation \(\rho_{U}(u_1,u_2)\) of the variable
\((U_1,U_2)\) (or the corresponding pseudo uniforms). However, fitting a
family of Gaussian density functions to finite support variables
requires special considerations of boundary effects which makes this
approach unpractical and illogical. The choice of Gaussian margins is
natural since we are dealing with local Gaussian approximations.

The LGC can be defined for a time series \(\{X_t\}\) as well simply by
taking \(X_1 = X_t\) and \(X_2 = X_{t+s}\) for a lag \(s\). The
asymptotic theory in Tjøstheim and Hufthammer (2013) is in fact carried
through for a stationary ergodic process \(\{X_t\}\), resulting in a
local autocorrelation function. Similarly, for a pair of time series
\(\{X_t,Y_t\}\) one can define the local cross correlation function by
taking \(X_1 = X_t\) and \(X_2 = Y_{t+s}\). The asymptotic estimation
theory can be found in Lacal and Tjøstheim (2017a). The asymptotic
distribution is fairly complicated and cannot be expected to work well
for a moderate sample size because of the presence of the bandwidth
parameter \(b\). Instead we have used the ordinary bootstrap and the
block bootstrap to obtain confidence intervals. The validity of the
bootstrapping procedures is demonstrated in Lacal and Tjøstheim (2017b;
2017a).

It should be noted that, in general, the local autocorrelation function
defined in this way lacks the important property of being positive
definite, since several Gaussians are involved in its definition (In the
Gaussian case it is positive definite). One has to bear this in mind in
applications to independence testing, density estimation and local
spectral estimation. See e.g.~Otneim and Tjøstheim (2017; 2018) and
Jordanger and Tjøstheim (2017b; 2017a).

In practical applications of the LGC an important point consists in
choosing the bandwidth parameter \(b\). Generally a cross-validation
procedure has been used for this; see e.g.Berentsen and Tjøstheim
(2014). Note that choosing the bandwidth by cross-validation generally
results in more stable results for \(\rho_Z(z_1,z_2)\). This is not
surprising in view of the standardized region. In a bootstrapping
situation the cross-validation procedure takes much time, so that for
the estimation of \(\rho_{Z}(z_1,z_2)\) a plug-in formula
\(b = 1.75n^{-1/6}\) has also been used, which seems to be working
fairly well, and where some empirical reasoning for its justification is
presented in Otneim (2016), who uses a simplified model for the
asymptotics of \(\widehat{\rho}_Z(z_1,z_2)\).

\begin{figure}[h!]
\subfloat[Gaussian data, $ \rho -0.5 $\label{fig:maps1}]{\includegraphics[width=0.49\linewidth,height=0.18\textheight]{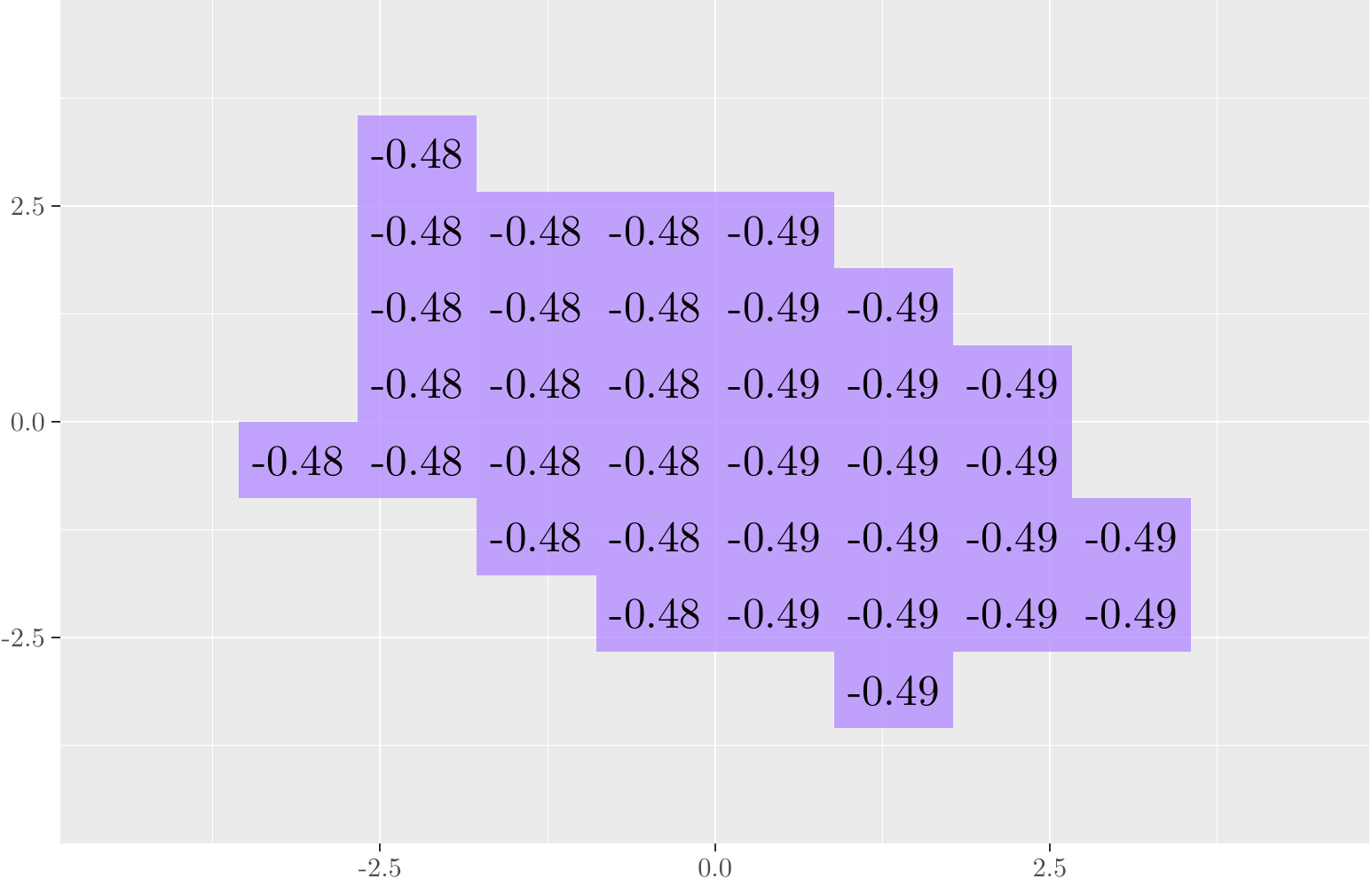} }\subfloat[Gaussian data, $ \rho = -0.5 $, normal scale\label{fig:maps2}]{\includegraphics[width=0.49\linewidth,height=0.18\textheight]{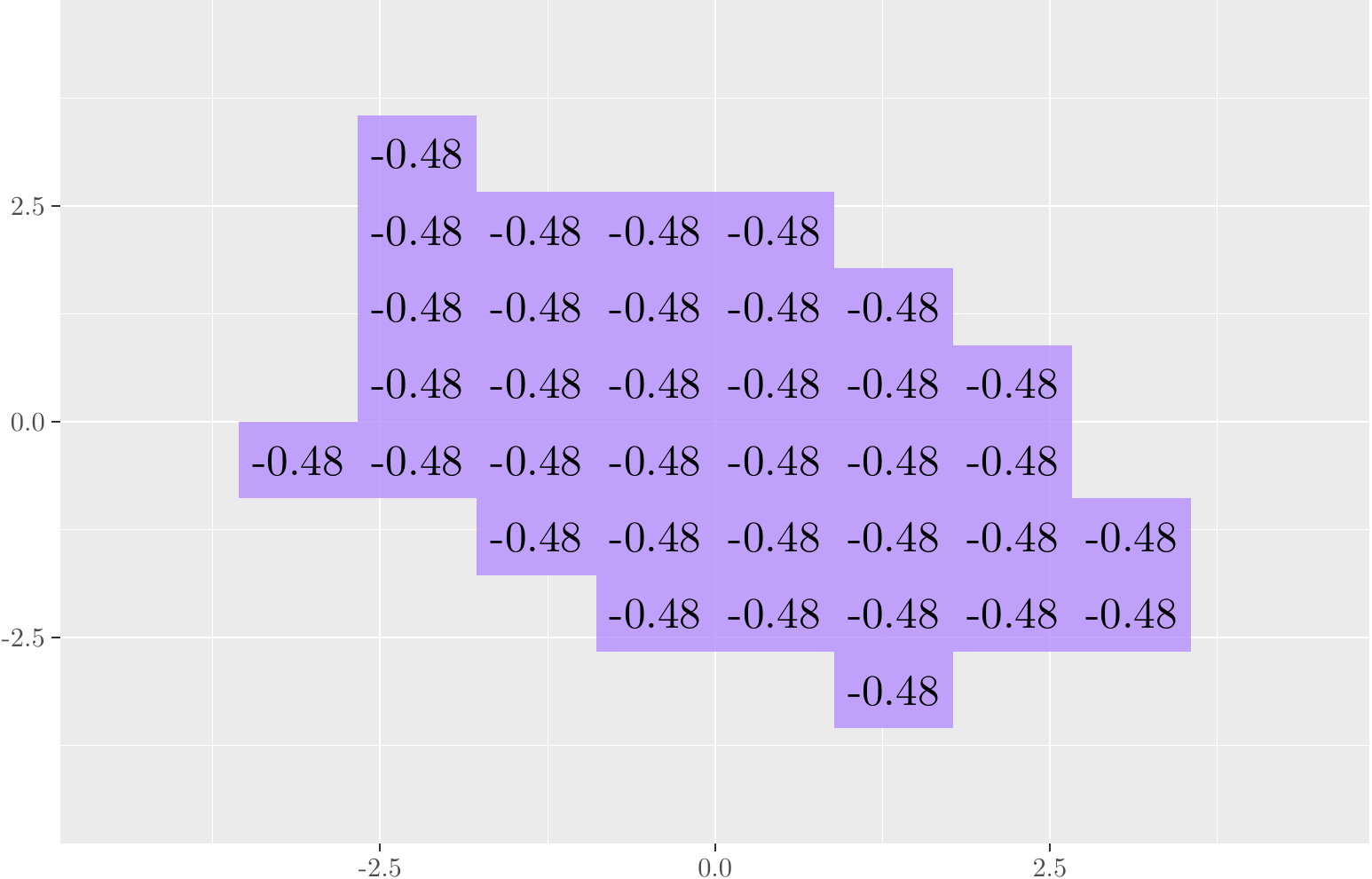} }\newline\subfloat[$t$-distributed data with 4 degrees \newline of  freedom\label{fig:maps3}]{\includegraphics[width=0.49\linewidth,height=0.18\textheight]{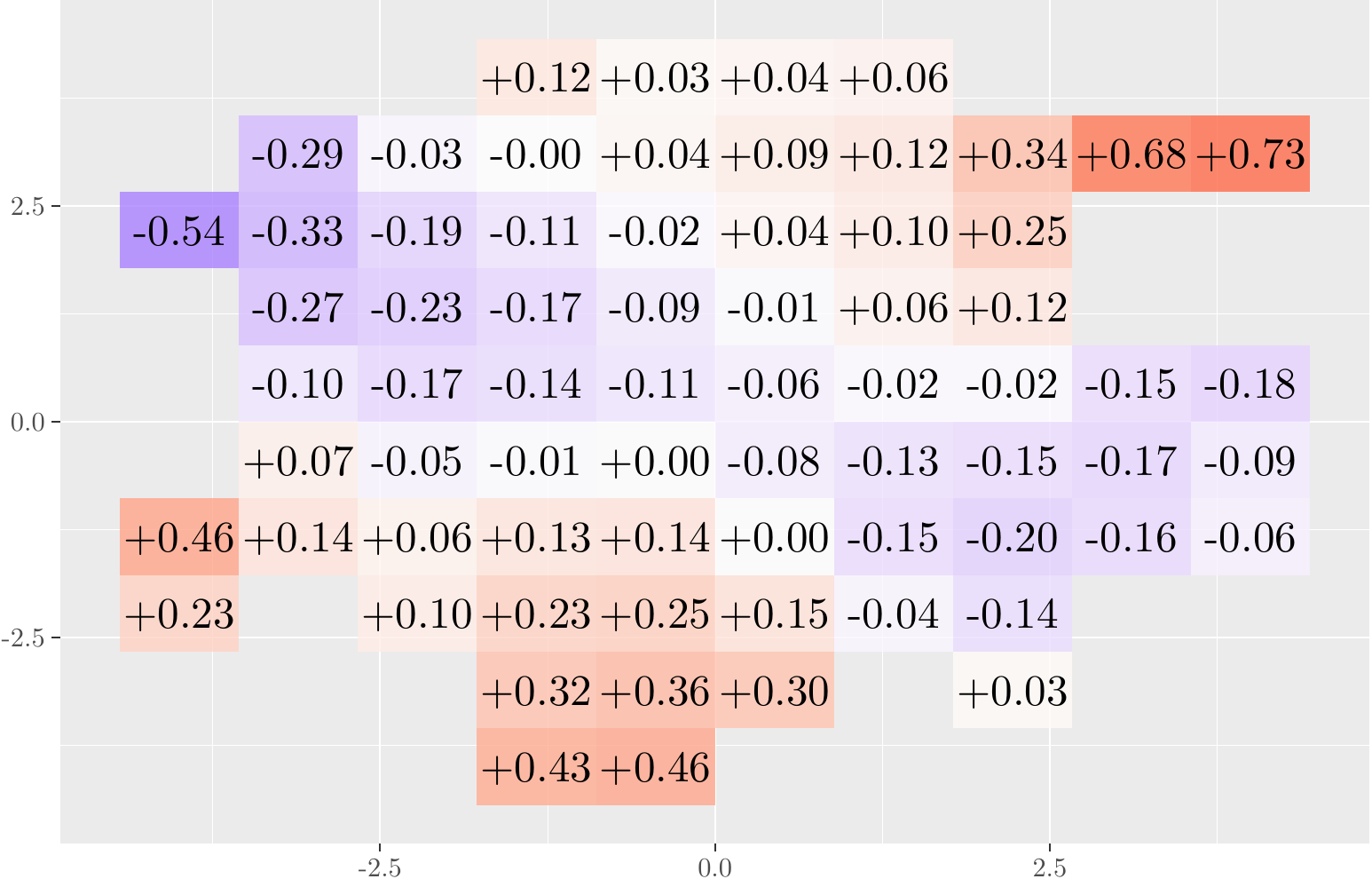} }\subfloat[$t$-distributed data with 4 degrees of \newline freedom, normal scale\label{fig:maps4}]{\includegraphics[width=0.49\linewidth,height=0.18\textheight]{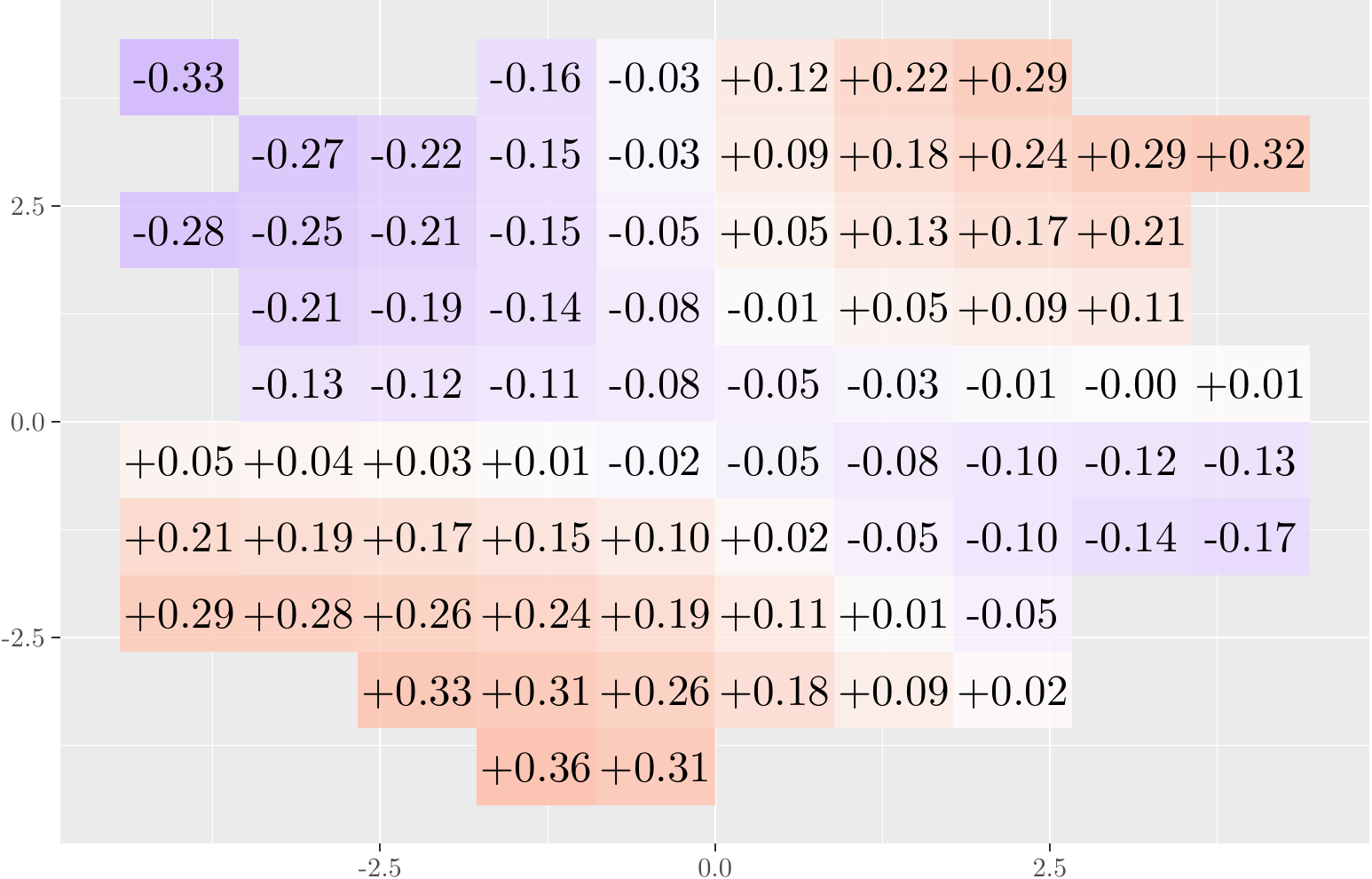} }\newline\subfloat[$X_t$ vs. $X_{t-1}$ where $X_t$ is generated from \newline a GARCH(1,1)-model\label{fig:maps5}]{\includegraphics[width=0.49\linewidth,height=0.18\textheight]{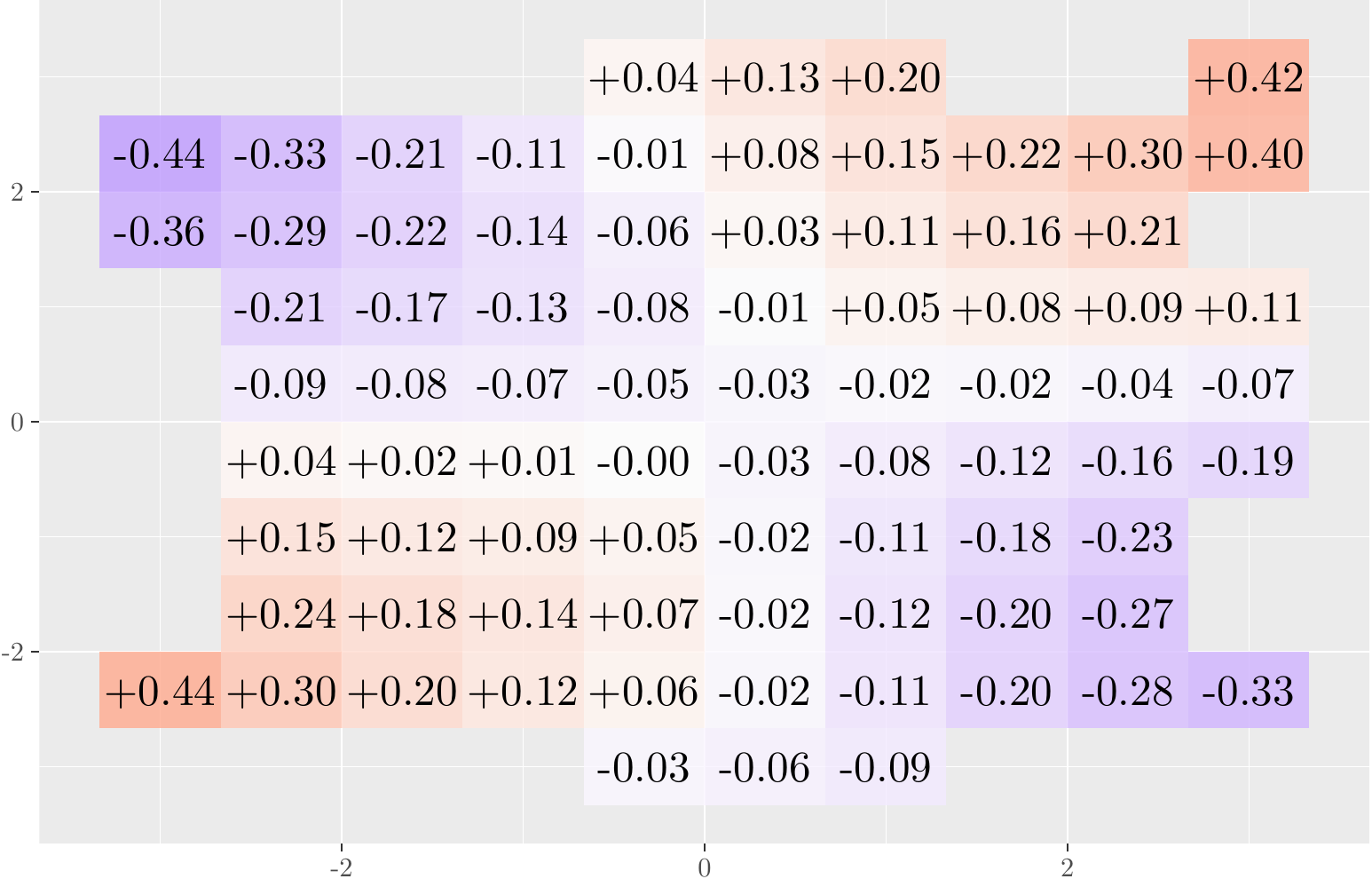} }\subfloat[$X_t$ vs. $X_{t-1}$ where $X_t$ is generated from \newline a GARCH(1,1)-model, normal scale\label{fig:maps6}]{\includegraphics[width=0.49\linewidth,height=0.18\textheight]{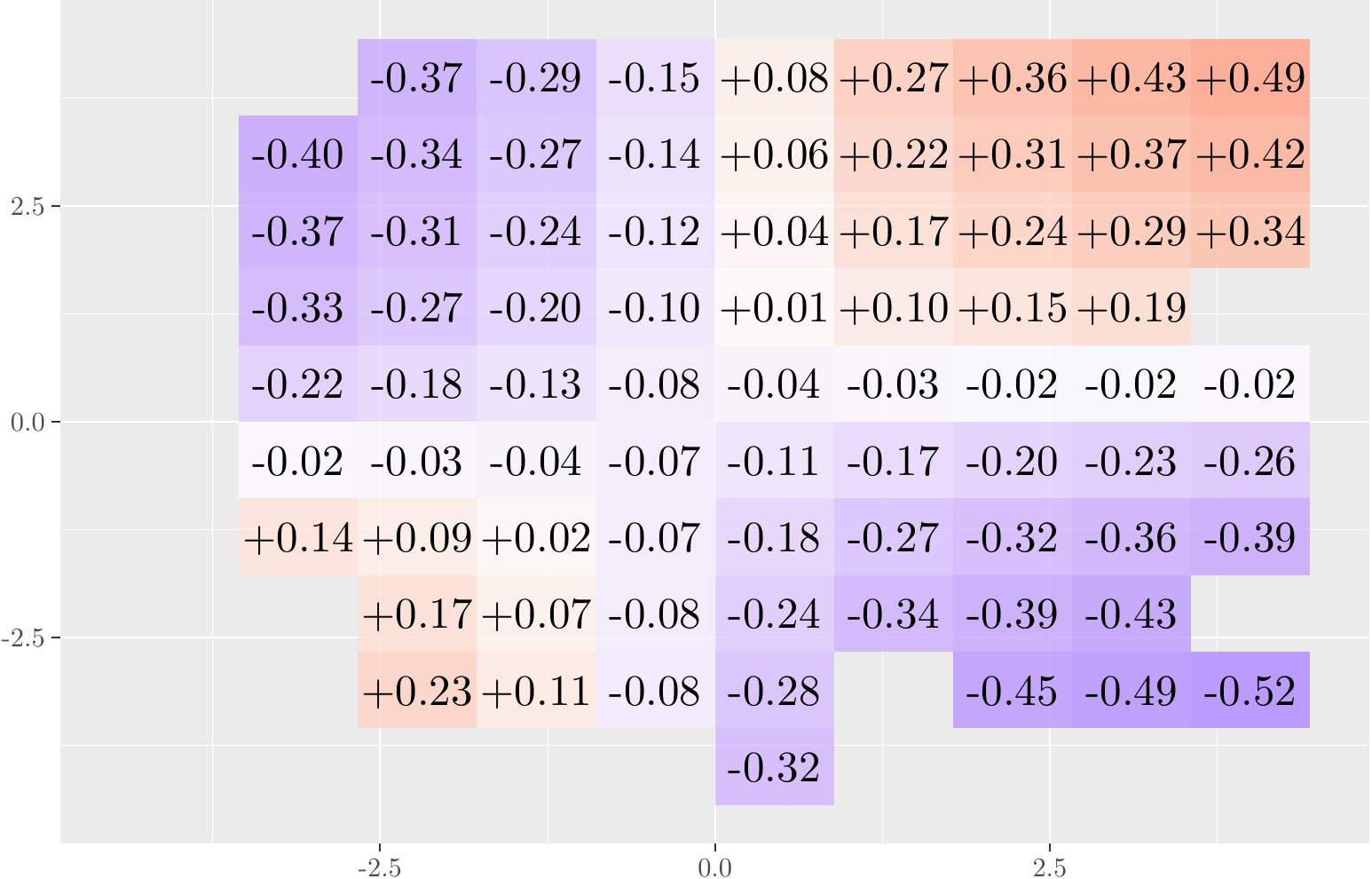} }\newline\subfloat[Log-returns data, S\&P 500 vs. FTSE 1000\label{fig:maps7}]{\includegraphics[width=0.49\linewidth,height=0.18\textheight]{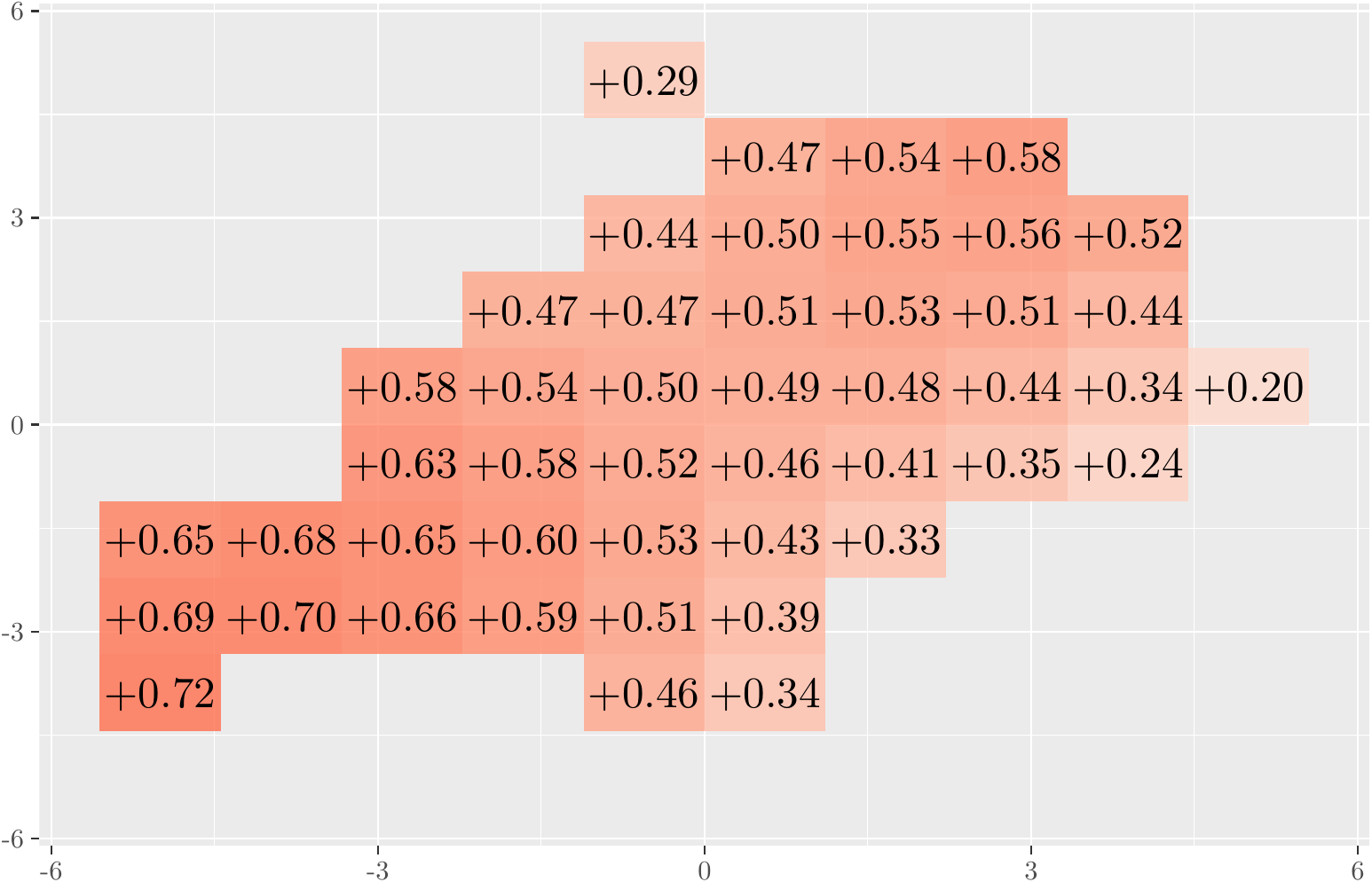} }\subfloat[Log-returns data, S\&P 500 vs. FTSE 1000, \newline normal scale\label{fig:maps8}]{\includegraphics[width=0.49\linewidth,height=0.18\textheight]{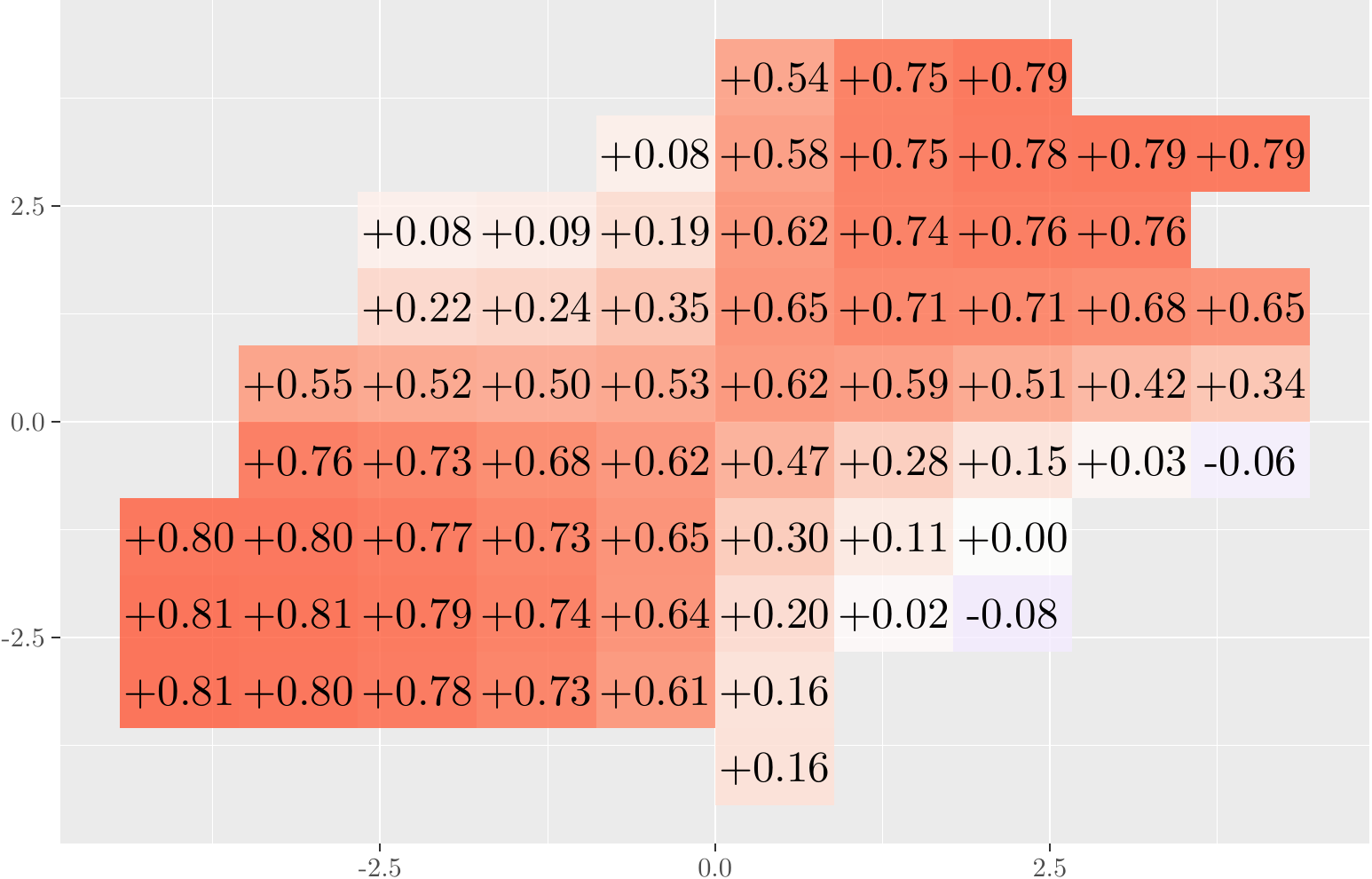} }\caption{Some examples of local correlation maps, $n = 784$ in all cases}\label{fig:maps}
\end{figure}

Figures \ref{fig:maps1} - \ref{fig:maps8} show some examples of
estimated LGC maps, using bandwidths determined by the cross-validation
procedure by Berentsen and Tjøstheim (2014), on both the original scale
(\(\widehat{\rho}(x)\)) and the Z-scale (\(\widehat{\rho}_Z(z))\) for
the following data:

\begin{enumerate}
\def\labelenumi{(\roman{enumi})}
\tightlist
\item
  Simulations from a bivariate standard normal distribution with
  correlation equal to -0.5.
\item
  Simulations from a bivariate \(t\)-distribution with 4 degrees of
  freedom and a global correlation of 0.
\item
  Simulations from the GARCH(1,1) model \eqref{eq:garch} with parameters
  \(\alpha = 0.1, \beta = 0.7\), \(\gamma = 0.2\) and
  \(\epsilon_t \sim \textrm{ iid } \mathcal{N}(0,1)\).
\item
  The daily log-return data as described in Section \ref{weaknesses}.
\end{enumerate}

We use the same sample size in the simulations as we have for the
log-return data: \(n = 784\).

For the first example in Figure \ref{fig:maps1} and \ref{fig:maps2}, the
original scale and the normal scale are almost identical (subject only
to the estimation error in the marginal distribution functions), and we
see the estimated LGC is close to the true value across all points (in
fact the sample Pearson correlation is \(\widehat\rho = -0.49\) for
these data). Radial and odd reflection symmetry of the LGC, to be
discussed in the next subsection, emerge clearly from the estimated
values in Figures \ref{fig:maps3} and \ref{fig:maps4}. Although \(X_t\)
and \(X_{t-1}\) are uncorrelated in the GARCH model, Figures
\ref{fig:maps5} and \ref{fig:maps6} show strong local dependence, with
much the same, although a bit stronger, pattern as for the
\(t\)-distribution, i.e.~positive dependence in the first and third
quadrant, and negative dependence in the second and fourth quadrant,
with increasing dependence away from the center.

For the return data in Figures \ref{fig:maps7} and \ref{fig:maps8} we
see clearly that the bivariate return distribution is not Gaussian,
since in particular there are large local correlations for both large
negative and large positive returns. The pattern remains the same on the
normal scale. Confidence intervals in such situations can be found in
Lacal and Tjøstheim (2017b; 2017a).

We will also illustrate the difference between \(\rho(x)\) and
\(\rho_Z(z)\) by looking at an exchangeable copula \(C\) with
\(C(u_1,u_2) = C(u_2,u_1)\). Several of the standard copulas such as the
Clayton, Frank and Gumbel copula are exchangeable. It is shown in
Berentsen et al. (2014) that it is possible to compute analytically
\(\mu(x)\), \(\Sigma(x)\) and in particular \(\rho(x)\) along the curve
defined by \(\{x = (x_1,x_2): F_1(x_1) = F_2(x_2)\}\). If
\(F_1 = F_2 = F\), then the curve reduces to the diagonal \(x=(d,d)\).
In such a case it is shown in Berentsen et al. (2014) that

\begin{equation}
\rho(d,d) = \frac{-C_{11}(F(d),F(d))\phi(\Phi^{-1}(F(d))}{\sqrt{\{\phi(\Phi^ {-1}(C_1(F(d),F(d))))\}^2+\{C_{11}(F(d),F(d))\phi(\Phi^{-1}(F(d)))\}^{2}}},
\label{eq:diagonal}
\end{equation}

where \(\phi\) is the standard normal density, where
\(C_1(u_1,u_2) = \partial /\partial u_1 C(u_1,u_2)\) and
\(C_{11} = \partial^2 /\partial u_1^2 C(u_1,u_2)\). Note that
\(C_1 = C_2\) and \(C_{11} = C_{22}\) due to exchangeability. It is seen
that in addition to the copula \(C\), the resulting local Gaussian
correlation \(\rho(x) = \rho(d,d)\) in \eqref{eq:diagonal} also depends on
the marginals \(F_1 = F_2 = F\).

For the \(Z\)-variable version in \eqref{eq:LGC-normal}, the cumulative
distribution function \(F_Z\) given by
\(F_Z(z_1,z_2) = C(\Phi(z_1),\Phi(z_2))\) is invariant to
transformations of the marginals as detailed above, and in particular by
inserting \(F_1 = F_2 = \Phi\) in \eqref{eq:diagonal}, it follows that
along the diagonal \(z_1=z_2=d\)

\begin{equation}
\rho_{Z}(d,d) = \frac{-C_{11}(\Phi(d),\Phi(d))\phi(d)}{\sqrt{\{\phi(\Phi^{-1}(C_1(\Phi(d),\Phi(d)))\}^2+\{C_{11}(\Phi(d),\Phi(d))\phi(d)\}^2}}.
\label{eq:canonical}
\end{equation}

For a parametric copula, the canonical local correlation \(\rho_Z(z)\)
will depend only on these parameters. In Berentsen et al. (2014) this is
used to pinpoint important dependence properties of exchangeable copulas
such as the Clayton, Gumbel and Frank copula and \(t\) copula. Among
other things it is shown for the Clayton and Frank copulas that
\(\rho_Z(d,d)\) tends to 1 as \(d \to -\infty\). Translated to the
language of returns of financial variables, this means that if variables
are described by the Clayton copula, then for very large negative
returns all guarding against risk in a portfolio of \(X_1\) and \(X_2\)
will disappear. Note that for the Gaussian copula, then
\(\rho(x_1,x_2) \equiv \rho_{Z}(z_1,z_2) \equiv \rho\).

\begin{figure}[t]
\includegraphics[width=1\linewidth]{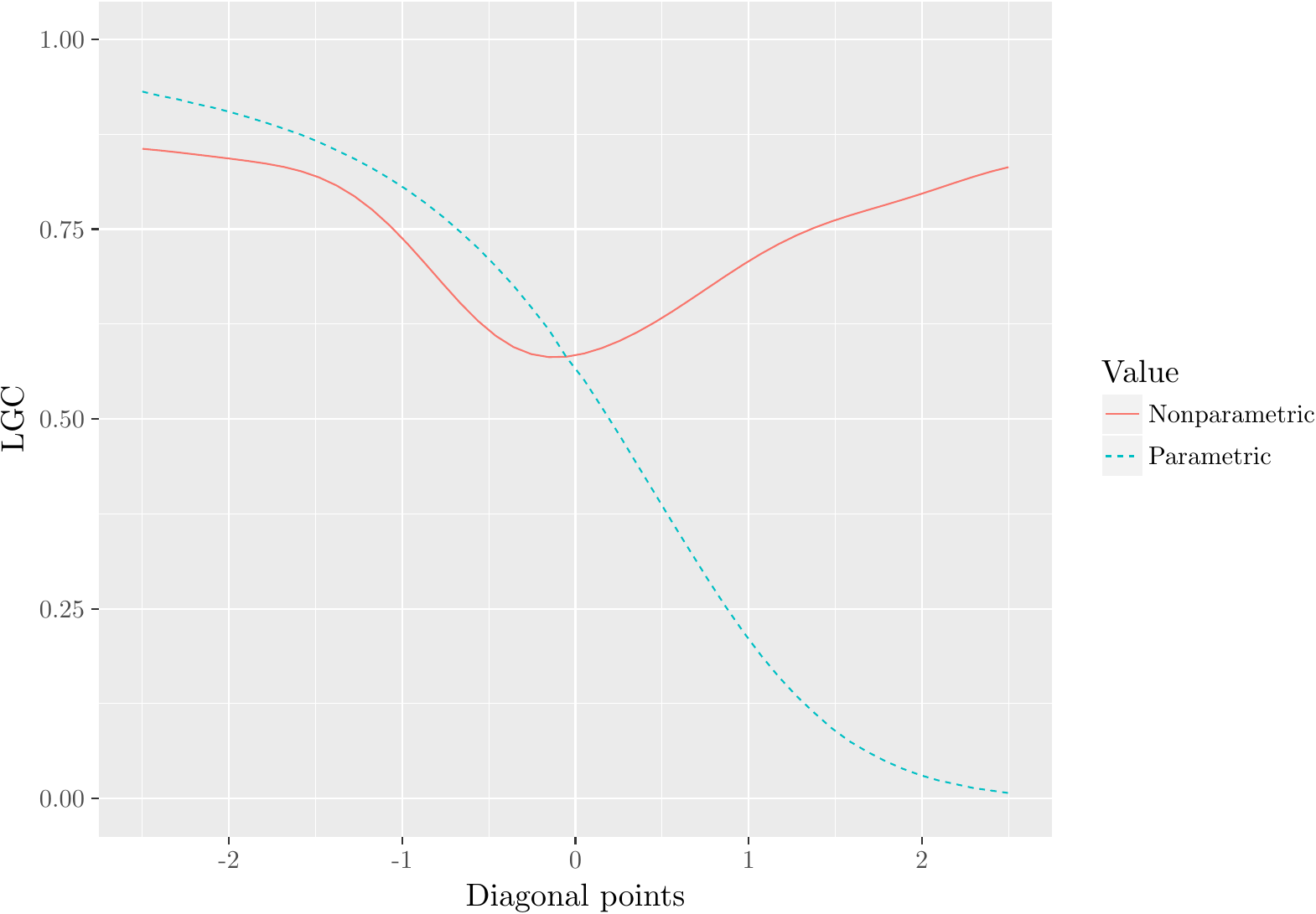} \caption{Parametric (Clayton) and nonparametric (LGC) correlation curves along the diagonal (on $Z$-scale) for the log-returns data set}\label{fig:diagonal}
\end{figure}

In Figure \ref{fig:diagonal} we have plotted the estimated LGC of the
returns (on \(Z\)-scale) on diagonal points, denoted by the label
«nonparametric». This curve corresponds to the values along the diagonal
through the first and third quadrant in the LGC map given in Figure
\ref{fig:maps8}. We have also estimated the canonical local correlation
\(\rho_Z(d,d)\) in~\eqref{eq:canonical}, by assuming a Clayton copula,
estimating the copula parameter. This curve is denoted by the label
«parametric» in Figure \ref{fig:diagonal}. As we see, this curve is
reasonably close to the estimated LGC in the left hand side of the plot
(i.e.~for negative returns), but as seen, the local Clayton-based
correlation goes towards zero on the right hand side of the plot
(i.e.~for positive returns), whereas the LGC is again increasing. Thus
the Clayton copula, which indeed is mostly intended to describe a
falling market, does not pick up the dependence patterns for positive
returns. The U-shaped LGC of daily returns often occurs. A possible
explanation for this, in addition to increasing dependence in a bull
market, is that during time periods with a bear market (falling prices),
even though the price trend is falling, for some days we actually will
see quite large positive daily returns, for example a bear market rally
(also known as ``sucker's rally''). These facts may explain why we
observe the high local correlation for positive returns too and not just
for negative returns. This effect is not typically seen when looking at
returns with longer time horizon, e.g.~one month. See Støve and
Tjøstheim (2014) for further examples.

Another example where an explicit formula for \(\rho(x)\) can be
obtained is in the parabola model \(X_2=X_1^2+\varepsilon\) that we have
mentioned several times before and where
\(\rho(x_1,x_1^2) = 2x_1/\sqrt{4x_1^2+\sigma_{\varepsilon}^2/\sigma_{X_1}^2}\),
which gives an intuitive picture of the local dependence properties of
this situation where the ordinary Pearson \(\rho\) fails so miserably.
It is seen that as \(\sigma_{\varepsilon} \to 0\),
\(\rho(x_1,x_2^2) \to 1\) if \(x_1 > 0\) and \(\rho(x_1,x_1^2) \to -1\)
as \(x_1 < 0\) and that \(\rho(0,0) = 0\), which all seem reasonable.

\subsection{Some additional properties of the local Gaussian
correlation}\label{lgc-properties}

Before stating particular properties of the LGC, it should be remarked
that the family of Gaussian distributions is especially attractive
because of its extensive and elegant theory in the multivariate case.
Most statisticians will agree that the Gaussian distribution is in a
class of itself when it comes to transparency and simplicity of
multivariate theory, and this is of course the main reason that it has
been so much used in applications, not only in finance but in a host of
multivariate problems. But sometimes, and this is certainly the case in
econometrics and finance, data do not follow a multivariate Gaussian,
and applications based on it can give very misleading results, see Taleb
(2007). The point of using the local Gaussian approximation is that one
can then move away from Gaussian distributions and describe much more
general situations, also multivariate thick tailed distributions like
those met in finance. At the same time one can exploit much of the
multivariate Gaussian theory locally. We have found this useful in a
number of papers extending in various directions: Støve, Tjøstheim, and
Hufthammer (2014), Støve and Tjøstheim (2014), Berentsen and Tjøstheim
(2014), Berentsen, Kleppe, and Tjøstheim (2014), Berentsen et al.
(2014), Berentsen et al. (2017), Lacal and Tjøstheim (2017b), Lacal and
Tjøstheim (2017a), Otneim and Tjøstheim (2017), Otneim and Tjøstheim
(2018). These results and further extensions will be collected in a
forthcoming book, Tjøstheim, Otneim, and Støve (2020).

We will now state some particulars for the local Gaussian correlation as
a measure of (local) dependence. We first look at the seven properties
(i) - (vii) stated in the beginning of Section \ref{global}. These were
first stated by Rényi (1959) as desirable properties for a dependence
measure to have.

\begin{enumerate}
\def\labelenumi{(\roman{enumi})}
\tightlist
\item
  \(\rho(x)\) is a local measure, but it is defined for any \(X_1\) and
  \(X_2\) under the regularity conditions discussed when defining the
  population value.
\item
  Because \(\rho(x)\) is a local measure we generally have
  \(\rho(x_1,x_2) \neq \rho(x_2,x_1)\) but
  \(\rho_{X_1,X_2}(x_1,x_2) = \rho_{X_2,X_1}(x_2,x_1)\) of course.
\item
  From definitions it follows that \(-1 \leq \rho(x) \leq 1\) and
  \(-1 \leq \widehat{\rho}(x) \leq 1\). We believe the possibility to
  measure negative as well as positive dependence to be one of the main
  assets of the local Gaussian correlation.
\item
  The condition \(\rho(x_1,x_2) \equiv 0\) implies independence of
  \(X_1\) and \(X_2\), but it is not sufficient. Then one must in
  addition require \(\mu_i(x_1,x_2) \equiv \mu_i(x_i)\) and
  \(\sigma_i(x_1,x_2) \equiv \sigma_i(x_i)\) for \(i=1,2\), Tjøstheim
  and Hufthammer (2013). In the global Gaussian case this is of course
  trivially fulfilled since all \(x\)-dependence disappears. Arguably,
  \(\rho_Z(z)\) is a better measure for the dependence between \(X_1\)
  and \(X_2\). (Note that even if \(Z_i \sim {\cal N}(0,1)\) we do not
  necessarily have \(\mu_i(z) = 0\), and \(\sigma_i(z)=1\) in a
  bivariate Gaussian approximation. For distributions where this is
  true, \(X_1\) and \(X_2\) are independent if and only if
  \(\rho(z) \equiv 0\)).
\item
  If \(X_1=f(X_2)\) or \(X_2=g(X_1)\), then, according to Tjøstheim and
  Hufthammer (2013), the limiting value \(\rho(x_1,x_2)\) as the
  neighborhood shrinks to the point \((x_1,x_2)\) is equal to 1 or -1,
  according to \(f\) or \(g\) having a positive or negative slope at
  \(x_2\) or \(x_1\). The same is true for a closed curve relationship
  \(g(X_1) + f(X_2) = c\) for a constant \(c\), then \(\rho(x)\) is
  equal to 1 or -1 along the curve, depending on whether the tangent at
  the point \(x\) has positive or negative slope.
\item
  Like the Pearson \(\rho\) the LGC \(\rho(x)\) depends on the
  marginals, but on the \(Z\)-scale, \(\rho_Z(z)\) is independent of
  marginals. Note that in the global case, transformations to standard
  normals of the margins, the Pearson's \(\rho\) reduces to the van der
  Waerden rank correlation, Waerden (1952).
\item
  In the Gaussian case we have \(\rho(x) \equiv \rho\) by construction.
  Note that this is also true for the Gaussian copula (under monotone
  transformations of the marginals from normal distributions).
\end{enumerate}

Concerning further properties of the LGC, first note the connection with
Lehmann's quadrant dependence. It is easily shown using the results in
that paper that if \(\rho(x) \geq 0\) or \(\rho(x) \leq 0\) for all
\(x\), then \((X_1,X_2)\) belongs to the class \({\cal F}\) defined in
Lehmann's paper. Moreover, if second moments exist, then
\(\rho(x) \geq 0\) for all \(x\), implies
\(\rho = \Corr(X_1,X_2) \geq 0\), and \(\rho(x) \leq 0\) for all \(x\)
implies \(\rho \leq 0\).

One main asset of Pearson's \(\rho\) is its scale invariance,
\(\rho_{\alpha_1+\beta_1X_1,\alpha_2+\beta_2 X_2} = \rho_{X_1,X_2}\).
There is a corresponding scale invariance for the LGC, but with the
proviso that the point \((x_1,x_2)\) is moved to the point
\((\alpha_1+\beta_1x_1,\alpha_2+\beta_2x_2)\). More generally, it is
shown by Tjøstheim and Hufthammer (2013) that for a vector \(\alpha\)
and a matrix \(A\) and for the stochastic variable \(Y = \alpha + AX\),
then for the local parameter vector \(\theta(y) = [\mu(y),\Sigma(y)]\)
at the point \(y=\alpha+Ax\), we have
\(\mu(y) = \alpha + \Sigma(x)\mu(x)\) and
\(\Sigma(y) = A \Sigma(x) A^{T}\). It follows that we have scale
invariance in the following sense
\(\rho_{Y_1,Y_2}(y_1,y_2) = \rho_{\alpha_1+\beta_1X_1,\alpha_2+\beta_2X_2}(\alpha_1+\beta_1x_1,\alpha_2+\beta_2x_2) = \rho_{X_1,X_2}(x_1,x_2)\).

In Tjøstheim and Hufthammer (2013) the transformation results for
\(\alpha+AX\) have been used to prove a number of symmetry properties.
In stating the results we have assumed \(\mu = \E(X)=0\), because
otherwise we may just center the density at \(\mu\), and make statements
about symmetry about \(\mu\). These symmetries are illustrated in
Figures \ref{fig:maps3}-\ref{fig:maps6}.

\begin{enumerate}
\def\labelenumi{(\roman{enumi})}
\tightlist
\item
  Radial symmetry: If \(f(x) = f(-x)\), then \(\Sigma(-x) = \Sigma(x)\),
  from which \(\rho(-x) = \rho(x)\), and \(\mu(-x) = -\mu(x)\).
\item
  Reflection symmetry: \(f(-x_1,x_2) = f(x_1,x_2)\) and/or
  \(f(x_1,-x_2) = f(x_1,x_2)\) imply \(\rho(-x_1,x_2) = -\rho(x_1,x_2\),
  \(\rho(x_1,-x_2)= -\rho(x_1,x_2)\),
  \(\mu_1(-x_1,x_2) = -\mu_1(x_1,x_2)\),
  \(\mu_2(-x_1,x_2) = \mu_2(x_1,x_2)\),
  \(\mu_1(x_1,-x_2) = \mu_1(x_1,x_2)\),
  \(\mu_2(x_1,-x_2) =-\mu_2(x_1,x_2)\).
\item
  Exchange symmetry: If \(f(x_1,x_2) = f(x_2,x_1)\), then
  \(\Sigma(x_1,x_2) = \Sigma(x_2,x_1)\) and hence
  \(\rho(x_1,x_2) = \rho(x_2,x_1)\).
\item
  Rotation symmetry: Then \(f(x) = \gamma(|x|)\) for a function
  \(\gamma\). If \(f\) is a spherical density, then \(f\) satisfies all
  the symmetry requirements mentioned above. It can be shown that in
  such a case \(\rho^2(x)\) takes its maximum along the lines
  \(x_1=x_2\) and \(x_1=-x_2\).
\end{enumerate}

In Tjøstheim and Hufthammer (2013) simulations are shown for
distributions satisfying these requirements.

\subsection{Testing for independence}\label{testing-lgc}

The LGC can be used for testing independence, and hence as a possible
supplement and competitor to the tests in Section \ref{global}. A
general estimated test functional for testing the independence between
two random variables \(X_1\) and \(X_2\) can be written

\begin{equation}
T_{n,b} = \int_S h(\widehat{\theta}_{n,b}(x)dF_n(x)
\label{eq:test-functional}
\end{equation}

where \(h\) is a measurable function, in general non-negative, and where
\(F_n\) is the empirical distribution function of \((X_1,X_2)\). This
functional estimates the functionals \[T_b = \int_S h(\theta_b(x)dF(x)\]
and \[T = \int_S h(\theta(x))dF(x).\] By choosing the set \(S\) one can
focus the test against special regions, e.g.~the tails. One can also
pre-test for symmetry properties discussed above and take advantage of
those as indicated in Berentsen and Tjøstheim (2014). Such a test can be
carried trough in increasing generality as demonstrated in Berentsen and
Tjøstheim (2014), and Lacal and Tjøstheim (2017b; 2017a).

Berentsen and Tjøstheim (2014) look at the case of iid pairs
\((X_1,X_2)\) and test for independence between \(X_1\) and \(X_2\).
Lacal and Tjøstheim (2017b; 2017a) consider the time series case where
for a time series \(\{X_t\}\) one can test for independence between
\(X_1=X_t\) and \(X_2 = X_{t+s}\). The corresponding LGC
\(\rho_{X_t,X_{t+s}}(x_1,x_2)\) is a function of \(s\) for a stationary
\(\{X_t\}\), and it defines a local autocorrelation function. Its
estimation theory is covered in Lacal and Tjøstheim (2017b). The theory
for the test functional is developed in Lacal and Tjøstheim (2017b;
2017a). From the reasoning in Section \ref{density}, the asymptotic
theory cannot be expected to be very accurate, and they have, as for
many (most) of the tests functionals of this type used bootstrapping in
tests for serial independence and the block bootstrap for the tests of
independence between two time series. In both cases the asymptotic
theory for the test functional and the validity of the bootstrap and the
block bootstrap have been established. The starting point for the theory
is a test functional based on a pairs \((X_t,X_{t+s})\) for the single
time series case and pairs \((X_t,Y_{t+s})\) for the two time series
case. These have been extended in Lacal and Tjøstheim (2017a) to
Box-Ljung type of functionals such as
\[T_{n,b}(\textrm{sum}) = \sum_{s_1}^{s_2}T_{n,b}^{(s)}\] where
\(T_{n,b}^{(s)}\) is of type \eqref{eq:test-functional} based on pairs
having a lag \(s\). Alternatively, they have also used a a test
functional
\[T_{n,b}(\mbox{max}) = \max_{s_1 \leq s \leq s_2} |T_{n,b}^{(s)}|.\]
These test functionals have been compared to several other test
functionals including Pearson's \(\rho\) and the dcov test. For the
simulation experiments tried in Lacal and Tjøstheim (2017b; 2017a) the
LGC test functional compares favorably with dcov. As is the case for the
dcov it beats the Pearson's \(\rho\) very clearly in non-Gaussian
situations, and for the sample sizes used it does not lose much compared
to \(\rho\) in a Gaussian situation, where \(\rho\) is optimal.

\section*{References}\label{references}
\addcontentsline{toc}{section}{References}

\hypertarget{refs}{}
\hypertarget{ref-aas:czad:frig:bakk:2009}{}
Aas, Kjersti, Claudia Czado, Arnoldo Frigessi, and Henrik Bakken. 2009.
``Pair-Copula Constructions of Multiple Dependence.'' \emph{Insurance:
Mathematics and Economics} 44 (2). Elsevier: 182--98.

\hypertarget{ref-abra:thom:1980}{}
Abrahams, Julia, and John B. Thomas. 1980. ``Properties of the Maximal
Correlation Function.'' \emph{Journal of the Franklin Institute} 310
(6): 317--23.

\hypertarget{ref-ande:1993}{}
Anderson, Theodore W. 1993. ``Goodness of Fit Tests for Spectral
Distributions.'' \emph{The Annals of Statistics} 21 (2): 830--47.

\hypertarget{ref-aron:1950}{}
Aronszajn, Nachman. 1950. ``Theory of Reproducing Kernels.''
\emph{Transactions of the American Mathematical Society} 68 (3):
337--404.

\hypertarget{ref-bart:1982}{}
Bartels, Robert. 1982. ``The Rank Version of von Neumann's Ratio Test
for Randomness.'' \emph{Journal of the American Statistical Association}
77 (377): 40--46.

\hypertarget{ref-bear:2010}{}
Beare, Brendan K. 2010. ``Copulas and Temporal Dependence.''
\emph{Econometrica} 78 (1): 395--410.

\hypertarget{ref-bera:bilo:lafa:2007}{}
Beran, Rudolf, Martin Bilodeau, and Pierre Lafaye de Micheaux. 2007.
``Nonparametric Tests of Independence Between Random Vectors.''
\emph{Journal of Multivariate Analysis} 98 (9): 1805--24.

\hypertarget{ref-bere:tjos:2014}{}
Berentsen, Geir Drage, and Dag Tjøstheim. 2014. ``Recognizing and
Visualizing Departures from Independence in Bivariate Data Using Local
Gaussian Correlation.'' \emph{Statistics and Computing} 24 (5):
785--801.

\hypertarget{ref-bere:cao:fran:tjos:2017}{}
Berentsen, Geir Drage, Ricardo Cao, Mario Francisco-Fernández, and Dag
Tjøstheim. 2017. ``Some Properties of Local Gaussian Autocorrelation and
Other Nonlinear Dependence Measures.'' \emph{Journal of Time Series
Analysis} 38 (2): 352--80.

\hypertarget{ref-bere:klep:tjos:2014}{}
Berentsen, Geir Drage, Tore Kleppe, and Dag Tjøstheim. 2014.
``Introducing \texttt{localgauss}, an R-Package for Estimating and
Visualizing Local Gaussian Corelation.'' \emph{Journal of Statistical
Software} 56 (12): 1--18.

\hypertarget{ref-bere:stov:tjos:nord:2014}{}
Berentsen, Geir Drage, Bård Støve, Dag Tjøstheim, and Tommy Nordbø.
2014. ``Recognizing and Visualizing Copulas: An Approach Using Local
Gaussian Approximation.'' \emph{Insurance: Mathematics and Economics}
57: 90--103.

\hypertarget{ref-berl:thom:2004}{}
Berlinet, Alain, and Christine Thomas-Agnan. 2004. \emph{Reproducing
Kernel Hilbert Spaces in Probability and Statistics}. Springer.

\hypertarget{ref-berr:samw:2017}{}
Berrett, Thomas B, and Richard J Samworth. 2017. ``Nonparametric
Independence Testing via Mutual Information.'' \emph{arXiv Preprint
arXiv:1711.06642}.

\hypertarget{ref-bick:rose:1973}{}
Bickel, Peter J., and Murray Rosenblatt. 1973. ``On Some Global Measures
of the Deviations of Density Function Estimators.'' \emph{The Annals of
Statistics} 1: 1071--95.

\hypertarget{ref-bill:2008}{}
Billingsley, Patrick. 2008. \emph{Probability and Measure}. John Wiley
\& Sons.

\hypertarget{ref-bjer:doks:1993}{}
Bjerve, Steinar, and Kjell Doksum. 1993. ``Correlation Curves: Measures
of Association as Function of Covariate Values.'' \emph{The Annals of
Statistics} 21 (2): 890--902.

\hypertarget{ref-blom:1950}{}
Blomqvist, Nils. 1950. ``On a Measure of Dependence Between Two Random
Variables.'' \emph{The Annals of Mathematical Statistics} 21 (4):
593--60.

\hypertarget{ref-blum:kief:rose:1961}{}
Blum, Julius R., Jack Kiefer, and Murray Rosenblatt. 1961.
``Distribution Free Tests of Independence Based on the Sample
Distribution Function.'' \emph{The Annals of Mathematical Statistics} 32
(2): 485--98.

\hypertarget{ref-boye:gibs:lore:1999}{}
Boyer, Brian H., Michael S. Gibson, and Mico Loretan. 1999. ``Pitfalls
in Tests for Changes in Correlation.'' Discussion Paper 597. Federal
Reserve Government Papers.

\hypertarget{ref-brad:1986}{}
Bradley, Richard C. 1986. ``Basic Properties of Strong Mixing
Conditions.'' Edited by E. Eberlein and M.S. Taqqu. Birkhäuser, Boston,
165--92.

\hypertarget{ref-brei:frie:1985}{}
Breiman, Leo, and Jerome H. Friedman. 1985. ``Estimating Optimal
Transformations for Multiple Regression and Correlation (with
Discussion).'' \emph{Journal of the American Statistical Association} 80
(391): 580--619.

\hypertarget{ref-broc:dech:sche:leba:1996}{}
Broock, William A., José A. Scheinkman, W. Davis Dechert, and Blake
LeBaron. 1996. ``A Test for Independence Based on the Correlation
Dimension.'' \emph{Econometric Reviews} 15 (3): 197--236.

\hypertarget{ref-carl:1988}{}
Carlstein, Edward. 1988. ``Degenerate \(U\)-Statistics Based on
Non-Independent Observations.'' \emph{Calcutta Statistical Association
Bulletin} 37 (1-2): 55--65.

\hypertarget{ref-chen:fan:2006}{}
Chen, Xiaohong, and Yanqin Fan. 2006. ``Estimation of Copula Based
Semiparametric Time Series Models.'' \emph{Journal of Econometrics} 130
(2): 307--35.

\hypertarget{ref-chen:zeng:luo:zhem:2016}{}
Chen, Yuan, Ying Zeng, Feng Luo, and Yuan Zheming. 2016. ``A New
Algorithm to Optimize Maximal Information Coefficient.'' \emph{PLoS ONE}
11 (6): e0157567.

\hypertarget{ref-chun:kann:ng:soho:1989}{}
Chung, J. K., P. L. Kannapan, C. T. Ng, and P. K. Sahoo. 1989.
``Measures of Distance Between Probability Distributions.''
\emph{Journal of Mathematical Analysis and Applications} 138 (1):
280--92.

\hypertarget{ref-chwi:gret:2014}{}
Chwialkowski, Kacper, and Arthur Gretton. 2014. ``A Kernel Independence
Test for Random Processes.'' In \emph{Proceedings of the 31st
International Conference on Machine Learning}, 32:1422--30.

\hypertarget{ref-csor:1985}{}
Csörgö, Sándor. 1985. ``Testing for Independence by the Empirical
Characteristic Function.'' \emph{Journal of Multivariate Analysis} 16:
290--99.

\hypertarget{ref-czak:fisc:1963}{}
Czáki, Péter, and János Fischer. 1963. ``On the General Notion of
Maximal Correlation.'' \emph{Magyar Tudományos Akad. Mat. Kutató
Intézetenk Közlemenényei (Publ. Math. Inst. Hungar: Acad. Sci.} 8:
27--51.

\hypertarget{ref-dars:nguy:1992}{}
Darsow, William F., Bao Nguyen, and Elwood T. Olsen. 1992. ``Copulas and
Markov Processes.'' \emph{Illinois Journal of Mathematics} 36 (4):
600--642.

\hypertarget{ref-datastream}{}
Datastream. 2018. Subscription service (Accessed June 2018).

\hypertarget{ref-davi:mats:miko:wan:2018}{}
Davis, Richard, Muneya Matsui, Thomas Mikosch, and Phyllis Wan. 2018.
``Applications of Distance Correlation to Time Series.''
\emph{Bernoulli} 24 (4A): 3087--3116.

\hypertarget{ref-dehe:1981b}{}
Deheuvels, Paul. 1981a. ``A Kolmogorov-Smirnov Type Test for
Independence and Multivariate Samples.'' \emph{Revue Roumaine de
Mathemátiques Pures et Appliquées} 26 (2): 213--26.

\hypertarget{ref-dehe:1981a}{}
---------. 1981b. ``An Asymptotic Decomposition for Multivariate
Distribution-Free Tests of Independence.'' \emph{Journal of Multivariate
Analysis} 11 (1): 102--13.

\hypertarget{ref-delg:1996}{}
Delgado, Miguel A. 1996. ``Testing Serial Independence Using the Sample
Distribution Function.'' \emph{Journal of Time Series Analysis} 17 (3):
271--86.

\hypertarget{ref-denk:kell:1983}{}
Denker, Manfred, and Gerhard L. Keller. 1983. ``On \(U\)-Statistics and
von Mises' Statistics for Weakly Dependent Processes.''
\emph{Zeitschrift Für Wahrscheinlichkeitstheorie Und Verwandte Gebiete}
64 (4): 505--22.

\hypertarget{ref-duec:edel:rich:2015}{}
Dueck, Johannes, Dominic Edelman, and Donald Richards. 2015. ``A
Generalization of an Integral Arising in the Theory of Distance
Correlation.'' \emph{Statistics and Probability Letters} 97: 116--19.

\hypertarget{ref-duec:edel:gnei:rich:2014}{}
Dueck, Johannes, Dominic Edelman, Tilmann Gneiting, and Donald Richards.
2014. ``The Affinely Invariant Distance Correlation.'' \emph{Bernoulli}
20 (4): 2305--30.

\hypertarget{ref-engl:1982}{}
Engle, Robert F. 1982. ``Autoregressive Conditional Heteroscedasticity
with Estimates of Varaiance of U.K. Inflation.'' \emph{Econometrica} 50
(4): 987--1008.

\hypertarget{ref-esca:vela:2006}{}
Escanciano, J. Carlos, and Carlos Velasco. 2006. ``Generalized Spectral
Tests for the Martingale Difference Hypothesis.'' \emph{Journal of
Econometrics} 134 (1): 151--85.

\hypertarget{ref-fan:demi:pene:sapo:2017}{}
Fan, Yanan, Pierre Lafaye de Micheaux, Spiridon Penev, and Donna
Sapolek. 2017. ``Multivariate Nonparametric Tests of Independence.''
\emph{Journal of Multivariate Analysis} 153: 189--210.

\hypertarget{ref-ferg:gene:hall:2000}{}
Ferguson, Thomas S., Christian Genest, and Marc Hallin. 2000.
``Kendall's Tau for Serial Dependence.'' \emph{Canadian Journal of
Statistics} 28 (3): 587--604.

\hypertarget{ref-ferr:vieu:2006}{}
Ferraty, Frédéric, and Philippe Vieu. 2006. \emph{Nonparametric
Functional Data Analysis: Theory and Practice}. Springer, New York.

\hypertarget{ref-fish:1915}{}
Fisher, Ronald A. 1915. ``Frequency Distribution of the Values of the
Correlation Coefficient in Samples of an Indefinitely Large
Population.'' \emph{Biometrika} 10 (4): 507--21.

\hypertarget{ref-fish:1921}{}
---------. 1921. ``On the Probable Error of a Coefficient of Correlation
Deduced from a Small Sample.'' \emph{Metron} 1: 3--32.

\hypertarget{ref-foki:pits:2017}{}
Fokianos, Konstantinos, and Maria Pitsillou. 2017. ``Consistent Testing
for Pairwise Dependence in Time Series.'' \emph{Technometrics} 59 (2):
262--70.

\hypertarget{ref-forb:rigo:2002}{}
Forbes, Kristin J., and Roberto Rigobon. 2002. ``No Contagion, Only
Interdependence: Measuring Stock Market Comovements.'' \emph{The Journal
of Finance} 57 (5): 2223--61.

\hypertarget{ref-fran:zako:2011}{}
Francq, Christian, and Jean-Michel Zakoian. 2011. \emph{GARCH Models:
Structure, Statistical Inference and Financial Applications}. John Wiley
\& Sons.

\hypertarget{ref-galt:1888}{}
Galton, Francis. 1888. ``Co-Relations and Their Measurement, Chiefly
from Anthropometric Data.'' \emph{Proceedings Royal Society, London} 45
(273--279): 135--45.

\hypertarget{ref-galt:1890}{}
---------. 1890. ``Kinship and Correlation.'' \emph{North American
Review} 150 (401): 419--31.

\hypertarget{ref-gebe:1941}{}
Gebelein, Hans. 1941. ``Das Statistische Problem Der Korrelation as
Variations - Und Egenwerthproblem Und Sein Zusammengeng Mit Der
Ausgleichsrechnung.'' \emph{ZAMM-Journal of Applied Mathematics and
Mechanics/Zeitschrift Für Angewandte Mathematik Und Mechanik} 21 (6):
364--79.

\hypertarget{ref-gene:nesl:2007}{}
Genest, Christian, and Johanna Nešlehová. 2007. ``A Primer on Copulas
for Count Data.'' \emph{ASTIN Bulletin} 37 (2): 475--515.

\hypertarget{ref-gene:remi:2004}{}
Genest, Christian, and Bruno Rémillard. 2004. ``Tests of Independence
and Randomness Based on the Empirical Copula Process.'' \emph{TEST} 13
(2): 335--69.

\hypertarget{ref-gene:ghou:remi:2007}{}
Genest, Christian, Kilani Ghoudi, and Bruno Rémillard. 2007.
``Rank-Based Extensions of the Brock, Dechert, Scheinkman Test.''
\emph{Journal of the American Statistical Association} 102 (480):
1363--76.

\hypertarget{ref-ghou:remi:2018}{}
Ghoudi, Kilani, and Bruno Rémillard. 2018. ``Serial Independence Tests
for Innovations of Conditional Mean and Variance Models.'' \emph{TEST}
27 (1): 3--26.

\hypertarget{ref-ghou:kulp:remi:2001}{}
Ghoudi, Kilani, Reg J. Kulperger, and Briuno Rémillard. 2001. ``A
Nonparametric Test of Serial Independence for Time Series and
Residuals.'' \emph{Journal of Multivariate Analysis} 79 (2): 191--201.

\hypertarget{ref-gorf:hell:hell:2018}{}
Gorfine, Malka, Ruth Heller, and Yair Heller. 2012. ``Comment on
``Detecting Novel Associations in Large Data Sets' by Reshef et al.",
Science Dec 16, 2011.''

\hypertarget{ref-gome:gome:mari:2003}{}
Gómez, Eusebio, Miguel A. Gómez-Villegas, and J.Miguel Mari'in. 2003.
``A Survey on Continuous Elliptical Vector Distributions.''
\emph{Revista Matemática Computense} 16: 345--61.

\hypertarget{ref-gran:lin:1994}{}
Granger, Clive, and Jin-Lung Lin. 1994. ``Using the Mutual Information
Coefficient to Identify Lags in Nonlinear Time Series Models.''
\emph{Journal of Time Series Analysis} 15 (4): 371--84.

\hypertarget{ref-gran:maas:raci:2004}{}
Granger, Clive, Esfandiar Maasoumi, and Jeffrey Racine. 2004. ``A
Dependence Metric for Possible Nonlinear Processes.'' \emph{Journal of
Time Series Analysis} 25 (5): 649--70.

\hypertarget{ref-gros:proc:1983}{}
Grassberger, Peter, and Itamar Procaccia. 1983. ``Measuring the
Strangeness of Attractors.'' \emph{Physica D: Nonlinear Phenomena} 9
(1-2): 189--208.

\hypertarget{ref-gret:2017}{}
Gretton, Arthur. 2017. ``Introduction to RKHS, and Some Simple Kernel
Algorithms.'' Lecture notes Gatsby Computational Neuroscience Unit.

\hypertarget{ref-gret:gyor:2010}{}
Gretton, Arthur, and László Györfi. 2010. ``Consistent Nonparametric
Tests of Independence.'' \emph{Journal of Machine Learning Research} 11
(Apr): 1391--1423.

\hypertarget{ref-gret:gyor:2012}{}
---------. 2012. ``Strongly Consistent Nonparametric Test of Conditional
Independence.'' \emph{Journal of Multivariate Analysis} 82 (6):
1145--50.

\hypertarget{ref-gret:bous:smol:scho:2005}{}
Gretton, Arthur, Olivier Bousquet, Alex Smola, and Bernhard Schölkopf.
2005. ``Measuring Statistical Dependence with Hilbert-Schmidt Norms.''
In \emph{International Conference on Algorithmic Learning Theory},
edited by S. Jain, U. Simon, and E. Tomita, 63--77. Springer, Berlin.

\hypertarget{ref-gret:herb:smol:bous:scho:2005}{}
Gretton, Arthur, Ralf Herbrich, Alexander Smola, Olivier Bousquet, and
Bernhard Schölkopf. 2005. ``Kernel Methods for Measuring Independence.''
\emph{Journal of Machine Learning Research} 6 (Dec): 2075--2129.

\hypertarget{ref-hall:1992}{}
Hall, Peter. 1992. \emph{The Bootstrap and Edgeworth Expansion}.
Springer, New York.

\hypertarget{ref-hall:mela:1988}{}
Hallin, Marc, and Guy Mélard. 1988. ``Rank Based Tests for Randomness
Against First Order Serial Dependence.'' \emph{Journal of the American
Statistical Association} 83 (404): 1117--29.

\hypertarget{ref-hall:inge:puri:1985}{}
Hallin, Marc, Jean-François Ingenbleek, and Madan L Puri. 1985. ``Linear
Serial Rank Tests for Randomness Against ARMA Alternatives.'' \emph{The
Annals of Statistics} 13 (3): 1156--81.

\hypertarget{ref-hast:tibs:1990}{}
Hastie, Trevor, and Rob Tibshirani. 1990. \emph{Generalized Additive
Models}. Chapman; Hall, London.

\hypertarget{ref-hell:hell:gorf:2013}{}
Heller, Ruth, Yair Heller, and Malka Gorfine. 2013. ``A Consistent
Multivariate Test of Association Based on Ranks of Distances.''
\emph{Biometrika} 100 (2): 503--10.

\hypertarget{ref-hell:hell:2016}{}
Heller, Ruth, Yair Heller, Shachar Kaufman, Barak Brill, and Malka
Gorfine. 2016. ``Consistent Distribution-Free K-Sample and Independence
Tests for Univariate Random Variables.'' \emph{Journal of Machine
Learning Research} 17 (1): 978--1031.

\hypertarget{ref-hini:1982}{}
Hinich, Melvin J. 1982. ``Testing for Gaussianity and Linearity of a
Stationary Time Series.'' \emph{Journal of Time Series Analysis} 3 (3):
169--76.

\hypertarget{ref-hjor:jone:1996}{}
Hjort, Nils Lid, and M. Chris Jones. 1996. ``Locally Parametric
Nonparametric Density Estimation.'' \emph{Annals of Statistics} 24 (4):
1619--47.

\hypertarget{ref-hoef:1940}{}
Hoeffding, Wassily. 1940. ``Mass-Stabinvariante Korrelationstheorie.''
\emph{Schriften Des Mathematischen Seminars Und Des Instituts Für
Angewandte Mathematik Der Universität Berlin} 5 (3): 179--233.

\hypertarget{ref-hoef:1948}{}
---------. 1948. ``A Nonparametric Test of Independence.'' \emph{Annals
of Mathematical Statistics} 19 (4): 546--57.

\hypertarget{ref-hoef:1963}{}
---------. 1963. ``Probability Inequalities for Sums of Bounded Random
Variables.'' \emph{Journal of the American Statistical Association} 58
(301): 13--30.

\hypertarget{ref-holl:wang:1987}{}
Holland, Paul W, and Yuchung J Wang. 1987. ``Dependence Functions for
Continuous Bivariate Densities.'' \emph{Communications in Statistics} 16
(3): 863--76.

\hypertarget{ref-hong:1998}{}
Hong, Yongmiao. 1998. ``Testing for Pairwise Serial Independence via the
Empirirical Distribution Function.'' \emph{Journal of the Royal
Statistical Society Series B} 60 (2): 429--60.

\hypertarget{ref-hong:1999}{}
---------. 1999. ``Hypothesis Testing in Time Series via the Empirical
Characteristic Function: A Generalized Spectral Density Approach.''
\emph{Journal of the American Statistical Association} 94 (448):
1201--20.

\hypertarget{ref-hong:2000}{}
---------. 2000. ``Generalized Spectral Tests for Serial Dependence.''
\emph{Journal of the Royal Statistical Society Series B} 62 (3):
557--74.

\hypertarget{ref-hong:lee:2003a}{}
Hong, Yongmiao, and Tae-Hwy Lee. 2003. ``Inference on Predictability of
Foreign Exchange Rates via Generalized Spectrum and Nonlinear Time
Series Models.'' \emph{Review of Economics and Statistics} 80 (4):
188--201.

\hypertarget{ref-hong:whit:2005}{}
Hong, Yongmiao, and Halbert White. 2005. ``Asymptotic Distribution
Theory for Nonparametric Entropy Measures of Serial Dependence.''
\emph{Econometrica} 73 (3): 837--901.

\hypertarget{ref-huan:2010}{}
Huang, Tzee-Ming. 2010. ``Testing Conditional Independence Using Maximal
Nonlinear Conditional Correlation.'' \emph{Annals of Statistics} 38 (4):
2047--91.

\hypertarget{ref-ibra:2009}{}
Ibragimov, Rustam. 2009. ``Copula Based Characterizations for
Higher-Order Markov Processes.'' \emph{Econometric Theory} 25 (3):
819--46.

\hypertarget{ref-inci:li:mcca:2011}{}
Inci, A Can, Hsi-Cheng Li, and Joseph McCarthy. 2011. ``Financial
Contagion: A Local Correlation Analysis.'' \emph{Research in
International Business and Finance} 25 (1): 11--25.

\hypertarget{ref-jent:leuc:meye:beer:2018}{}
Jentsch, Carsten, Anne Leucht, Marco Meyer, and Carina Beering. 2018.
``Empirical Characteristic Functions-Based Estimation and Distance
Correlation for Locally Stationary Processes.'' Working paper,
University of Mannheim.

\hypertarget{ref-joe:2014}{}
Joe, Harry. 2014. \emph{Dependence Modeling with Copulas}. Chapman;
Hall, London.

\hypertarget{ref-jone:1996}{}
Jones, M. Chris. 1996. ``The Local Dependence Function.''
\emph{Biometrika} 83 (4): 899--904.

\hypertarget{ref-jone:1998}{}
---------. 1998. ``Constant Local Dependence.'' \emph{Journal of
Multivariate Analysis} 64 (2): 148--55.

\hypertarget{ref-jone:koch:2003}{}
Jones, M. Chris, and Inge Koch. 2003. ``Dependence Maps: Local
Dependence in Practice.'' \emph{Statistics and Computing} 13 (3):
241--55.

\hypertarget{ref-jord:tjos:2018b}{}
Jordanger, Lars Arne, and Dag Tjøstheim. 2017a. ``Nonlinear
Cross-Spectrum Analysis via the Local Gaussian Correlation.''
\emph{arXiv Preprint arXiv:1708.02495}.

\hypertarget{ref-jord:tjos:2018a}{}
---------. 2017b. ``Nonlinear Spectral Analysis via the Local Gaussian
Correlation.'' \emph{arXiv Preprint arXiv:1708.02166v2}.

\hypertarget{ref-kend:1938}{}
Kendall, Maurice G. 1938. ``A New Measure of Rank Correlation.''
\emph{Biometrika} 30 (1/2): 81--89.

\hypertarget{ref-kend:1970}{}
---------. 1970. \emph{Rank Correlation Methods}. 4th ed. Griffin,
London.

\hypertarget{ref-king:1987}{}
King, Maxwell L. 1987. ``Testing for Autocorrelation in Linear
Regression Models: A Survey.'' In \emph{Specification Analysis in the
Linear Regression Model}, edited by M.L. King and D.E.A. Giles, 19--73.
Rutledge; Kegan Paul, London.

\hypertarget{ref-kinn:atwa:2014}{}
Kinney, Justin B, and Gurinder S Atwal. 2014. ``Equitability, Mutual
Information, and the Maximal Information Coefficient.''
\emph{Proceedings National Academy of Science USA} 111: 3354--9.

\hypertarget{ref-klaa:well:1997}{}
Klaassen, Chris AJ, and Jon A Wellner. 1997. ``Efficient Estimation in
the Bivariate Normal Copula Model: Normal Margins Are Least Favorable.''
\emph{Bernoulli} 3 (1): 55--77.

\hypertarget{ref-knok:1977}{}
Knoke, James D. 1977. ``Testing for Randomness Against Autocorrelation:
Alternative Tests.'' \emph{Biometrika} 64 (3): 523--29.

\hypertarget{ref-koja:2009}{}
Kojadinovic, Ivan, and Mark Holmes. 2009. ``Tests of Independence Among
Continuous Random Vectors Based on Cramér--von Mises Functionals of the
Empirical Copula Process.'' \emph{Journal of Multivariate Analysis} 100
(6): 1137--54.

\hypertarget{ref-koya:1987}{}
Koyak, R. 1987. ``On Measuring Internal Dependence in a Set of Random
Variables.'' \emph{Annals of Statistics} 15 (3): 1215--28.

\hypertarget{ref-laca:tjos:2018}{}
Lacal, Virginia, and Dag Tjøstheim. 2017a. ``Estimating and Testing
Nonlinear Local Dependence Between Two Time Series.'' \emph{Journal of
Business and Economic Statistics, to Appear}.

\hypertarget{ref-laca:tjos:2017}{}
---------. 2017b. ``Local Gaussian Autocorrelation and Tests of Serial
Dependence.'' \emph{Journal of Time Series Analysis} 38 (1): 51--71.

\hypertarget{ref-land:vald:2003}{}
Landsman, Zinoviy M, and Emiliano A Valdez. 2003. ``Tail Conditional
Expectations for Elliptical Distributions.'' \emph{North American
Actuarial Journal} 7 (4). Taylor \& Francis: 55--71.

\hypertarget{ref-lehm:1966}{}
Lehmann, Erich Leo. 1966. ``Some Concepts of Dependence.'' \emph{Annals
of Mathematical Statistics} 37 (5): 1137--53.

\hypertarget{ref-load:1996}{}
Loader, Clive R. 1996. ``Local Likelihood Density Estimation.''
\emph{Annals of Statistics} 24 (4): 1602--18.

\hypertarget{ref-low:alco:brai:faff:2013}{}
Low, Rand Kwong Yew, Jamie Alcock, Robert Faff, and Timothy Brailsford.
2013. ``Canonical Vine Copulas in the Context of Modern Portfolio
Management: Are They Worth It?'' \emph{Journal of Banking and Finance}
37 (8): 3085--99.

\hypertarget{ref-lyon:2013}{}
Lyons, Russell. 2013. ``Distance Covariance in Metric Spaces.''
\emph{Annals of Probability} 41 (5): 3284--3305.

\hypertarget{ref-mang:2017}{}
Mangold, Benedikt. 2017. ``A Multivariate Rank Test of Independence
Based on a Multiparametric Polynomial Copula.'' IWQW Discussion Papers
10/2015, University of Erlangen-Nürnberg.

\hypertarget{ref-mark:1952}{}
Markowitz, Harry. 1952. ``Portfolio Selection.'' \emph{Journal of
Finance} 7 (1): 77--91.

\hypertarget{ref-mcle:wkli:1983}{}
McLeod, Allan I, and William K Li. 1983. ``Diagnostic Checking ARMA Time
Series Models Using Squared Residuals and Autocorrelations.''
\emph{Journal of Time Series Analysis} 4 (4): 269--73.

\hypertarget{ref-meuc:2011}{}
Meucci, A. 2011. ``A New Breed of Copulas for Risk and Portfolio
Management.'' \emph{Risk} 24 (9): 122--26.

\hypertarget{ref-murr:murr:murr:2014}{}
Murrell, Ben, Daniel Murrell, and Hugh Murrell. 2014.
``\(R^2\)-Equitability Is Satisfiable.'' \emph{Proceedings National
Academy of Science USA} 111 (21): E2160.

\hypertarget{ref-musc:2014}{}
Muscat, Joseph. 2014. \emph{Functional Analysis: An Introduction to
Matric Spaces, Hilbert Spaces and Banach Algebras}. Springer, New York.

\hypertarget{ref-nels:1999}{}
Nelsen, Roger B. 1999. \emph{An Introduction to Copulas}. Springer, New
York.

\hypertarget{ref-vonn:1941}{}
Neumann, John von. 1941. ``Distribution of the Ratio of Mean Square
Successive Differences to the Variance.'' \emph{Annals of Mathematical
Statistics} 12 (4): 367--95.

\hypertarget{ref-vonn:1942}{}
---------. 1942. ``A Further Remark Concerning the Distribution of the
Ratio of Mean Square Difference to the Variance.'' \emph{Annals of
Mathematical Statistics} 13 (1): 86--88.

\hypertarget{ref-oh:patt:2018}{}
Oh, Dong Hwan, and Andrew J. Patton. 2018. ``Time-Varying Systemic Risk
from a Dynamic Copula Model of Cds Spreads.'' \emph{Journal of Business
and Economic Statistics} 36 (2): 181--95.

\hypertarget{ref-otne:2016}{}
Otneim, Håkon. 2016. ``Multivariate and Conditional Density Estimation
Using Local Gaussian Approximations.'' PhD thesis, Doctoral thesis,
University of Bergen.

\hypertarget{ref-otne:tjos:2017}{}
Otneim, Håkon, and Dag Tjøstheim. 2017. ``The Locally Gaussian Density
Estimator for Multivariate Data.'' \emph{Statistics and Computing} 27
(6). Springer: 1595--1616.

\hypertarget{ref-otne:tjos:2018}{}
---------. 2018. ``Conditional Density Estimation Using the Local
Gaussian Correlation.'' \emph{Statistics and Computing} 28 (2).
Springer: 303--21.

\hypertarget{ref-patt:2012}{}
Patton, Andrew J. 2012. ``A Review of Copula Models for Economic Time
Series.'' \emph{Journal of Multivariate Analysis} 110: 4--18.

\hypertarget{ref-pear:1896}{}
Pearson, Karl. 1896. ``Mathematical Contributions to the Theory of
Evolution. Iii. Regression, Heredity and Panmixia.'' \emph{Philosophical
Transactions of the Royal Society of London} 187: 253--318.

\hypertarget{ref-pear:1922}{}
---------. 1922. \emph{Francis Galton: A Centenary Appreciation}.
Cambridge University Press, Cambridge.

\hypertarget{ref-pear:1930}{}
---------. 1930. \emph{The Life, Letters and Labors of Francis Galton}.
Cambridge University Press, Cambridge.

\hypertarget{ref-pfis:2017}{}
Pfister, Niklas, and Jonas Peters. 2017. \emph{dHSIC: Independence
Testing via Hilbert Schmidt Independence Criterion}.
\url{https://CRAN.R-project.org/package=dHSIC}.

\hypertarget{ref-pfis:2018}{}
Pfister, Niklas, Peter Bühlmann, Bernhard Schölkopf, and Jonas Peters.
2018. ``Kernel-Based Tests for Joint Independence.'' \emph{Journal of
the Royal Statistical Society Series B} 80 (1): 5--31.

\hypertarget{ref-pink:1998}{}
Pinkse, Joris. 1998. ``Consistent Nonparametric Testing for Serial
Independence.'' \emph{Journal of Econometrics} 84 (2): 205--31.

\hypertarget{ref-prud:bryc:mari:1986}{}
Prudnikov, A.P., A. Brychkov, and O.I Marichev. 1986. \emph{Integrals
and Series}. Gordon Breach Science Publisher, New York.

\hypertarget{ref-resh:resh:etal:2011}{}
Reshef, David N, Yakir A Reshef, Hilary K Finucane, Sharon R Grossman,
Gilean McVean, Peter J Turnbaugh, Eric S Lander, Michael Mitzenmacher,
and Pardis C Sabeti. 2011. ``Detecting Novel Associations in Large
Datasets.'' \emph{Science} 334 (6062): 1518--24.

\hypertarget{ref-resh:resh:etal:2014}{}
Reshef, David N, Yakir A Reshef, Michael Mitzenmacher, and Pardis C
Sabeti. 2014. ``Cleaning up the Record on the Maximal Information
Coefficient and Equitability.'' \emph{Proceedings National Academy of
Science USA} 111 (33): E3362--E3363.

\hypertarget{ref-resh:resh:etal:2013}{}
Reshef, David, Yakir Reshef, Michael Mitzenmacher, and Pardis Sabeti.
2013. ``Equitability Analysis of the Maximal Information Coefficient,
with Comparisons.''

\hypertarget{ref-reny:1959}{}
Rényi, Alfréd. 1959. ``On Measures of Dependence.'' \emph{Acta
Mathematica Hungarica} 10 (3-4): 441--51.

\hypertarget{ref-reny:1961}{}
Rényi, Alfréd. 1961. ``On Measures of Information and Entropy.'' In
\emph{Proceedings of the 4th Berkeley Symposium on Mathematical
Statistics and Probability 1960}, 547--61.

\hypertarget{ref-rizz:szek:2018}{}
Rizzo, Maria L., and Gabor J. Szekely. 2018. \emph{Energy: E-Statistics:
Multivariate Inference via the Energy of Data}.
\url{https://CRAN.R-project.org/package=energy}.

\hypertarget{ref-robi:1991}{}
Robinson, Peter M. 1991. ``Consistent Nonparametric Entropy-Based
Testing.'' \emph{Review of Economic Studies} 58 (3): 437--53.

\hypertarget{ref-rose:1975}{}
Rosenblatt, Murray. 1975. ``A Quadratic Measure of Deviation of
Two-Dimensional Density Estimates and a Test of Independence.''
\emph{Annals of Statistics} 3 (1): 1--14.

\hypertarget{ref-sank:gupt:2004}{}
Sankaran, P.G., and R.P. Gupta. 2004. ``Characterizations Using Local
Dependence Function.'' \emph{Communications in Statistics. Theory and
Methods} 33 (12): 2959--74.

\hypertarget{ref-sejd:srip:gret:fuku:2013}{}
Sejdinovic, Dino, Bharath Sriperumbudur, Arthur Gretton, and Kenji
Fukumizu. 2013. ``Equivalence of Distance-Based and RKHS-Based
Statistics in Hypothesis Testing.'' \emph{Annals of Statistics} 41 (5):
2263--91.

\hypertarget{ref-silv:gran:2001}{}
Silvapulle, Param, and Clive W.J. Granger. 2001. ``Large Returns,
Conditional Correlation and Portfolio Diversification. a Value-at-Risk
Approach.'' \emph{Quantative Finance} 1 (5): 542--51.

\hypertarget{ref-simo:tibs:2014}{}
Simon, Noah, and Robert Tibshirani. 2014. ``Comment on `Detecting Novel
Associations in Large Data Sets' by Reshev et Al Science Dec 16, 2011.''
\emph{arXiv Preprint: 1401.7645v1}.

\hypertarget{ref-skau:1993}{}
Skaug, Hans Julius. 1993. ``The Limit Distribution of the Hoeffding
Statistic for Tests of Serial Independence.'' Manuscript, Department of
Mathematics, University of Bergen.

\hypertarget{ref-skau:tjos:1993b}{}
Skaug, Hans Julius, and Dag Tjøstheim. 1993a. ``A Nonparametric Test of
Serial Independence Based on the Empirical Distribution Function.''
\emph{Biometrika} 80 (3): 591--602.

\hypertarget{ref-skau:tjos:1993a}{}
---------. 1993b. ``Nonparametric Tests for Serial Independence.'' In
\emph{Developments in Time Series Analysis, the Priestley Birthday
Volume}, edited by T. Subba Rao, 207--30. Chapman; Hall, London.

\hypertarget{ref-skau:tjos:1996}{}
---------. 1996. ``Testing for Serial Independence Using Measures of
Distance Between Densities.'' In \emph{Athens Conference on Applied
Probability and Time Series Vol. Ii, in Memory of E.J. Hannan}, edited
by P. M. Robinson and M. Rosenblatt, 115:363--78. Springer Lecture Notes
in Statistics. Springer, Berlin.

\hypertarget{ref-skla:1959}{}
Sklar, A. 1959. \emph{Fonctions de Répartition à N Dimensions et Leurs
Marges}. Université Paris 8.

\hypertarget{ref-smit:2015}{}
Smith, Michael Stanley. 2015. ``Copula Modelling of Dependence in
Multivariate Time Series.'' \emph{International Journal of Forecasting}
31 (3): 815--33.

\hypertarget{ref-smit:min:alme:czad:2010}{}
Smith, Michael, Aleksey Min, Carlos Almeida, and Claudia Czado. 2010.
``Modeling Longitudinal Data Using a Pair-Copula Decomposition of Serial
Dependence.'' \emph{Journal of the American Statistical Association} 61
(492): 1467--79.

\hypertarget{ref-spea:1904}{}
Spearman, Charles. 1904. ``The Proof and Measurement of Association
Between Two Things.'' \emph{American Journal of Psychology} 15 (1):
72--101.

\hypertarget{ref-stan:2001}{}
Stanton, Jeffrey M. 2001. ``Galton, Pearson, and the Peas: A Brief
History of Linear Regression for Statistics Instructors.'' \emph{Journal
of Statistical Education} 9 (3): 1--13.

\hypertarget{ref-stei:1999}{}
Stein, Michael L. 1999. \emph{Interpolation of Spatial Data: Some Theory
for Kriging}. Springer, New York.

\hypertarget{ref-stig:1989}{}
Stigler, Stephen M. 1989. ``Francis Galton's Account of the Invention of
Correlation.'' \emph{Statistical Science} 4 (2): 73--86.

\hypertarget{ref-stov:tjos:2014}{}
Støve, Bård, and Dag Tjøstheim. 2014. ``Measuring Asymmetries in
Financial Returns: An Empirical Investigation Using Local Gaussian
Correlation.'' In \emph{Essays in Nonlinear Time Series Econometrics},
edited by M. Meitz N. Haldrup and P. Saikkonen, 307--29. Oxford
University Press, Oxford.

\hypertarget{ref-stov:tjos:huft:2014}{}
Støve, Bård, Dag Tjøstheim, and K. Hufthammer. 2014. ``Using Local
Gaussian Correlation in a Nonlinear Re-Examination of Financial
Contagion.'' \emph{Journal of Empirical Finance} 25: 785--801.

\hypertarget{ref-su:whit:2007}{}
Su, Liangjun, and Halbert White. 2007. ``A Consistent
Characteristic-Function-Based Test for Conditional Independence.''
\emph{Journal of Econometrics} 141 (2): 807--37.

\hypertarget{ref-subb:gabr:1980}{}
Subba Rao, T., and M. M. Gabr. 1980. ``A Test for Linearity of
Stationary Time Series.'' \emph{Journal of Time Series Analysis} 1 (2):
145--58.

\hypertarget{ref-szek:rizz:2005}{}
Szekely, Gabor J, and Maria L Rizzo. 2005. ``Hierarchical Clustering via
Joint Between-Within Distances: Extending Ward's Minimum Variance
Method.'' \emph{Journal of Classification} 22 (2): 151--83.

\hypertarget{ref-szek:rizz:2014}{}
---------. 2014. ``Partial Distance Correlation with Methods for
Dissimilarities.'' \emph{Annals of Statistics} 42 (6): 2382--2412.

\hypertarget{ref-szek:2002}{}
Székely, Gábor J. 2002. ``\({\mathcal E}\)-Statistics: The Energy of
Statistical Samples.'' Technical report 02-16, Bowling Green State
University.

\hypertarget{ref-szek:rizz:2009}{}
Székely, Gábor J, and Maria L Rizzo. 2009. ``Brownian Distance
Covariance.'' \emph{Annals of Applied Statistics} 3 (4): 1236--65.

\hypertarget{ref-szek:rizz:2012}{}
---------. 2012. ``On the Uniqueness of Distance Correlation.''
\emph{Statistics and Probability Letters} 82 (12): 2278--82.

\hypertarget{ref-szek:rizz:2013}{}
---------. 2013. ``Energy Statistics: A Class of Statitics Based on
Distances.'' \emph{Journal of Statistical Planning and Inference} 143
(8): 1249--72.

\hypertarget{ref-szek:rizz:baki:2007}{}
Székely, Gábor J, Maria L Rizzo, and Nail K Bakirov. 2007. ``Measuring
and Testing Dependence by Correlation of Distances.'' \emph{Annals of
Statistics} 35 (6): 2769--94.

\hypertarget{ref-tale:2007}{}
Taleb, Nassim Nicholas. 2007. \emph{The Black Swan: The Impact of the
Highly Improbable}. Random house.

\hypertarget{ref-tera:tjos:gran:2010}{}
Teräsvirta, Timo, Dag Tjøstheim, and Clive W.J. Granger. 2010.
\emph{Modelling Nonlinear Economic Time Series}. Oxford University
Press.

\hypertarget{ref-tjos:otne:stov:2020}{}
Tjøstheim, D., H. Otneim, and B. Støve. 2020. \emph{Statistical Modeling
Using Local Gaussian Approximation}. Elsevier, Amsterdam, to appear.

\hypertarget{ref-tjos:1996}{}
Tjøstheim, Dag. 1996. ``Measures of Dependence and Tests of
Independence.'' \emph{Statistics} 28 (3): 249--84.

\hypertarget{ref-tjos:huft:2013}{}
Tjøstheim, Dag, and Karl O. Hufthammer. 2013. ``Local Gaussian
Correlation: A New Measure of Dependence.'' \emph{Journal of
Econometrics} 172 (1): 33--48.

\hypertarget{ref-tuke:1958}{}
Tukey, John W. 1958. ``A Problem of Berkson, and Minimum Variance
Orderly Estimators.'' \emph{Annals of Mathematical Statistics} 29 (2):
588--92.

\hypertarget{ref-vand:1952}{}
Waerden, Bartel Leendert van der. 1952. ``Order Tests for the Two-Sample
Problem and Their Power.'' \emph{Idagationes Mathematicae} 55: 453--58.

\hypertarget{ref-wang:li:etal:2015}{}
Wang, Yi, Yi Li, Hongbao Cao, Momiao Xiong, Yin Yao Shugart, and Li Jin.
2015. ``Efficient Test for Nonlinear Dependence of Two Continuous
Variables.'' \emph{BMC Bioinformatics} 16 (1): 260.

\hypertarget{ref-wang:li:etal:2017}{}
Wang, Yi, Yi Li, Xiaoyu Liu, Weilin Pu, Xiaofeng Wang, Jiucun Wang,
Momiao Xiong, Yin Yao Shugart, and Li Jin. 2017. ``Bagging
Nearest-Neighbor Prediction Independence Test: An Efficient Method for
Nonlinear Dependence of Two Continuous Variables.'' \emph{Scientific
Reports} 7 (1).

\hypertarget{ref-wilc:2005}{}
Wilcox, Rand R. 2005. ``Estimating the Conditional Variance of \(Y\),
Given \(X\), in a Simple Regression Model.'' \emph{Journal of Applied
Statistics} 32 (5): 495--502.

\hypertarget{ref-wilc:2007}{}
---------. 2007. ``Local Measures of Association: Estimating the
Derivative of the Regression Line.'' \emph{British Journal of
Mathematical and Statistical Psychology} 60 (1): 107--17.

\hypertarget{ref-yao:zhan:shao:2018}{}
Yao, Shun, Xianyang Zhang, and Xiaofeng Shao. 2018. ``Testing Mutual
Independence in High Dimension via Distance Covariance.'' \emph{Journal
of the Royal Statistical Society Series B} 80 (3): 455--80.

\hypertarget{ref-yeni:rizz:2014}{}
Yenigün, C Deniz, and Maria L. Rizzo. 2014. ``Variable Selection in
Regression Using Maximal Correlation and Distance Correlation.''
\emph{Journal of Statistical Computation and Simulation} 85 (8):
1692--1705.

\hypertarget{ref-yeni:szek:rizz:2011}{}
Yenigün, C Deniz, Gabor J. Székely, and Maria L. Rizzo. 2011. ``A Test
of Independence in Two-Way Contingency Tables Based on Maximal
Correlation.'' \emph{Communications in Statistics - Theory and Methods}
40 (12): 2225--42.

\hypertarget{ref-zhan:pete:janz:scho:2012}{}
Zhang, Kun, Jonas Peters, Dominik Janzing, and Bernhard Schölkopf. 2012.
``Kernel-Based Conditional Independence Test and Applications in Causal
Discovery.'' In \emph{Proceedings of the Uncertainty in Artificial
Intelligence}, 804--13. AUAI Press, Corvallis Oregon.

\hypertarget{ref-zhou:2012}{}
Zhou, Zhou. 2012. ``Measuring Nonlinear Dependence in Time Series, a
Distance Correlation Approach.'' \emph{Journal of Time Series Analysis}
33 (3): 438--57.

\end{document}